\documentclass[11pt]{article}%
\usepackage{makeidx}
\usepackage{amssymb}
\usepackage[all]{xy}
\usepackage[latin1]{inputenc}
\usepackage{amsfonts}
\usepackage{amsmath}
\usepackage{amssymb}
\usepackage{url}
\usepackage{hyperref}
\usepackage{graphicx}%
\usepackage{color}
\usepackage{xcolor}
\usepackage{amsbsy}
\usepackage{yfonts}
\usepackage[all]{xy}
\usepackage{makeidx}
\usepackage{color}
\usepackage{mathrsfs}

\usepackage{pgfplots}

\usetikzlibrary{matrix}
\usepackage{fancybox}

\usepackage{color}

\providecommand{\U}[1]{\protect\rule{.1in}{.1in}}

\numberwithin{equation}{section}
\providecommand{\U}[1]{\protect\rule{.1in}{.1in}}
\providecommand{\U}[1]{\protect\rule{.1in}{.1in}}
\newtheorem{theo}{Theorem}[section]
\newtheorem{prop}[theo]{Proposition}
\newtheorem{lem}[theo]{Lemma}
\newtheorem{cor}[theo]{Corollary}

\newtheorem{rem}[theo]{Remark}
\newtheorem{defi}[theo]{Definition}

\newcommand{\CC}{\mathbb{C}}

\newcommand{\EE}{\mathbb{E}}

\newcommand{\KK}{\mathbb{K}}
\newcommand{\LL}{\mathbb{L}}

\newcommand{\NN}{\mathbb{N}}
\newcommand{\PP}{\mathbb{P}}

\newcommand{\RR}{\mathbb{R}}

\newcommand{\Aa}{ {\cal A }}
\newcommand{\Ba}{ {\cal B }}
\newcommand{\Ca}{ {\cal C }}

\newcommand{\Ka}{ {\cal K }}

\newcommand{\Ya}{ {\cal Y }}

\newcommand{\Va}{ {\cal V }}

\newcommand{\Fa}{ {\cal F }}

\newcommand{\Fo}{ \mbox{\it Forest}}

\newcommand{\Ga}{ {\cal G }}
\newcommand{\Qa}{ {\cal Q }}
\newcommand{\Ia}{ {\cal I }}
\newcommand{\Xa}{ {\cal X }}
\newcommand{\Ma}{ {\cal M }}
\newcommand{\Ta}{ {\cal T}}
\newcommand{\Sa}{ {\cal S}}

\newcommand{\Ja}{ {\cal J }}
\newcommand{\Ju}{ \mbox{\it Jungle}}
\newcommand{\Pa}{ {\cal P }}
\newcommand{\Za}{ {\cal Z }}
\newcommand{\Wa}{ {\cal W }}

\newcommand{\point}{\mbox{\LARGE .}}
\newcommand{\proof}{\noindent\mbox{\bf Proof:}\\}
\newcommand{\cqfd}{\hfill\blbx \\}
\def\blbx{\hbox{\vrule height 5pt width 5pt depth 0pt}\medskip}

\def \PP{\mathbb{P}}
\def \RR{\mathbb{R}}
\def \EE{\mathbb{E}}

\def \CC{\mathbb{C}}
\def \LL{\mathbb{L}}

\def \AA{\mathbb{A}}
\def \BB{\mathbb{B}}

\def \XX{\mathbb{X}}

\def \YY{\mathbb{Y}}

\usepackage{bbm}

\def \KK{\mathbbm{K}}

\let\oldchi\chi
\renewcommand{\chi}[1][1.3pt]{%
  \mathrel{\raisebox{#1}{\scalebox{1.2}{$\oldchi$}}}%
}

\def \XXb{\boldsymbol{\mathbbm{X}}}

\def \XX{\mathbbm{X}}

\def \MM{\boldsymbol{\mathbbm{M}}}

\def \MMb{\mathbbm{M}}

\begin{document}

\title{On particle Gibbs Markov chain Monte Carlo models}
\author{P. Del Moral\thanks{School of Mathematics and Statistics, University of New South Wales, p.del-moral@unsw.edu.au}, 
R. Kohn\thanks{School of Economics, University of New South Wales, r.kohn@unsw.edu.au}, F. Patras\thanks{Universit\'e de Nice et CNRS, patras@unice.fr}
}
\maketitle

\begin{abstract}
   This article analyses a new class of 
   advanced particle Markov chain Monte Carlo algorithms recently 
   introduced by Andrieu, Doucet,  and Holenstein (2010).
   We present a natural interpretation of these methods in terms of well known unbiasedness properties
of Feynman-Kac particle measures, and a new duality with Feynman-Kac models.
 
This perspective sheds a new light on the foundations and the mathematical analysis of this class of methods. A key consequence is the equivalence between
the backward and ancestral particle Markov chain Monte Carlo methods, with the Gibbs sampling of a (many-body) Feynman-Kac target distribution. 
 Our approach also presents a new stochastic differential calculus
based on geometric combinatorial techniques
to derive explicit non-asymptotic Taylor type series of the semigroup of a class of particle Markov chain Monte Carlo models around their invariant measures with respect to the population
size of the auxiliary particle sampler.  
These results provide sharp quantitative estimates of the convergence properties of conditional particle Markov chain models with respect to the time horizon and the size of the systems. We illustrate the implication of these results with sharp estimates of the contraction coefficient and the Lyapunov exponent of conditional particle samplers, and
explicit and non-asymptotic $\LL_p$-mean error decompositions of the law of the random states around the limiting invariant measure. The abstract framework developed
in the article also allows the design of natural extensions to island (also called SMC$^2$) type particle methodologies.

\end{abstract}
\section{Introduction}
In the last two decades, particle simulation techniques have become one of the most active
contact points between Bayesian statistical inference and applied probability. Their range of applications
goes from statistical machine learning to information theory, theoretical chemistry
and quantum physics, financial mathematics, signal processing, risk analysis,
and several other domains in engineering and computer sciences. 
 In contrast to conventional Markov chain Monte Carlo
methodologies, particle methods  are not based on sampling long runs
of a judiciously chosen Markov chain with a prescribed target probability measure, 
but on the mean field particle simulation of nonlinear Markov chain models. 

The seminal article~\cite{adh-2010} by Andrieu, Doucet and Holenstein introduced a new way to combine Markov chain Monte Carlo methods ({\em
MCMC}) with Sequential Monte Carlo methodologies ({\em 
SMC}). A variant of this particle Gibbs type method, where ancestors are resampled in a forward pass, is developed in
 Sch\"on and Jordan in ~\cite{lindsten}, and Lindsten and Sch\"on~\cite{lindsten2}.  

This new class of Monte Carlo samplers is termed particle Markov chain Monte Carlo methods ({\em 
PMCMC}).   These emerging particle sampling technologies are particularly important in signal processing and in Bayesian statistics.
In these application areas, they are used to estimate posterior distributions of unknown parameters
when the likelihood functions are unknown or computationally untractable. 
The central idea here is to run a MCMC sampler and compute the likelihood function using 
an auxiliary particle sampler. In this situation, the updates of the resulting particle MCMC samplers are defined on extended state spaces so that the marginal of their invariant measure coincides with the desired posterior distribution.

Recently, these powerful PMCMC methodologies have attracted considerable attention in a variety of application domains, including 
in statistical machine leaning~\cite{chopin,kantas,lindsten,whiteley}, finance and econometrics~\cite{creal,dejong,flury,lopes,pitt}, biology~\cite{golin,mish,rasmussen}, computer science~\cite{soren}, environmental statistics~\cite{down,down2,gareth}, social networks analysis~\cite{everitt},
signal processing~\cite{lindsten2,olsson},
 forecasting and data assimilation~\cite{launay,morad,vrugt}. 
 
The convergence analysis of the PMCMC methods was started in a series of 
articles~\cite{alv-2013,sumeet,ldm-2014,lindsten,lindsten2}.
The $\phi$-irreducibility and aperiodicity of PMCMC models was already discussed 
in the pioneering article~\cite{adh-2010}.
The uniform ergodicity results with quantitative estimates of the convergence properties of PMCMC models were presented by Chopin and Singh in~\cite{sumeet},
using a coupling technology of ancestral particle paths.
More refined contraction estimates have been obtained by Andrieu, Lee, and Vihola~\cite{alv-2013}
using a conditional type analysis of the normalizing particle constants, and in the article
by Lindsten, Douc, Moulines~\cite{ldm-2014} which provides similar quantitative estimates
using lower bound estimates of PMCMC transition probabilities based on the stability of Feynman-Kac semigroups.

In all of these studies, the validity of PMCMC samplers is assessed by interpreting these models as a traditional MCMC sampler
on a sophisticated and extended state space in which all the random variables generated by some particle model are seen as auxiliary variables. The target measure of these MCMC models is expressed in terms of a density involving 
compositions of random mappings encoding the full ancestral lineages of all the
genetic type particle, from the origin up to the final time horizon.

This article is concerned with an alternative probabilistic foundation of PMCMC methodology. It is well known that
Metropolis-Hasting type PMCMC methods reduce to
standard Metropolis-Hasting with a target 
Feynman-Kac distribution that encapsulates the distribution of the auxiliary particle model
(see, for instance, the article~\cite{adh-2009}, and Section 2.2.5 in~\cite{dhw-12}, in the context of Hidden Markov Chain problems with fixed unknown parameters).

In the first part of the article, we show that the conditional SMC type PMCMC method reduces to a standard
Gibbs sampler (see e.g. Section~\ref{sub:frozen} for a description of the latter). 
The proof of this result is based on
a new duality relation between Feynman-Kac measures on path spaces and their many-body version (hereafter, many-body Feynman-Kac models refer to the extension of usual Feynman-Kac models to collective motions of particles, see Section~\ref{mbody-pm}  for precise definitions). This duality relation can be seen as an extension of
the well known unbiasedness properties
of unnormalized particle measures to these many-body Feynman-Kac models. 

This natural viewpoint simplifies considerably the design and the convergence analysis 
of this class of particle models. 

Finally, the 
new formulation also allows the design of new and natural classes of PMCMC samplers based on island type models and particle Gibbs methodology.

The second part of the article is concerned with  the propagation of chaos properties of PMCMC samplers based on the sampling of a particle model
with a frozen trajectory. We  design explicit Taylor type expansions of the law of a finite block of particles in terms of the population size of the auxiliary particle model. These expansions are naturally 
 parametrized by 
 decorated ("infected") forests.
  Their accuracy at any order is related naturally to the number of coalescent edges and
 the number of infections. To the best of our knowledge, these propagation of chaos series are the first results of this type for this class of particle  
 Markov chain Monte Carlo.

These expansions provide Taylor decompositions of 
the semigroup of conditional PMCMC chains {\it around their invariant target measures} 
w.r.t. the precision parameter $1/N$, where $N$ stands for the size of the auxiliary particle system. It is known that the many-body target Feynman-Kac distribution is the equilibrium measure of the particle Gibbs sampler, for any choice of the population size $N$. Nevertheless the rate of convergence to the equilibrium of this particle MCMC method strongly depends on the parameter $N$. Under some stability properties on Feynman-Kac
semigroups, a direct consequence of these Taylor series expansions is that the Dobrushin contraction coefficient is of order $n/N$, where $n$ stands for the length of the trajectories. This shows that the particle Gibbs method converge to the desired target measure as the number of iterations {\em or} the number of particles tends to infinity. 

The linear dependency of the rate w.r.t. the time horizon $n$ is not surprising. The same type of linear scaling arise in the propagation of chaos and the fluctuation analysis of Feynman-Kac particle models on path spaces (see, for instance, Corollary 8.5.1 in~\cite{d-2004}, as well as Corollary 15.2.5 and Theorem 16.5.1 in \cite{d-2013}). 
 Conditional SMC methods can also be interpreted as mean field particle models with a given frozen path which also has a negligible impact of order $1/N$ upon the auxiliary particle system. In this context, we derive the rate to equilibrium $n/N$ for particle Gibbs samplers using propagation of chaos expansions of the distribution of these particle models, as soon as the Feynman-Kac semigroup of the marginal target measures forgets its initial condition (see, for instance Theorem~\ref{theo-key-introduction}, and the regularity condition (\ref{H})).
 
  Under stronger mixing conditions on the reference Markov chain of Feynman-Kac models, such as the minorization conditions discussed in the end of Section~\ref{classical-FK}, similar decays rates to equilibrium with linear scaling have also been derived in~\cite{alv-2013,ldm-2014} using different lower bound estimation techniques on the minorizing constant for the particle Gibbs kernel.  As noticed in~\cite{ldm-2014} (Section 4.2), these strong mixing conditions are stronger than the Doeblin type condition on Feynman-Kac semigroups since they typically require that the state space is compact. Using partially heuristic arguments, the linear scaling of the number of particles w.r.t. the time horizon is also discussed in the articles~\cite{dpk-2012,pitt}. 

The impact of the Taylor series expansions developed in the present article is also illustrated with sharp 
and non-asymptotic expansions of the Dobrushin contraction coefficient
of any iterated conditional PMCMC transitions. 
We also provide an explicit decomposition of the
 $\LL_p$-distance between the law of the random states of a class of PMCMC methods
 around the limiting invariant measure. 
These results can also be used to estimate the 
bias and the variance of the random states of the occupation measures of the auxiliary particle systems. This duality between Feynman-Kac models and their many-body
versions allows us to apply these Taylor expansions to the original Feynman-Kac particle models.

 Finally, the duality relation and differential calculus developed in this article also open an avenue of research problems in the field of Feynman-Kac particle models  and PMCMC methodologies. One important problem is extending the result developed in the article to continuous time Feynman-Kac models.
 In this context, it is also important to analyze the  effect of the time discretization of the models. Another important issue is comparing the
 stability properties of PMCMC methods with Metropolis-Hasting moves with the one based on Gibbs sampling. Our stability analysis is restricted to particle Gibbs methods on the space of ancestral lines. We plan to analyze the convergence to equilibrium
 of particle Metropolis-Hasting models and the one of the Gibbs sampler with a target many-body Feynman-Kac measure in a forthcoming article.

The article is organized as follows. Section~\ref{sec-statements} provides a brief description of Feynman-Kac models and their particle interpretations. We also state some of the main results of the article and provide a series of signal processing and quantum physics illustrations.

Section~\ref{sec-mean-field-intro} reviews some well know results on Feynman-Kac models and their mean field particle interpretation, 
including path space models and backward particle Markov chain measures. Section~\ref{mbody-pm}
introduces many-body Feynman-Kac models aimed at describing the collective motion of particles in the usual Feynman-Kac models. These models are particularly well suited to the analysis of PMCMC samplers.

Section~\ref{conditiona-sec} considers conditional particle MCMC methodology and proves that the conditional SMC samplers reduce to the Gibbs sampling 
of a many-body Feynman-Kac target measure.
Section~\ref{islands-pm} provides a transport equation and a new duality relation between many-body Feynman-Kac models and a conditional Feynman-Kac  particle model with a frozen trajectory. Section~\ref{historical-sec} considers historical particle methods and their dual frozen particle method. For instance, we show that the conditional distribution of {\em the ancestral lines} of the Feynman-Kac particle model w.r.t. its complete ancestral tree
coincides with the backward particle method.  
Section~\ref{g+b-ref},
 presents two equivalent classes of PMCMC methods: genealogical tree based samplers and 
backward sampling methods. 
Section~\ref{taylor-invariant} 
presents a 
basic description of the Taylor expansions  of conditional PMCMC transitions
around their invariant measures.
We also derive important consequences of these expansions, including 
quantitative estimates of the stability properties of these methods, and sharp estimates of the bias and the variance of the random states of 
the PMCMC Markov chain.

Section~\ref{general-second-order} considers the propagation of chaos properties of a conditional PMCMC particle model.
Section~\ref{scp-sec} collects some preliminary combinatorial results on tensor products 
of empirical measures. Section~\ref{un-norm-sec} considers non-asymptotic Taylor series of 
$q$-tensor products of unnormalized
particle measures. 
Section~\ref{norm-FK-ref} discusses 
the propagation of chaos properties and related Taylor expansions 
of frozen particle models. 
Section~\ref{ifexp-ref} describes the Taylor series decompositions in terms of infected and coalescent forest expansions.

Section~\ref{last-section} 
presents a new class of island PMCMC samplers, and
discusses some extensions and open questions. 

\section{Presentation of the models and statement of some results}\label{sec-statements}
\subsection{Classical Feynman--Kac methods}\label{classical-FK}

We first briefly  survey the classical Feynman--Kac method as well as classical exemples of applications --one of the reasons being that all these strategies can potentially be involved in the more complex and refined Feynman--Kac methods with frozen trajectories and PMCMC methods to be 
studied later in this paper. 

We consider a Markov chain $X_n^{\prime}$ evolving in some measurable state space $S^{\prime}_n$ with some Markov transitions kernels $M^{\prime}_n$, and a collection of non-negative bounded and measurable  functions $G_n^{\prime}$
on $S_n^{\prime}$. We let $\eta_n$ be the Feynman-Kac probability measures on the path space $S_n:=\prod_{0\leq p\leq n}S_p^{\prime}$ defined for any bounded measurable function $f_n$ on  $S_n$ by
\begin{equation}\label{FK-intro-statements}
\eta_n(f_n)\propto \EE\left(f_n(X_n)~Z_n(X)\right) \quad\mbox{\rm with}\quad Z_n(X):=\prod_{0\leq p<n}G_p(X_p)
\end{equation}
In (\ref{FK-intro-statements}), $X_n=(X_0^{\prime},\ldots,X_n^{\prime})$ stands for the historical process, and the potential functions are given by $G_p(X_p):=G^{\prime}_p(X^{\prime}_p)$, with $0\leq p\leq n$. 
It is implicitly assumed that the normalizing constants $\EE\left(Z_n(X)\right)$ are positive, so that $\eta_n$ are well defined probability measures.

Notice that the $n$-th time marginal $\eta_n^{\prime}$ of the measure $\eta_n$ on $S^{\prime}$ is given as in (\ref{FK-intro-statements}) by replacing
$X_n$ and $G_p$ by $X^{\prime}_n$ and $G_p^{\prime}$, with $0\leq p\leq n$.
From the pure mathematical viewpoint, the $n$-th time marginal measure have exactly the same form as the one on path space. This structural stability property is pivotal in the design and the mathematical analysis of Feynman-Kac model on path-spaces and their genealogical tree based particle interpretation (see, for instance, the article~\cite{dm-2001} and the path-particle models discussed in (\ref{FK-path-particle})). 
Our article uses this property to define consistently Gibbs type PMCMC transitions on path spaces.
In the illustrations discussed below, it also allows us to describe filtering and smoothing problems, as well as rare event importance distributions and quantum integration models  in a single path-space Feynman-Kac model.

Feynman--Kac models appear in numerous scientific fields including signal processing, statistics, mathematical finance, rare event analysis, chemistry and statistical physics; see
 \cite{cappe-moulines}, \cite{cdho}, \cite{d-2004}, \cite{dm-2000}, \cite{d-2013}, \cite{dg-2005,cdg-2011} and \cite{doucet2001}. Their interpretation 
depends on the application domain. We now briefly give some examples.

\textit{Nonlinear filtering}.  Let $(\Xa_{n},\Ya_{n})_{n\geq0}$ be a Markov chain on some
product state space $(E_{1}\times E_{2})$ whose transition mechanism takes the
form
\begin{equation}\label{ex-filter-gen}
\mathbb{P}\left(  (\Xa_{n},\Ya_{n})\in d(\textsl{x},\textsl{y})~|~(\Xa_{n-1},\Ya_{n-1})\right)
=K_{n}(\Xa_{n-1},d\textsl{x})~g_{n}(\textsl{y},\textsl{x})~\nu_{n}(d\textsl{y}),
\end{equation}
where $\left(  \nu_{n}\right)  _{n\geq0}$ is a sequence of positive measures
on $E_{2}$, $\left( K_{n}\right)  _{n\geq0\text{ }}$ is a sequence of Markov
kernels from $E_{1}$ into itself, and $\left(g_{n}(\point,\textsl{x})\right)_{n\geq0}$ is a
sequence of conditional density functions on $E_{2}$. These filtering models
are often described by a partially observed dynamical random system. For instance,
when $E_1=E_2=\RR$, these systems may take the common form
\begin{equation}\label{ex-filter-1d}
\Xa_n=a(\Xa_{n-1})+\Wa_n\quad\mbox{\rm and}\quad \Ya_n=b(\Xa_n)+\Va_n
\end{equation}
with a sequence of independent random variables $\Wa_n$ and $\Va_n$, and some bounded functions
$a,b$ on $\RR$. When the random variables $\Va_n$ have some density $h_n(v)$ w.r.t. the Lebesgue
measure $dv$ on $\RR$, (\ref{ex-filter-gen}) holds with $g_n(y,x)=h_n(y-b(x))$.
The aim of 
filtering is to infer the trajectories of the hidden Markov process $\Xa_{n}
$ given a series of observations $\Ya_{k}$ from the origin $k=0$, up to the current time $k=n$. Choosing
$$
S^{\prime}_n=E_1\qquad G^{\prime}_{k}(\textsl{x}):=g_{k}(\textsl{y}_{k},\textsl{x}) \quad\mbox{\rm and}\quad M^{\prime}_k=K_k
$$
 it is easily
checked that
\[
\eta_{n}=\mbox{\rm Law}\left(  (\Xa_{0},\ldots,\Xa_n)~|~\Ya_{k}=\textsl{y}_{k},~0\leq k<n\right)  ,
\]
\emph{Importance sampling and rare event analysis.}
Feynman-Kac models are also closely related to importance sampling Monte Carlo methods. For a given target Feynman-Kac measure (\ref{FK-intro-statements}), 
the choice of the Markov chain $X^{\prime}_n$ and
the potential functions $G^{\prime}_n$ is rather flexible. For instance,
 if we set
\begin{equation}\label{ex-filter}
G^{\prime}_{k}(\textsl{x})=\int~K^{\prime}_{k+1}(\textsl{x},d\textsl{z})g_{k+1}(\textsl{y}_{k+1},\textsl{z})\quad\mbox{\rm and}\quad M^{\prime}_k(\textsl{x},d\textsl{z})\propto K^{\prime}_k(\textsl{x},d\textsl{z})g_{k}(\textsl{y}_{k},\textsl{z})
\end{equation}
then 
$
\eta_{n}=\mbox{\rm Law}\left(  (\Xa_{0},\ldots,\Xa_n)~|~\Ya_{k}=y_{k},~0\leq k\leq n\right),
$
as soon as the initial random variable $X^{\prime}_0$ is distributed with the  conditional distribution of $\Xa_0$ given $\Ya_0$. In other situations, such as in rare 
event analysis~\cite{dg-2005} and Markov bridge type sampling problems~\cite{dmurray}, the target distribution has the form
$$
\PP^W_n(d(\textsl{x}_0,\ldots,\textsl{x}_n))\propto W_n(\textsl{x}_n)~\PP_n\left(d(\textsl{x}_0,\ldots,\textsl{x}_n)\right)
$$
where $\PP_n$ stands for the distribution of the random trajectories of a Markov chain $\Xa_k$ evolving in some measurable state space $E$, and $W_n$ is some
importance function that is non-negative. If we choose 
\begin{equation}\label{ex-is}
S^{\prime}_k=E\qquad X^{\prime}_k=\left(\Xa_k,\Xa_{k+1}\right)\quad\mbox{\rm and}\quad G^{\prime}_k\left(\Xa_k,\Xa_{k+1}\right):=W_{k+1}(\Xa_{k+1})/W_{k}(\Xa_{k}).
\end{equation}
then, it can be readily checked that
$$
\eta_n(f_n)\propto\EE\left(F_n\left(\Xa_0,\ldots,\Xa_n\right)~W_n(\Xa_n)\right)
$$
as soon as $W_0=1$ and $f_n((\textsl{x}_0,\textsl{x}_1),\ldots,(\textsl{x}_n,\textsl{x}_{n+1}))=\Fa_n(\textsl{x}_0,\ldots,\textsl{x}_n)$.

\emph{Particle absorption models}. Consider a particle in an absorbing random medium, whose
successive states $\left(  X_{n}^{\prime}\right)  _{n\geq0}$ evolve according to a
Markov kernel $M^{\prime}$ on some state space $E$. At time $n$, the particle is absorbed with probability
$1-G^{\prime}\left(  X^{\prime}_{n}\right)$, where $G^{\prime}$ is a $\left(  0,1\right)  $-valued
potential function. Let $G^{\prime}_{n}:=G^{\prime}$ for all $n\geq0$ and $M_{n}^{\prime}:=M^{\prime}$ the Markov transitions of the chain $X^{\prime}_n$ for all
$n\geq1$. Then, the connection with the Feynman--Kac formalism is the following.
Let $T$ be the absorption time of the particle. Then, formula (\ref{FK-intro-statements}) has the following interpretation  
$$\eta_{n}=\mbox{\rm Law}\left(  X_{n}~|~T\geq n\right).$$

\emph{Physics and Chemistry}.
In these contexts, Feynman--Kac models are
widely used to describe molecular systems. The discrete generation
Feynman-Kac measures (\ref{FK-intro-statements}) can be interpreted as the solution of a discrete-time approximation of an
imaginary time Schrödinger equation. If we set $Id$ as the identity operator, then the Markov kernel $M^{\prime}\simeq_{\Delta t\downarrow 0}
Id+L~\Delta t$ of the chain $X^{\prime}_n$ corresponds to the discretization of a continuous-time stochastic
process $X^{\prime}_{t}$ with infinitesimal generator $L^{\prime}$, $G^{\prime}_{n}=e^{-V~\Delta t}$,
where $V$ is a potential energy, and $t_{n+1}-t_n:=\Delta t\ll1$, is a discretization time-step associated with some time mesh $t_n=n\lfloor t/n\rfloor $.
Replacing  the chain $X^{\prime}_n$ in (\ref{FK-intro-statements}) by the random state of the discrete time approximation model $X^{\prime}_{t_n}$, we have
\begin{eqnarray*}
\eta_{t_n}(f)&\propto&\mathbb{E}\left(  f(X^{\prime}_{t_0},\ldots,X^{\prime}_{t_n})~\exp{\left\{  -\sum_{0\leq t_k<t_n}V(X^{\prime}_{t_k}) (t_{k+1}-t_k)\right\}  }\right)\\&\simeq_{\Delta t\downarrow 0}&\mathbb{E}\left(  f(X^{\prime}_{s}, s\leq t)~\exp{\left\{  -\int_{0}^{t}%
V(X^{\prime}_{s})ds\right\}  }\right)  
\end{eqnarray*}
The marginal $\gamma_t$ w.r.t. the terminal time $t$ of the above measures
is often defined, in a weak sense, by
the imaginary time Schrödinger equation
\[
\frac{d}{dt}\gamma_{t}(f)=\gamma_{t}(L^{V}(f))\quad\mbox{\rm with}\quad L^V(f)=L(f)-Vf
\]

 For a more thorough
discussion of these continuous-time models and their applications in chemistry
and physics, see~\cite{caffarel,cances-tony,djlmp-2013,tony3,tony,tony-rs,rousset}, the recent monograph~\cite{d-2013}, and the
references therein.

To approximate these measures we 
run a mutation-selection genetic type particle model with $N$ individuals.   The $n$-th selection is associated with the fitness functions $G_n^{\prime}$, and the mutation transition is dictated by the transitions of the chain $X_n^{\prime}$. The empirical occupation measures of the genealogical trees associated with the resulting 
genetic population model approximate $\eta_n$ as the size of the population tends to infinity. Section~\ref{sec-mean-field-intro} gives a more precise description  of these particle models, including their mean field particle interpretation. We also mention that for unit potential functions $G_n^{\prime}=1$
the particle model reduces to the so-called neutral genetic model~\cite{dmpr-2009}.

Let  $\eta^{\prime}_n$ be the $n$-th marginal of the measure $\eta_n$. The sequence of measures $\eta^{\prime}_n$ satisfies a nonlinear updating-correction evolution equation (a.k.a. the filtering equation in signal processing, or the infinite population model in the scientific computing literature). We let $\eta^{\prime,\textsl{x}}_{p,n}$ be the solution of these equations starting
at the Dirac measure $\delta_{\textsl{x}}$ at time $p\leq n$, for some state $\textsl{x}\in S^{\prime}_p$. 
Let us introduce the regularity condition (well-suited, as it appears from its very definition, to the study of ergodic properties of Feynman--Kac models)
\begin{equation}
(\mbox{\rm H})\qquad \sup_{p\geq 0}\sum_{n\geq p}\left(\overline{g}^{\prime}_n-1\right)\beta_{p,n}<\infty,
\label{H}%
\end{equation}
with
$$
\beta_{p,n}:=\sup_{\textsl{x},\textsl{y}}\left\Vert \eta^{\prime,\textsl{x}}_{p,n}-\eta^{\prime,\textsl{y}}_{p,n}\right\Vert_{\tiny tv} \quad\mbox{\rm and}\quad\overline{g}^{\prime}_n=\sup_{\textsl{x},\textsl{y}}(G_n^{\prime}(\textsl{x})/G_n^{\prime}(\textsl{y}))<\infty,
$$
where $\Vert\point\Vert_{\tiny tv}$ stands for the total variation norm (cf. Section~\ref{notation-section}).

Condition (H) is met as soon as $\sup_{n\geq 0}\overline{g}^{\prime}_n<\infty$ and $\beta_{p,n}\leq a~e^{-\lambda(n-p)}$, for some finite constants
$0<a,\lambda<\infty$. 
It is related to the stability properties of the limiting Feynman-Kac measures. 
It ensures that local errors 
do not propagate w.r.t. the time horizon. In addition, in this situation  the bias and variance of the genealogical tree occupation measures on a time horizon $n$ are of order $n/N$ (see for example Corollary 14.3.6 and Corollary 15.2.2 in~\cite{d-2013}). 
For example,
in nonlinear filtering problems, the condition ensures that the optimal filter
forgets any erroneous initial condition. In quantum physics and molecular chemistry, the ground state energies of quantum systems are described by the limiting
 measures $\eta^{\prime}_{\infty}$ of the Feynman-Kac flow $\eta^{\prime}_n$ as $n\rightarrow\infty$. In computational physics and mathematical biology literatures
these limiting measures are sometimes called the quasi-invariant measures of the Yaglom limits. In this context, the regularity condition  (H) ensures the existence and uniqueness of these measures as well as the rate of convergence of the flow of measures towards their equilibrium.

For time homogeneous models, condition (\ref{H}) is met as soon as 
the Markov transition $M$ of the chain $X^{\prime}_n$ satisfies the minorisation condition $M^m(\textsl{x},d\textsl{z})\geq \epsilon M^m(\textsl{y},d\textsl{z})$ for some $m\geq 1$ and $\epsilon>0$, and for any $\textsl{x},\textsl{y}\in S^{\prime}$. This condition is met for regular nonlinear filtering and particle absorption models in compact spaces (see for instance~\cite{dm-2000,dm-toulouse}). It is met for the importance sampling target distributions 
(\ref{ex-is}) as soon as the transition probabilities of $\Xa_k$ satisfy the minorisation condition and $W_k$ are lower and upper bounded. 
It is also met for (\ref{ex-filter}) {\em for any integrable likehood function $\textsl{x}\mapsto g_{k}(\textsl{y}_{k},\textsl{x})$},
as soon as the Markov transition of the signal are chosen so that $K_k(\textsl{x},d\textsl{z})\geq \epsilon~K_k(\textsl{y},d\textsl{z})$ for some $\epsilon>0$ and for any $x,y\in E$. For instance, this condition is met for the nonlinear filtering model (\ref{ex-filter-1d}) for bi-Laplace random variables $\Wa_n$~\cite{dm-2000}. It is also met for Gaussian random variables $\Wa_n$ as soon as the drift function $a$ is constant outside some compact interval. 

\subsection{Frozen trajectories and PMCMC}\label{sub:frozen}

We fix now the size of the system $N$ as well as the time horizon $n\geq 0$. 
Let us run the $N$-genetic particle model defined as above, but with a given \it frozen \rm trajectory $\YY_0^{(N)}:=\textsl{x}=(\textsl{x}_p)_{0\leq p\leq n}\in S_n$ up to a given time horizon $n$. Details will be given later on the proper way to define this frozing process. For the time being, let us describe it informally:
during the $p$-th selection the particles may select the $p$-th coordinate $x_p$ of the frozen path, with $0\leq p\leq n$.  
At the moment, we simply defined a biased version of the Feynman--Kac method which may be expected to behave properly statistically only under some asumptions on the distribution of the frozen trajectory. 

The key idea underlying particle Markov chain MC methods is that this frozing technique allows to define a new Markov process on the space of trajectories that (under some reasonable assumptions) has the targeted Feynman--Kac measure $\eta_n$ as its invariant distribution. This allows to use this new Markov chain to simulate $\eta_n$ in a MCMC way.

The definition of the new Markov process runs as follows.
At the terminal time $n$, we select uniformly
one the $N$ ancestral lines $\YY_1^{(N)}=\textsl{y}$ of the corresponding genealogical tree. More generally, we let $(\YY_k^{(N)})_{k\geq 0}$ be a Markov chain on $S_n$ with transition probabilities
\begin{equation}\label{def-KK-intro}
\KK^{(N)}_n(f_n)(\textsl{x}):=\EE\left(f_n(\YY^{(N)}_1)~|~\YY^{(N)}_0=\textsl{x}\right)
\end{equation}
To simplify the notation, we suppress the index $(\point)^{(N)}$ and write $\YY_k$ and $\KK_n$, instead of 
 $\YY^{(N)}_k$ and $\KK^{(N)}_n$. 
 
 The diagram below provides a realization of the transition $\YY_0\leadsto \YY_1$
 for $N=3$ particles and a time horizon $n=3$. 
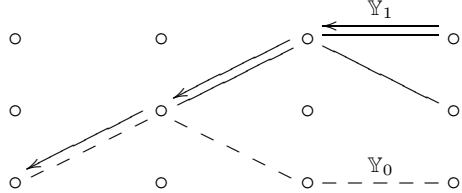
\begin{center}
\vskip-.3cm\hskip.2cm
\xymatrix@C=4em@R=1.5em{ 
\circ&\circ&\circ\ar@<.2ex>@{-}[ld]\ar@<-.5ex>@{->}[ld]\ar@{-}[rd]&\circ\ar@<-.25ex>@{-}[l]\ar@<-1ex>@{->}[l]_{\YY_1}\\ 
\circ&\circ\ar@<.2ex>@{--}[ld]\ar@<-.6ex>@{->}[ld] &\circ&\circ\\  
\circ&\circ&\circ\ar@{--}[r]^{\YY_0}\ar@<.25ex>@{--}[lu]&\circ
                                 }
\end{center}

The  goal of the present article is to build on the previously obtained results in the literature (as observed in the introduction to which we refer for references, in spite of a wide range of applications, the litterature on the \it theoretical \rm properties of these models is still very limited !) and analyse this Markov process from the point of view of Feynman--Kac models.
 
Since the theoretical analysis of PMCMC models under strong regularity assumptions of the potential functions is already a challenging task from the Feynman--Kac perspective, we have decided to refrain from proving results in the most general possible framework. Namely,
to avoid unnecessary technical discussions, we will often assume that the potential functions
$G_n$ are upper and lower bounded by some finite positive constant.  
However, in view of the discussion above on condition (H), we expect the boundedness assumption on the potential functions to be 
relaxed, at the cost of increased length and complexity of the proofs.

Let us mention, for example, that the algebraic polynomial developments presented in 
this article involve the positivity of certain integral operators (indexed by infected and coalescent forests). In this context, these rather strong regularity properties of the potential functions could be relaxed so that these integral operators still have a finite norm. For instance, Lemma~\ref{letruc} allows to extend the analysis to bounded potential functions which are not necessarily lower bounded. Besides, the extension of the results presented in this article
to more general models, including indicator type 
functions and unbounded potential functions, can also be analyzed using the methodologies developed in~\cite{d-2004}
(see for instance sections 2.3, 2.4, 3.5.2, and Section 7.2.2). 

The article will express the rate of convergence to equilibrium of the PMCMC chain $\YY_k$ presented above in terms of powers $(c(n)/N)^k$ of some ratio 
depending on some constant $c(n)$ that, in turns, depends on the time horizon $n$ of the target Feynman-Kac measures (\ref{FK-intro-statements}). The non-asymptotic estimates
derived in the article are valid for any bounded potential functions
$G_k$ and for any Markov chain $X_k$ s.t. $\eta_n(G_n)>0$. The constant $c(n)$ is expressed in terms of the norm of the potential functions and the quantities $\eta_k(G_k)$ (see, for instance the statement of Theorem~\ref{theo-taylor-intro}). 

Without any additional regularity condition these constants $c(n)$ grow exponentially w.r.t. the parameter $n$. Besides duality results on Feynman--Kac models of general (and, we believe, fundamental) theoretical interest, one of the main purposes of the article will be to show that $n\mapsto c(n)$ grow linearly w.r.t. the time $n$ as soon as the Feynman-Kac model
satisfy some natural stability conditions.

To describe the main results of the article with some precision, recall the notion of the differential for sequences of measures introduced in~\cite{dpr-2009}. We let $\mu^N$ be a uniformly bounded sequence of measures on some measurable state space $S$ in the sense that $\sup_{N\geq 1}\Vert\mu^N\Vert_{\tiny tv}<\infty$. The sequence $\mu^N$ is said to converge 
strongly to some measure $\mu$, as $N\uparrow\infty$ if we have $\lim_{N\uparrow\infty}\mu^N(f)=\mu(f)$, for any bounded measurable function $f$.
In this case, the discrete derivative of $\mu^N$ is defined by
$$
\partial \mu^N:=N~\left(\mu^N-\mu\right)
$$
We say that $\mu^N$ is differentiable whenever $\partial\mu^N$ is uniformly bounded and it strongly converges to some signed
measure $d^{(1)}\mu$, as $N\uparrow\infty$. When $\partial\mu^N$ is differentiable, with a discrete derivative writtem $\partial^{(2)}\mu^N$ we can define its derivative,
 denoted by $d^{(2)}\mu$, and so on. A mapping $N\mapsto \mu^N$ that is differentiable up to some order $l$ can be written as
$$
\mu^N=\sum_{0\leq k\leq l}~\frac{1}{N^k}~d^{(k)}\mu+\frac{1}{N^{l+1}}~\partial^{(l+1)}\mu^N
$$
with 
the convention $d^{(0)}\mu=\mu$. We easily extend these definitions to a sequence of integral operators $Q^N$ and a sequence of functions $f^N$. Here, we denote the corresponding differentials by $d^{(l)}Q$ and $d^{(l)}f$. 

We consider now a genetic model with $N$ individuals on $S^{\prime}_n$, with mutation transitions $M^{\prime}_n$ and selection fitness function $G^{\prime}_n$. 
We let  $\chi^{\prime}_k:=\left(\chi^{\prime i}_k\right)_{1\leq i\leq N}$ be the population after the $k$-th mutation. 
 (initially we start with $N$ independent copies of $X^{\prime}_0$).
We also denote by
$
\boldsymbol{\chi_n}=\left(\chi_k\right)_{0\leq k\leq n}\in \boldsymbol{\Sa_n}:=\prod_{0\leq k\leq n}S_k^N
$,  the ancestral lines $\chi_n=\left(\chi_n^i\right)_{1\leq i\leq N}$
 of the $N$ individuals $\chi^{\prime}_n:=\left(\chi^{\prime i}_n\right)_{1\leq i\leq N}$. We also  
let  $\XX_n$ be a randomly chosen ancestral line with the uniform distribution $\frac{1}{N}\sum_{1\leq i\leq N}\delta_{\chi_n^i}$.

 Finally, we consider the probability distributions $\pi_n$ on $\left(S_n\times \boldsymbol{\Sa_n}\right)$ given
  for any bounded measurable function ${f_n}$ on $\left(S_n\times\boldsymbol{\Sa_n} \right)$ by 
the formula
$$
 \pi_n(f_n)\propto\EE\left({f_n}(\XX_n,\boldsymbol{\chi_n})~{\Za_n(\chi)}\right)
 $$
with
 $$
 \Za_n(\chi)=\prod_{0\leq k<n}\Ga_k(\chi_k)\quad\mbox{\rm with}\quad \Ga_k(\chi_k)=\frac{1}{N}\sum_{1\leq i\leq N}G_k(\chi^i_k) .
 $$ 
The transition probabilities of the Gibbs sampling of the multivariate distribution $\pi_n$ on the product space $\left(S_n\times\boldsymbol{\Sa_n} \right)$
are described by the synthetic diagram
$$
\left\{
\begin{array}{rcl}
\XX_n&=&\textsl{x}\\
\boldsymbol{\chi_n}&=&x
\end{array}
\right\}\rightarrow \left\{
\begin{array}{rcl}
\overline{\XX}_n&=&\overline{\textsl{x}}\sim \left(\XX_n~|~\boldsymbol{\chi_n}=x\right)\\
\boldsymbol{\chi_n}&=&x
\end{array}
\right\}\rightarrow
\left\{
\begin{array}{rcl}
\overline{\XX}_n&=&\overline{\textsl{x}}\\
\boldsymbol{\overline{\chi}}_{\bf n}&=&\overline{x}\sim \left(\boldsymbol{\chi_n}~|~\XX_n=\overline{\textsl{x}}\right)
\end{array}
\right\}
$$

In the above display, $\left(\XX_n~|~\boldsymbol{\chi}_n\right)$ and $\left(\boldsymbol{\chi}_n~|~\XX_n\right)$ 
is a shorthand notation for the $\pi_n$-conditional distributions of $\XX_n$ given $\boldsymbol{\chi}_n$, and $\boldsymbol{\chi}_n$ given $\XX_n$.

The first transition of the Gibbs sampler reduces to the uniform sampling of an ancestral line. Next, we present a duality formula that shows that
the second transition amounts of sampling a genetic particle model with a frozen ancestral line.

Using this notation, one of our main result can be stated as follows.

\begin{theo}\label{theo-key-introduction}
 For any bounded measurable function $f_n$ on the product space $\left(S_n\times \boldsymbol{\Sa_n}\right)$ (symmetric on the product spaces $S^{N}_k$) we have the duality formula
\begin{equation}\label{duality-intro-ref}
\EE\left(f_n(\XX_n,\boldsymbol{\chi_n})~{\Za_n(\chi)}\right)=\EE\left(f_n(X_n,
\boldsymbol{\Xa_{n}})~{Z_n(X)}\right)
\end{equation}
where
 $\boldsymbol{\Xa_{n}}=\left({\Xa_{k}}\right)_{0\leq k\leq n}$
stands for all the ancestral lines at each level $k$, of a conditional $N$-genetic particle model with a given frozen trajectory $X_n$. Consequently, the PMCMC chain 
$\YY_k$ with ancestral sampling defined in (\ref{def-KK-intro}) coincides with the first coordinate of the Gibbs sampler of the target distribution $\pi_n$. In addition, the PMCMC chain $\widetilde{\YY}_k$ avoiding the last frozen state 
coincides with the first coordinate of the Gibbs sampler targeting the marginal $\widetilde{\pi}_n$ of $\pi_n$ on $\left(S_n\times \boldsymbol{\Sa_{n-1}}\right)$. 
 
Furthermore, the Markov transitions $\KK_n$ are differentiable at any order $l\geq 1$ with $d^{(0)}\KK_n(f_n)=\eta_n(f_n)$.
If we further assume that the regularity condition  $(H)$ stated in
(\ref{H}) is satisfied, we have, for any bounded measurable function $f$ on the path space $S_n$ s.t. $\Vert f_n\Vert\leq 1$, the non-asymptotic estimates
$$
\forall 1\leq k\leq l\qquad \left\Vert d^{(k)}\KK_n(f)\right\Vert\leq (cnk^2)^{k}~~\mbox{and}~~
\left\Vert \partial^{(l+1)} \KK_n(f)\right\Vert
\leq  (cn(l+1)^2)^{l+1}
$$
as soon as $N>cn(l+1)^2$, for some finite constant $c<\infty$. 

Consequently, there exists some finite constant $c<\infty$ such that
for any $m\geq 1$, $\textsl{x}\in S_n$, and any $N>cn$ we have
\begin{equation}\label{key-illustration}
\left\vert\EE\left(f(\YY_m)~|~\YY_0=\textsl{x}\right)-\eta_n(f)\right\vert\leq (cn/N)^m\quad\left(\longrightarrow_{\min{(N,m)}\rightarrow\infty}0\right)
\end{equation}
and for any $p\geq 1$, we have a sharp non-asymptotic estimate of the $\LL_p(\eta_n)$-mean error norms
$$
\left\vert\left\Vert \KK_n^m(f)-\eta_n(f)\right\Vert_{\LL_p(\eta_n)}-N^{-m}~\left\Vert  \left[d^{(1)}\KK_n\right]^m\!\!\!(f)\right\Vert_{\LL_p(\eta_n)}\right\vert\leq (cn/N)^{m+1}
$$
with the $m$-th iterates $\KK^m_n=\KK^{m-1}_n\KK_n$ and $ \left[d^{(1)}\KK_n\right]^m:=\left[d^{(1)}\KK_n\right]^{m-1}d^{(1)}\KK_n$ of the operators $\KK_n$ and $d^{(1)}\KK_n$.

\end{theo}

 The first assertion is proved in Section~\ref{conditiona-sec} (see, for instance, Thm.~\ref{theo-transport} and
Corollary~\ref{cor-pGibbs} for Feynman-Kac models (\ref{FK-intro-statements})
 on general state spaces $S_n$,
  and Section~\ref{historical-sec} for models associated with an historical process). This result shows that the chains $\YY_k$ and $\widetilde{\YY}_k$ reduce to the first coordinate
 of a couple of Gibbs samplers with target measure $\pi_n$ and $\widetilde{\pi}_n$. By (\ref{duality-intro-ref}) we conclude that $\YY_k$ and $\widetilde{\YY}_k$ are reversible w.r.t. the target measure $\eta_n$. For instance, using the duality formula (\ref{duality-intro-ref}), the transition probabilities of the Gibbs sampling of the multivariate distribution $\pi_n$ on the product space 
$\left(S_n\times\boldsymbol{\Sa_n} \right)$
are described by the synthetic diagram
 $$
\left\{
\begin{array}{rcl}
\XX_n&=&\textsl{x}\\
\boldsymbol{\chi_n}&=&x
\end{array}
\right\}\rightarrow \left\{
\begin{array}{rcl}
\overline{\XX}_n&=&\overline{\textsl{x}}\sim \left(\XX_n~|~\boldsymbol{\chi_n}=x\right)\\
\boldsymbol{\chi_n}&=&x
\end{array}
\right\}\rightarrow
\left\{
\begin{array}{rcl}
\overline{\XX}_n&=&\overline{\textsl{x}}\\
\boldsymbol{\overline{\chi}}_{\bf n}&=&\overline{x}\sim \left(\boldsymbol{\Xa_n}~|~X_n=\overline{\textsl{x}}\right)
\end{array}
\right\}
$$

The end of the Theorem is a direct consequence of a more general theorem, Theorem~\ref{theo-taylor-intro}, describing non asymptotic Taylor series without the regularity condition
$(H)$.

The $k$-th derivative integral operators $d^{(k)}\KK_n$
will be described explicitly in Section~\ref{ifexp-ref} in terms of Feynman-Kac semigroups parametrized by coalescent and decorated forests of length $n$ with $(2k+1)$ edges, and less than $k$ coalescences and infections.

Notice that (\ref{key-illustration}) can be used to estimate the Lyapunov exponent of the distribution semigroup of the PMCMC chain $\YY_m$; that is, 
$$
\liminf_{m\rightarrow\infty}-\frac{1}{m}\log{\left\Vert\mbox{\rm Law}(\YY_m)-\eta_n 
\right\Vert}_{\tiny tv}\geq \log{(N/(cn))}
$$

A more precise description of these non-asymptotic Taylor expansions, including a series of illustrations of the impact of these results
in the estimation of the variance and the Dobrushin contraction coefficient of these models 
is provided in Section~\ref{taylor-invariant} (see, for instance, Theorem~\ref{theo-taylor-intro} and the discussion that follows).

\subsection{Backward Sampling}

We end this section with a discussion on PMCMC models based on backward particle samplers. To describe with some precision these models, 
we further assume that the integral operators $G^{\prime}_{k}(x)M^{\prime}_{k+1}(x,dy)$ have a density $H^{\prime}_{k+1}(x,y)$ w.r.t. some reference measure. In this situation, the Feynman-Kac measure $\eta_n$ defined in (\ref{FK-intro-statements}) can be interpreted as the distribution of a nonlinear backward Markov chain model. 
The particle interpretation of these distributions is defined as follows (a more detailed description of these nonlinear backward models is provided in section~\ref{sec-pathspace} (see, for instance, (\ref{backward-2}) and (\ref{unbiased-estimators-ref})).

The historical process 
 of all populations of ancestors $$\boldsymbol{\chi^{\prime}_n}:=\left(\chi_k^{\prime}\right)_{0\leq k\leq n}\in \boldsymbol{\Sa_n^{\prime}}:=\prod_{0\leq k\leq n}S^{\prime N}_k, \quad\mbox{\rm with}\quad \chi^{\prime}_k:=\left(\chi^{\prime i}_k\right)_{1\leq i\leq N}\in S^{\prime N}_k$$ at every level $k$ can be interpreted as the complete ancestral tree 
(without the tree structure) of the genetic model discussed above. The backward particle model is a Markov chain running backward in time with the state spaces  $\left\{\chi^{\prime i}_k; 1\leq i\leq N\right\}$ at each level  $k$.  The initial state of the chain takes the value $\chi^{\prime i}_n$ with probability $1/N$, with $1\leq i\leq N$.
 Then, at each level $0\leq k<n$ the (conditional) probability to go from state $\chi^{\prime i}_{k+1}$ to state $\chi^{\prime j}_k$  is proportional to $H^{\prime}_{k+1}\left(\chi^{\prime j}_k,\chi^{\prime i}_{k+1}\right)$.   We denote by ${\XX}_n^{\flat}$ a backward randomly chosen ancestral line.

Running the $N$-genetic particle model with a given frozen path $\YY^{\flat}_0:=\textsl{x}=(\textsl{x}_p)_{0\leq p\leq n}\in S_n$ up to a given time horizon $n$,
we let $\YY^{\flat}_1:=\textsl{y}$ be an ancestral line randomly chosen with the backward Markov chain model discussed above. The initial value of the chain
is one of the states at time $n$ of the genetic model with the frozen trajectory (including $\textsl{x}_n$).  
Iterating this transition, we define  a Markov chain  $(\YY^{\flat}_k)_{k\geq 0}$ on $S_n$.   Section~\ref{g+b-ref} provides a more detailed description of the Markov transition of this backward particle MCMC model. 

The diagrams below provide realizations of the transitions $\YY_0^{\flat}\leadsto \YY_1^{\flat}$
 for $N=3$ particles and a time horizon $n=3$. 
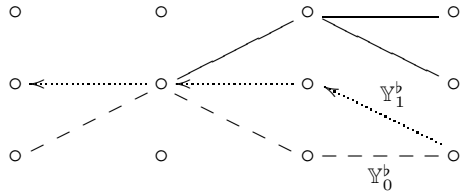
\begin{center}
\hskip.2cm
\xymatrix@C=4em@R=1.5em{ 
\circ&\circ&\circ\ar@<.2ex>@{-}[ld]\ar@{-}[rd]&\circ\ar@<.4ex>@{-}[l]\\ 
\circ&\circ\ar@<.2ex>@{--}[ld]\ar@{.>}[l]&\circ\ar@{.>}[l]&\circ\\  
\circ&\circ&\circ\ar@{--}[r]_{\YY_0^{\flat}}\ar@<.25ex>@{--}[lu]&\circ\ar@<-.4ex>@{.>}[lu]_{\YY_1^{\flat}} 
                                 }
\end{center}

In the context of backward PMCMC models, one of our main results can be stated as follows.

\begin{theo}
Given the complete ancestral tree $\boldsymbol{\chi^{\prime}_n}$, the ancestral lines $\chi_n^i=(\chi^i_{k,n})_{0\leq k\leq n}$ are copies of the 
backward trajectory $\XX^{\flat}_n$ starting at the terminal state $\chi^{\prime i}_n=\chi^i_{n,n}$, with $1\leq i\leq N$. That is,
\begin{equation}\label{equivalence-intro-key}
\mbox{\rm Law}\left(\XX_n~|~\boldsymbol{\chi^{\prime}_n}\right)=\mbox{\rm Law}\left(\XX_n^{\flat}~|~\boldsymbol{\chi^{\prime}_n}\right)
\end{equation}

Consequently, the PMCMC model  with backward sampling coincide with the Gibbs sampler targeting the $(S_n\times\boldsymbol{\Sa_n^{\prime}})$-marginal distribution $\pi_n^{\prime}$ 
of the measure $\pi_n$ defined in (\ref{duality-intro-ref}). 
\end{theo}

This result is proved in Section~\ref{historical-sec}. Theorem~\ref{theo-ancestral-backward} also provides an interpretation
 of the conditional genealogical trees $(\chi_k)_{0\leq k\leq n}$ and $(\Xa_k)_{0\leq k\leq n}$ given the complete tree of all ancestors in terms of a Markov chain
with elementary transitions defined by  backward ancestor sampling.
 The diagram below provides a realization of the transitions $\chi_2\leadsto \chi_3$ given $\chi^{\prime}_3$
 for $N=3$ particles. 
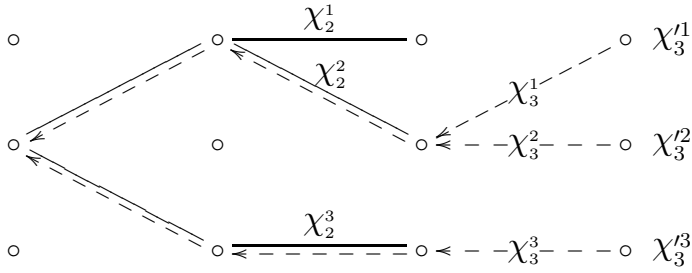
\begin{center}
\vskip-.3cm\hskip.2cm
\xymatrix@C=6em@R=2em{ 
\circ&\circ\ar@<.3ex>@{-}[rd]\ar@<-.3ex>@{-}[ld]\ar@<.3ex>@{-->}[ld]&\circ\ar@{-}[l]_{\chi^1_2}&\circ\ar@{-->}[ld]|--{\chi^1_3}&\hskip-4cm\chi^{\prime 1}_3\\ 
\circ&\circ&\circ\ar@<.3ex>@{-->}[lu]_{\chi^2_2}&\circ\ar@{-->}[l]|--{\chi^2_3}&\hskip-4cm\chi^{\prime 2}_3\\  
\circ&\circ\ar@<.3ex>@{-->}[lu]\ar@<-.2ex>@{-}[lu]&\ar@<.3ex>@{-->}[l]\ar@<-.4ex>@{-}[l]_{\chi^3_2}\circ&\circ\ar@{-->}[l]|--{\chi^3_3}&\hskip-4cm\chi^{\prime 3}_3
                                 }
\end{center}

The equivalence formulae (\ref{equivalence-intro-key})  between the ancestral and the backward samplers are proved in Section~\ref{historical-sec}. Combining (\ref{duality-intro-ref}) and (\ref{equivalence-intro-key}), we prove
the duality formula
\begin{equation}\label{key-form-introduction}
\pi_n^{\prime}(f_n)\propto \EE\left(f_n\left({\XX_n^{\flat}},\boldsymbol{\chi^{\prime}_n}\right)~\Za_n^{\prime}(\chi^{\prime})\right)=\EE\left(f_n(X_n,
\boldsymbol{\Xa_{n}^{\prime}})~{Z_n^{\prime}(X^{\prime})}\right)
\end{equation}
for any bounded measurable function $f_n$ on the product space (symmetric on the product spaces $S^{\prime N}_k$), where
 $\boldsymbol{\Xa^{\prime}_{n}}=\left({\Xa^{\prime}_{k}}\right)_{0\leq k\leq n}$
stands for the populations of the $N$-genetic particle model with a given frozen trajectory $X_n$ and $\Za_n^{\prime},\ Z_n^{\prime}$ are defined as usual.

The theorem also shows that the Markov chain $\YY^{\flat}_k$ reduces to the first coordinate
 of a Gibbs sampler with target measure $\pi_n^{\prime}$.  By (\ref{key-form-introduction}) 
  and (\ref{duality-intro-ref}) we conclude that $\YY^{\flat}_k$  is reversible w.r.t. the target measure $\eta_n$.

By the duality formula (\ref{key-form-introduction}), the transition probabilities of the Gibbs sampling of the multivariate distribution $\pi_n^{\prime}$ on the product space 
$\left(S_n\times\boldsymbol{\Sa_n^{\prime}} \right)$
are described by the synthetic diagram
$$
\left\{
\begin{array}{rcl}
\XX_n^{\flat}&=&\textsl{x}\\
\boldsymbol{\chi^{\prime}_n}&=&x
\end{array}
\right\}\rightarrow \left\{
\begin{array}{rcl}
\overline{\XX}^{\flat}_n&=&\overline{\textsl{x}}\sim \left(\XX_n^{\flat}~|~\boldsymbol{\chi_n^{\prime}}=x\right)\\
\boldsymbol{\chi^{\prime}_n}&=&x
\end{array}
\right\}\rightarrow
\left\{
\begin{array}{rcl}
\overline{\XX}_n^{\flat}&=&\overline{\textsl{x}}\\
\boldsymbol{\overline{\chi}_n^{\prime}}&=&\overline{x}\sim \left(\boldsymbol{\Xa^{\prime}_n}~|~X_n=\overline{\textsl{x}}\right)
\end{array}
\right\}
$$

\section{Feynman-Kac models: old and new}\label{sec-mean-field-intro}

This section collects first some basic notations used in this article. We recall then the definition and main properties of Feynman-Kac measures on their usual state and path spaces. The last paragraph introduces a particular Feynman-Kac model \cite{vddm-2014} well-suited to the mathematical analysis of PMCMC samplers. Although we will not develop further this point of view, the statistically-minded reader will note the analogy of the model with the ones familiar in $U$-statistics, in that it relies strongly on properties of symmetric functions on the space of samples of a target distribution.

\subsection{Notation}\label{notation-section}

Given some measurable space $S$ we denote respectively by $\Ma(S)$,
$\Pa(S)$ and $\Ba(S)$, the set of finite signed measures on $S$, the convex subset of probability measures, and the Banach space of bounded measurable functions equipped with the uniform norm $\Vert f\Vert=\sup_{\textsl{x}\in S}\vert f(\textsl{x})\vert$.  

The total variation norm on measures $\mu\in \Ma(S)$ 
 is defined by $$
 \Vert\mu\Vert_{\tiny tv}:=\sup_{f\in\Ba(S)~:~\|f\|\leq 1}\vert \mu(f)\vert
\quad\mbox{\rm with the Lebesgue integral}\quad
\mu(f):=\int~\mu(d\textsl{x})~f(\textsl{x}) $$
 We also denote by  $\delta_a$ the Dirac measure at some state $a$, so that $\delta_a(f)=f(a)$. We say that $\nu\leq \mu$ as soon
 as $\nu(f)\leq \mu(f)$ for any non-negative function $f$.

A bounded integral operator $Q(x,dy)$ between the measurable spaces $S$ and $S^{\prime}$ is defined for any $f\in\Ba(S^{\prime})$
by the measurable function $Q(f)\in\Ba(S)$ defined by
$$
Q(f)(\textsl{x}):=\int~Q(\textsl{x},d\textsl{y})~f(\textsl{y}) 
$$
The operator $Q$ generates a dual operator $\mu\in \Ma(S)\mapsto \mu Q\in\Ma(S^{\prime})$ by the dual formula
$
(\mu Q)(f)=\mu(Q(f))
$. 

When a bounded integral operator  $M$ from a state space $S$ into a possibly different state space $S^{\prime}$ has a constant mass, that is, when $M(1)\left(  x\right)  =M(1)\left(
y\right)  $ for any $(x,y)\in S^{2}$, the operator $\mu\mapsto\mu M$ maps the set
$\mathcal{M}_{0}(S)$ of measures $\mu$ on $S$ with null mass $\mu(1)=0$ into $\mathcal{M}_{0}(S^{\prime})$. In this situation, we let
$\beta(M)$ be the Dobrushin coefficient of a bounded integral operator $M$
defined by the formula
$$\beta(M):=\sup{\ \{\mbox{\rm osc}(M(f))\;;\;\;f~\mbox{\rm s.t.}~\mbox{\rm osc}(f)\leq 1\}} $$ where $\mbox{\rm osc}(f):=\sup_{x,y}|f(x)-f(y)|$ stands for the oscillation of some
function.

When $M$ is a Markov transition, $\beta(M)$ coincides with the Dobrushin contraction parameter (a.k.a. the Dobrushin ergodic coefficient)  defined by 
$$\beta(M)=\sup_{\mu,\nu}{\left({\Vert \mu M-\nu M\Vert_{\tiny tv}}/{\Vert\mu-\nu\Vert_{\tiny tv}}\right)}=2^{-1}\sup_{x,y}{\Vert M(x,\cdot)-M(y,\cdot)\Vert_{\tiny tv}}$$

The $q$-tensor product of $Q$ is the integral operator defined for any $f\in \Ba(S^q)$ by
$$
Q^{\otimes q}(f)(\textsl{x}^1,\ldots,\textsl{x}^q):=\int~\left\{\prod_{1\leq i\leq q}Q(\textsl{x}^i,d\textsl{y}^i)\right\}~f(\textsl{y}^1,\ldots,\textsl{y}^q)
$$
We also denote by $Q_1Q_2$ the composition of two operators defined by
$$
(Q_1Q_2)(\textsl{x},d\textsl{z}):=\int Q_1(\textsl{x},d\textsl{y})Q_2(\textsl{y},d\textsl{z})
$$

The Boltzmann-Gibbs transformation $\Psi_G~:~\eta\in\Pa(S)\mapsto \Psi_{G}(\eta)\in \Pa(S)$ 
associated with some positive function $G$ on some state space $S$ is defined by
$$
\Psi_{G}(\eta)(d\textsl{x}):=\frac{1}{\eta(G)}~G(\textsl{x})~\eta(d\textsl{x})
$$
We also denote by $\#(E)$ the cardinality of a finite set
and we use the standard conventions $\left(
\sup_{\emptyset},\inf_{\emptyset}\right)=\left(-\infty,+\infty\right)$, and $\left(\sum_{\emptyset},\prod_{\emptyset}\right)=(0,1)$.

\subsection{Mean field particle models}

We consider a collection of non-negative bounded potential functions
$G_n$ on some  measurable state spaces $S_n$, with $n\in\NN$.
We also let $X_n$ be a Markov chain on $S_n$ with initial distribution $\eta_0\in\Pa(S_0)$ and
some Markov transitions $M_n$ from $S_{n-1}$ into $S_n$. The Feynman-Kac measures $(\eta_n,\gamma_n)$  associated with the parameters
$(G_n,M_n)$ are defined for any $f_n\in\Ba(S_n)$ by $\eta_n(f_n):={\gamma_n(f_n)}/{\gamma_n(1)}$ with
\begin{equation}\label{FK-def-intro-ref}
\gamma_n(f_n)=\EE\left(f_n(X_n)~Z_n(X)\right) \quad\mbox{\rm and}\quad Z_n(X)=\prod_{0\leq p<n}G_p(X_p)
\end{equation} 
It is implicitly assumed that $\gamma_n(1)>0$.
The evolution equations associated with these measures are given by
\begin{equation}\label{evolve-eq-intro-ref}
\gamma_{n+1}=\gamma_{n}Q_{n+1}\quad\mbox{and}\quad
\eta_{n+1}=\Phi_{n+1}(\eta_n):=\Psi_{G_{n}}(\eta_{n})M_{n+1}
\end{equation} 
with the integral operators $$
Q_{n+1}(\textsl{x}_{n},d\textsl{x}_{n+1})=G_{n}(\textsl{x}_{n})~
M_{n+1}(\textsl{x}_{n},d\textsl{x}_{n+1})$$

The unnormalized measures $\gamma_n$ can be expressed in terms of the normalized ones
using the well known product formula
$$
\gamma_n(f_n)=\eta_n(f_n)~\prod_{0\leq p<n}\eta_p(G_p)
$$
We also recall the semigroup decompositions
$$
\forall 0\leq p\leq n\qquad \gamma_n=\gamma_pQ_{p,n}\qquad\mbox{\rm and}\quad \eta_n=\eta_p\overline{Q}_{p,n}
$$
with the integral operators $Q_{p,n}=Q_{p+1}\ldots Q_n$, and the normalized semigroups
$$
\overline{Q}_{p,n}(f_n)(\textsl{x}_p)={{Q}_{p,n}(f_n)(\textsl{x}_p)}/{\eta_p{Q}_{p,n}(1)}=(\overline{Q}_{p+1}\ldots \overline{Q}_n)(f_n)(\textsl{x}_p)
$$
In the above display, $\overline{Q}_{p+1}$ stands for the collection of integral operators defined as $Q_{p+1}$ by replacing $G_{p}$
with the normalized potential functions $\overline{G}_p=G_p/\eta_p(G_p)$.  {We further assume that
\begin{equation}\label{Hprime}
(H^{\prime})\qquad\overline{g}_n=\sup_{0\leq p\leq q\leq n}\left\Vert\overline{Q}_{p,q}(1)\right\Vert<\infty
\end{equation}
This condition is clearly met when $\sup_{x,y}(G_n(x)/G_n(y))<\infty$. In addition, in the context of the Feynman-Kac models (\ref{FK-intro-statements})
 we have $\sup_{n\geq 0}\overline{g}_n:=\overline{g}<\infty$ as soon as the condition ({\rm H}) stated in (\ref{H}) is satisfied. A proof of this result is provided in Chapter 12.2.1 in~\cite{d-2013}.}

The mean field particle interpretation of the measures $(\eta_n,\gamma_n)$
starts with $N$ independent random variables $\xi_0:=(\xi^{i}_0)_{1\leq i\leq N}\in S_0^N$ with common law $\eta_0$.
The simplest way to evolve the population of $N$ individual (a.k.a. samples, particle, or walkers)  $\xi_n:=(\xi^{i}_n)_{1\leq i\leq N}\in S_n^N$ is to consider $N$ conditionally independent individuals $\xi_{n+1}:=(\xi^{i}_{n+1})_{1\leq i\leq N}\in S_{n+1}^N$ with common distribution 
\begin{equation}\label{mean-field-intro}
\Phi_{n+1}(m(\xi_n)
)\quad\mbox{\rm with}\quad m(\xi_{n}):=\frac{1}{N}\sum_{1\leq i\leq N}\delta_{\xi^{i}_{n}}
\end{equation}
This particle model (\ref{mean-field-intro}) is a genetic type particle model with a selection and a mutation transition dictated by the potential function $G_n$ and the Markov transition $M_{n}$. 

Loosely speaking, the model functions recursively as follows: starting from a sample $\xi_0^{(N)}$ at $t=0$ of the initial distribution $\eta_0$ (so that $m(\xi_0)\simeq_{N\uparrow\infty}\eta_0$), and assuming
$m(\xi_n)\simeq_{N\uparrow\infty}\eta_n$, then the population at time $(n+1)$ 
is formed with  $N$ "almost" independent samples w.r.t. $\eta_{n+1}$ so that $m(\xi_{n+1})\simeq_{N\uparrow\infty}\eta_{n+1}$. The reader is refered to \cite{d-2004} for details.

\subsection{Path space models}\label{sec-pathspace}

To illustrate the generality of the Feynman-Kac models discussed above, let us replace the 5-tuple $(G_n,M_n,Q_n,S_n,X_n)$ by its path-space analog $(\boldsymbol{G_n,M_n,Q_n,S_n,X_n})$. 
That is, in the constructions of the previous paragraph, each item of the first 5-tuple is going to be replaced by its path space analog:
$\boldsymbol{X}_n$ is the historical process associated to $X_n$,
\begin{equation}\label{path-space-intro-ref}
\boldsymbol{X}_n:=(X_0,\ldots,X_n)\in \boldsymbol{S}_n:=(S_0\times\ldots\times S_n).
\end{equation} 
We write $\boldsymbol{M}_n$ for the Markov transition of $\boldsymbol{X}_n$. 
The functions $\boldsymbol{G}_n$ on $\boldsymbol{S}_n$ only depend on the last coordinate and are defined by
$\boldsymbol{G}_n(\boldsymbol{X}_n):=G_n(X_n)$.

In general, in the article, a bold symbol will denote an element, function, measure... on a path space, even when the latter is considered as a state space --as in the present paragraph.
In particular, we let $(\boldsymbol{\gamma_n},\boldsymbol{\eta_n},\boldsymbol{\xi_n})$ be the Feynman-Kac measures and the particle model
defined as $(\gamma_n,\eta_n,\xi_n)$,
by replacing  $(G_n,M_n,Q_n,S_n,X_n)$ by 
$(\boldsymbol{G_n,M_n,Q_n,S_n,X_n})$. The two measures on the state space $\boldsymbol{S}_n$ are
given
for any $\boldsymbol{f_n}\in\Ba(\boldsymbol{S_n})$ by $\boldsymbol{\eta_n(f_n):={\gamma_n(f_n)}/{\gamma_n(1)}}$, with
\begin{equation}\label{FK-path-intro}
\boldsymbol{\gamma_n(f_n)}=\EE\left(\boldsymbol{f_n(X_n)}~Z_n(X)\right).
\end{equation}

By construction,  $(\gamma_n,\eta_n)$ are the $S_n$ marginals of the 
measures $(\boldsymbol{\gamma_n},\boldsymbol{\eta_n})$. 
The same property holds at the level of the particles of the two models.
To be more precise, we observe that the $i$-th path space particle 
\begin{equation}\label{FK-path-particle}
\boldsymbol{\xi}_n^i=\left(\xi^i_{0,n},\xi^i_{1,n},\ldots,\xi^i_{n,n}\right)\in \boldsymbol{S}_n:=(S_0\times\ldots\times S_n)
\end{equation}
of the particle model $\boldsymbol{\xi}_n$ can be interpreted as
 the  line of ancestors $\xi^i_{p,n}$ of the $i$-th individual $\xi^i_{n,n}$ 
at time $n$, at every level $0\leq p\leq n$, with $1\leq i\leq N$. 
The path space model $\boldsymbol{\xi_n}$
is called the genealogical tree model associated with the particle system $\xi_n$.

To distinguish these two Feynman-Kac models we adopt the following terminology.
The $3$-tuple $(\eta_n,\gamma_n,\xi_n)$ is called 
the Feynman-Kac particle model associated with the potential functions $G_n$ and the Markov transitions $M_n$ on the state spaces $S_n$.
The path space model $(\boldsymbol{\gamma_n},\boldsymbol{\eta_n},\boldsymbol{\xi_n})$ is called the historical version of $(\gamma_n,\eta_n,\xi_n)$.

Whenever the integral operators $Q_{n}$ have some densities $H_n$ w.r.t. some reference distributions  $\upsilon_n$
 on $S_n$, the path space measure $\boldsymbol{\eta}_n$ can be expressed in terms of the marginal measures 
 $(\eta_p)_{0\leq p\leq n}$
using the well known  backward formula
    \begin{equation}\label{backward-2}
\boldsymbol{\eta_n(\boldsymbol{d{\textsl{\bf x}}_n})}=\eta_n(d\textsl{x}_n)~\prod_{1\leq k\leq n}\LL_{k,\eta_{k-1}}(\textsl{x}_k,
d\textsl{x}_{k-1})
    \end{equation}
 with  the collection of Markov transitions $\LL_{n+1,\eta_n}$ from $S_{n+1}$ into $S_n$ defined by
\begin{equation}\label{def-LL}
\LL_{n+1,\eta_n}(\textsl{x}_{n+1},d\textsl{x}_{n})=\eta_n(d\textsl{x}_n)~H_{n+1}(\textsl{x}_n,\textsl{x}_{n+1})/{\eta_n\left(H_{n+1}(\point,\textsl{x}_{n+1})\right)}
\end{equation}
In the above displayed formula, $\boldsymbol{{\textsl{\bf x}}_n}$ stands for the
 trajectory 
$\boldsymbol{{\textsl{\bf x}}_n}=(\textsl{x}_0,\ldots,\textsl{x}_n)\in \boldsymbol{S}_n:=(S_0\times\ldots\times S_n)$.

In this setting, the two unbiased estimates of $\boldsymbol{\gamma}_n$ are defined by
\begin{equation}\label{unbiased-estimators-ref}
\forall i=1,2\qquad \boldsymbol{\gamma^{(N,i)}_n}=\left\{\prod_{0\leq p<n} m(\xi_p)(G_p)\right\}~\boldsymbol{\eta^{(N,i)}_n}
\end{equation}
with the couple of random measures $\left(\boldsymbol{\eta^{(N,1)}_n},\boldsymbol{\eta^{(N,2)}_n}\right)$ on $\boldsymbol{S}_n$ defined by
$$
\boldsymbol{\eta^{(N,1)}_n(\boldsymbol{d{\textsl{\bf x}}_n})}:=m(\boldsymbol{\xi}_n)(\boldsymbol{d{\textsl{\bf x}}_n})~~\mbox{\rm and}~~
\boldsymbol{\eta^{(N,2)}_n(\boldsymbol{d{\textsl{\bf x}}_n})}:=m(\xi_n)(d\textsl{x}_n)~\prod_{1\leq k\leq n}\LL_{k,m(\xi_{k-1})}(\textsl{x}_k,d\textsl{x}_{k-1}).
$$

\subsection{Many body Feynman-Kac models}\label{mbody-pm}

\subsubsection{Some terminology}
We fix the size $N$ of the particle model, and set $\Sa_n:=S_n^{[N]}$ for the $N$-th symmetric power of $S_n$: $S_n^{[N]}:=S_n\times ...\times S_n/\Sigma_N=S_n^N/\Sigma_N$, where we write $\Sigma_N$ for the symmetric group of order $N$.
The image in $\Sa_n$ of an ordered sequence $(\textsl{x}_1,...,\textsl{x}_n)\in S_n^N$ will be sometimes written with the set-theoretical notation 
$\{\textsl{x}_1,...,\textsl{x}_n\}$ to emphasize that the order of the $\textsl{x}_i$ does not matter, although we will also often identify $(\textsl{x}_1,...,\textsl{x}_n)$ with its image in $S_n^{[N]}$ without further notice when no confusion can arise.

For example with this slight abuse of notation, noticing for further use that the particle model $\xi_n$ can be viewed as a $\Sa_n$-valued Markov chain (since the distribution of the $\xi_n^i,\ i=1...N$ is $\Sigma_N$-invariant) we will have, for a function $f$ on $S_n^{[N]}$, $$f(\xi_n):=f(\{\xi_n^1,...,\xi_n^N\})=:f(\xi_n^1,...,\xi_n^N).$$

 In the further development of this section we use
 calligraphic letters such as $x_n$ and $y_n=\{\textsl{y}_n^i\}_{1\leq i\leq N}$ to denote states in the product spaces $\Sa_n=S_n^{[N]}$, and slanted roman letters such as $\textsl{x}_n$,
 $\textsl{y}_n$, $\textsl{z}_n$ to denote states in $S_n$.  
 The path sequences in the product spaces $\boldsymbol{\Sa}_n:=\prod_{0\leq p\leq n}\Sa_p$ and $\boldsymbol{S}_n:=\prod_{0\leq p\leq n}
 S_p$ are denoted by bold letters such as $\boldsymbol{x_n}=(x_p)_{0\leq p\leq n}\in \boldsymbol{\Sa}_n$ and $\boldsymbol{\textsl{\bf x}_n}=(\textsl{x}_p)_{0\leq p\leq n}\in \boldsymbol{S}_n$.

\subsubsection{Description of the models}\label{descrip-ref}

In the further development of this article, it is implicitly assumed that the functions on
population models are symmetric (due to the symmetry properties of Feynman-Kac models, this does not result in a loss of generality).
We write $\Ma_n$ for the Markov transitions of the particle model $\chi_n:=\xi_n$ viewed now as a \it Markov chain \rm on $\Sa_n$, and introduce the potential functions $\Ga_n(\chi_n)=m(\chi_n)(G_n)$.
We let $(\Pi_n,\Gamma_n)$ be the Feynman-Kac measures on $\Sa_n$ defined for any $\Fa_n\in \Ba(\Sa_n)$ by $\Pi_n(\Fa_n):={\Gamma_n(\Fa_n)}/{\Gamma_n(1)}$,
with
\begin{equation}\label{def-Gamma-Pi}
\Gamma_n(\Fa_n)=\EE\left(\Fa_n(\chi_n)~\Za_n(\chi)\right)
\quad\mbox{\rm and}\quad \Za_n(\chi)=\prod_{0\leq p<n}\Ga_p(\chi_p).
\end{equation}
 Equivalently, the Feynman-Kac model on  $\Sa_n$ can be turned
into a Feyman-Kac model on $\Pa(S_n)$, by replacing in (\ref{def-Gamma-Pi}) the chain $\chi_n$
by the occupation measure $m(\chi_n)$.

Also notice that
the unbiasedness properties of $\boldsymbol{\gamma^{(N,1)}_n(1)}$ ensures that
$\Gamma_n(1)=\gamma_n(1)$.
Using (\ref{evolve-eq-intro-ref}) it is readily checked that
\begin{equation}\label{evol-Gamma-Pi}
\Gamma_{n+1}=\Gamma_{n}\Qa_{n+1}\quad\mbox{and}\quad
\Pi_{n+1}:=\Psi_{\Ga_{n}}(\Pi_{n})\Ma_{n+1}
\end{equation}
with the integral operators 
$$
\Qa_{n+1}(x_{n},dx_{n+1})=\Ga_{n}(x_{n})~\Ma_{n+1}(x_{n},dx_{n+1})
$$

We denote by $(\boldsymbol{\Pi_n,\Gamma_n})$ the Feynman-Kac measures 
 associated with the historical process $\boldsymbol{\chi_n}=\left(\chi_0,\ldots,\chi_n\right)
 $, and the potential functions $\boldsymbol{\Ga_n(\chi_n)}:=\Ga_n(\chi_n)$
  on the path space $\boldsymbol{\Sa}_n$.
These measures are defined for any $\boldsymbol{\Fa_n}\in\Ba(\boldsymbol{\Sa_n})$ by
$\boldsymbol{\Pi_n(\Fa_n):={\Gamma_n(\Fa_n)}/{\Gamma_n(1)}}$, with
\begin{equation}\label{FK-path-intro-MF}
\boldsymbol{\Gamma_n(\Fa_n)}=\EE\left(\boldsymbol{\Fa_n(\chi_n)}~{\Za_n(\chi)}\right).
\end{equation}

Whenever the integral operators $Q_{n}$ have some densities $H_n$ w.r.t. some reference distributions  $\upsilon_n$
 on $S_n$, given  $\boldsymbol{\chi_n}$ we let ${\bf \XXb_n^{\flat}}$ be a random path with conditional distributions
 \begin{eqnarray}\label{def-Kab}
\boldsymbol{\Ka_n^{\flat}}(\boldsymbol{\chi_n},\boldsymbol{d\textsl{\bf x}_n})&:=&m(\chi_n)(d\textsl{x}_n)~\prod_{1\leq k\leq n}\LL_{k,m\left(\chi_{k-1}\right)}(\textsl{x}_k,d\textsl{x}_{k-1})
\end{eqnarray}
In the above displayed formula $\boldsymbol{\textsl{\bf x}_n}$ stands for  the path 
$\boldsymbol{\textsl{\bf x}_n}=(\textsl{x}_p)_{0\leq p\leq n}\in  \boldsymbol{S_n}$, and $\LL_{k,m(x_{k-1})}$ are the Markov transitions defined in (\ref{def-LL}).
We also denote by $\XX_n$ a random variable with conditional distribution given $\chi_n$ defined by
 \begin{eqnarray}\label{def-Kab-2}
\Ka_n(\chi_n,d\textsl{x}_n)=m(\chi_n)(d\textsl{x}_n)
\end{eqnarray}

 The unbiasedness properties of the measures
$\boldsymbol{\gamma^{(N,i)}_n}$ are equivalent to the fact that
 for any $\boldsymbol{f_n}\in\Ba(\boldsymbol{S_n})$ and $f_n\in\Ba(S_n)$, we have
\begin{equation}
\EE\left(\boldsymbol{f_n(\XXb^{\flat}_n)}~{\Za_n(\chi)}\right)=
\EE\left(\boldsymbol{f_n(X_n)}~Z_n(X)\right)\quad\mbox{\rm and}\quad
\EE\left(f_n(\XX_n)~{\Za_n(\chi)}\right)=\EE\left(f_n(X_n)~Z_n(X)\right)\label{ident}
\end{equation}
We emphasize that the last assertion in (\ref{ident}) holds for general Feynman-Kac models (i.e. without any regularity on $Q_n$).

\begin{defi}
The measures $(\Pi_n,\Gamma_n)$ and their path space versions  $(\boldsymbol{\Pi_n,\Gamma_n})$ are called the many body Feynman-Kac measures associated with the particle interpretation
(\ref{mean-field-intro}) of the measures $(\eta_n,\gamma_n)$. 
\end{defi}
As the name ``many-body'' suggests, these Feynman-Kac models encode properly the collective motion under mean field constraints of the system of particles associated to a standard Feynman-Kac particle system.

From an abstract point of view, in view of (\ref{ident}), all of these  measures are of course essentially equivalent to the  abstract Feynman-Kac
model introduced in (\ref{FK-def-intro-ref}).

\section{Conditional particle Markov chain models}\label{conditiona-sec}

This section aims at understanding PMCMC samplers from the point of view of many body Feynman-Kac models.

\subsection{Transport equation for many body Feynman-Kac models}
\label{islands-pm}

We start the section with a
pivotal duality formula between the Feynman-Kac integral operators $(Q_n,\Qa_n)$.

\begin{lem}\label{lem-key-duality}
We have the duality formula between integral operators on $\Sa_n\times S_n$
\begin{equation}\label{key-duality}
\Qa_{n}(x_{n-1},dx_n)~m(x_{n})(d\textsl{x}_n)~=~(m(x_{n-1})Q_n)(d\textsl{x}_n)~\Ma_{\textsl{x}_n,n}(x_{n-1},dx_n)
\end{equation}
and
$$
\eta^{\otimes N}_0(dx_0)~m(x_0)(d\textsl{x}_0)=\eta_0(d\textsl{x}_0)~\mu_{\textsl{x}_0}(dx_0)
$$
with the collection of Markov transitions 
$$
\Ma_{\textsl{x}_n,n}(x_{n-1},dx_n)=\frac{1}{N}\left[\sum_{i=0}^{N-1}\Phi_{n}(m(x_{n-1})
)^{\otimes (i)}\otimes \delta_{\textsl{x}_{n}}\otimes\Phi_{n}(m(x_{n-1})
)^{\otimes (N-i-1)}\right](dx_n)
$$
and the distribution
$$
\mu_{\textsl{x}_0}:=\frac{1}{N}\sum_{i=0}^{N-1}\left(\eta_0
^{\otimes (i)}\otimes \delta_{\textsl{x}_{0}}\otimes\eta_0^{\otimes (N-i-1)}\right)
$$
\end{lem}
\proof
To prove (\ref{key-duality}) 
we use the symmetry properties of the Markov transitions $\Ma_n$ to check that
for any function $H_n\in\Ba(S_n\times\Sa_n)$ (extended by right composition with the canonical projection from $S_n^N$ to $\Sa_n$ to a function still written $H_n$ in $\Ba(S_n\times S_n^N)$), we have
$$
\begin{array}{l}
\int~\Qa_{n}(x_{n-1},dx_n)~m(x_{n})(d\textsl{z}_n)~H_n(\textsl{z}_n,x_n)\\
\\
=\Ga_{n-1}(x_{n-1})~\int~\Phi_{n}(m(x_{n-1}))^{\otimes N}(dx_n)~H_n(\textsl{x}^1_n,x_n)\\
\\
=
m(x_{n-1})(G_{n-1})~\int~\Phi_{n}(m(x_{n-1}))(d\textsl{x}^1_n)~\left[\delta_{\textsl{x}^1_n}\otimes
\Phi_{n}(m(x_{n-1}))^{\otimes (N-1)}\right](dy_n)~H_n(\textsl{x}^1_n,y_n)
\end{array}
$$
The end of the proof comes
from the fact that
$$
m(x_{n-1})(G_{n-1})~\Phi_{n}(m(x_{n-1}))(d\textsl{x}^1_n)=(m(x_{n-1})Q_n)(d\textsl{x}^1_n)
$$
The proof of the lemma is now completed.
\cqfd

\begin{defi}\label{def-Xa}
Given a random path $(X_n)_{n\geq 0}$ we let $\Xa_n=\{\Xa_{n}^i\}_{i=1...N}\in \Sa_n$ be the Markov chain with 
the transitions 
$
\Ma_{X_n,n}
$,
and the initial distribution
$
\mu_{\textsl{X}_0}
$ introduced in lemma~\ref{lem-key-duality}.
We denote by $\boldsymbol{\MM_{n}}(\boldsymbol{\textsl{\bf X}_n},\boldsymbol{dx_n})$
the conditional distributions of the random path $\boldsymbol{\Xa_{n}}=\left(\Xa_{p}\right)_{0\leq p\leq n}$ on $\boldsymbol{\Sa_n}$.
The process $\Xa_n$ is called the dual mean field model associated with the Feynman-Kac particle model $\chi_n$ and the 
frozen path $\boldsymbol{\textsl{\bf X}_n}$. 
\end{defi}

The justification of the "duality" terminology between the processes $\Xa_n$ and $\chi_n$ is discussed in the end of the section.
The Feynman-Kac measures $\left(\boldsymbol{\gamma_n},\boldsymbol{\eta_n}\right)$ 
and their many body version $\left(\boldsymbol{\Gamma_n},\boldsymbol{\Pi_n}\right)$ are connected
by the following duality theorem which can be seen as an extended version of the unbiasedness properties (\ref{ident}).

\begin{theo}\label{theo-transport}

 For any $\boldsymbol{F_n}\in\Ba(S_n\times\boldsymbol{\Sa_n})$ we have the duality formula
\begin{equation}\label{key-decomposition}
\EE\left(\boldsymbol{F_n}(\XX_n,\boldsymbol{\chi_n})~{\Za_n(\chi)}\right)=\EE\left(\boldsymbol{F_n}(X_n,
\boldsymbol{\Xa_{n}})~{Z_n(X)}\right)
\end{equation}

When the integral operators $Q_{n}$ have some densities $H_n$ w.r.t. some reference distributions  $\upsilon_n$,
 for any $\boldsymbol{F_n}\in\Ba(\boldsymbol{S_n\times\Sa_n})$ we also have the duality formulae
\begin{equation}\label{key-decomposition-I}
\EE\left(\boldsymbol{F_n}({\bf \XXb_n^{\flat}},\boldsymbol{\chi_n})~{\Za_n(\chi)}\right)=\EE\left(\boldsymbol{F_n}(\boldsymbol{X_n},
\boldsymbol{\Xa_n})~Z_n(X)\right)
\end{equation}

\end{theo}
\proof
The proof of (\ref{key-decomposition}) is a direct consequence of (\ref{key-duality}).
Indeed, using this formula, we find that 

\begin{eqnarray*}
\Qa_{n}(x_{n-1},dx_n)m(x_n)(d\textsl{z}_n)
&=&\left[m(x_{n-1})Q_n\right](d\textsl{z}_n)~\Ma_{\textsl{z}_n,n}(x_{n-1},dx_n)
\\
&=&\int~m(x_{n-1})(d\textsl{z}_{n-1})~Q_n(\textsl{z}_{n-1},d\textsl{z}_n)~\Ma_{\textsl{z}_n,n}(x_{n-1},dx_n)
\end{eqnarray*}
and therefore
$$
\begin{array}{l}
\Qa_{n-1}(x_{n-2},dx_{n-1})\Qa_{n}(x_{n-1},dx_n)m(x_n)(d\textsl{z}_n)\\
\\
\displaystyle=\int~m(x_{n-2})(d\textsl{z}_{n-2})~Q_{n-1}(\textsl{z}_{n-2},d\textsl{z}_{n-1})~Q_n(\textsl{z}_{n-1},d\textsl{z}_n)\\
\\
\displaystyle\hskip7cm\times \Ma_{\textsl{z}_{n-1},n-1}(x_{n-2},dx_{n-1})\Ma_{\textsl{z}_n,n}(x_{n-1},dx_n)
\end{array}
$$
Iterating backward in time we prove (\ref{key-decomposition}) we find that
$$
\begin{array}{l}
\boldsymbol{\Gamma_n}(\boldsymbol{dx_n})~m(x_n)(d\textsl{z}_n)=
\int~\eta_0(d\textsl{z}_0) Q_1(\textsl{z}_{0},d\textsl{z}_n)\ldots Q_n(\textsl{z}_{n-1},d\textsl{z}_n)~\boldsymbol{\MM_{n}}(\boldsymbol{\textsl{\bf z}_n},\boldsymbol{dx_n})\end{array}
$$
This ends the proof of the first assertion.

The proof of (\ref{key-decomposition-I}) is a also direct consequence of (\ref{key-duality}).
Indeed, using this formula, we find that
$$
\begin{array}{l}
\boldsymbol{\Gamma_n}(\boldsymbol{dx_n})~\prod_{0\leq p\leq n}m(x_p)(d\textsl{x}_p)\\
\\
=
\Za_n(x)~
\eta_0^{\otimes N}(dx_0)~m(x_0)(d\textsl{x}_0)~\left\{\prod_{1\leq p\leq n}\Ma_p(x_{p-1},dx_p)~m(x_p)(d\textsl{x}_p)\right\}\\
\\
=\Za_n(x)~\left\{\eta_0(d\textsl{x}_0)~\prod_{1\leq p\leq n}\Phi_p(m(x_{p-1}))(d\textsl{x}_p)\right\}
\boldsymbol{\MM_{n}}(\boldsymbol{\textsl{\bf x}_n},\boldsymbol{dx_n})\\
\\
=\left\{\eta_0(d\textsl{x}_0)~\prod_{1\leq p\leq n}m(x_{p-1})(H_p(\point,\textsl{x}_p))~\upsilon_p(d\textsl{x}_p)\right\}
\boldsymbol{\MM_{n}}(\boldsymbol{\textsl{\bf x}_n},\boldsymbol{dx_n})
\end{array}
$$
with $\Za_n(x):=\prod_{0\leq p<n}m(x_p)(G_p)$.
The last assertion comes from the fact that
$$
m(x_{p-1})(G_{p-1})~\Phi_p(m(x_{p-1}))(d\textsl{\bf x}_p)=m(x_{p-1})(H_p(\point,\textsl{z}_p))~\upsilon_p(d\textsl{\bf x}_p)
$$
On the other hand, we have
we have
$$
\boldsymbol{\Ka_n^{\flat}}(\boldsymbol{x_n},\boldsymbol{d\textsl{\bf x}_n}):=m(x_n)(d\textsl{x}_n)~\prod_{1\leq p\leq n}\frac{m(x_{p-1})(d\textsl{x}_{p-1})~H_p(\textsl{x}_{p-1}, \textsl{x}_p)}{m(x_{p-1})(H_p(\point,\textsl{x}_p))}
$$
where $\boldsymbol{\textsl{\bf x}_n}$ stands for the path 
$\boldsymbol{\textsl{\bf x}_n}=(\textsl{x}_p)_{0\leq p\leq n}\in  \boldsymbol{S_n}$. Recalling that
$$
Q_p(\textsl{x}_{p-1}, d\textsl{x}_p)=G_p(\textsl{x}_{p-1})~M_p(\textsl{x}_{p-1}, d\textsl{x}_p)=H_p(\textsl{x}_{p-1}, \textsl{x}_p)~\upsilon_p(d\textsl{x}_p)
$$
This implies that
$$
\begin{array}{l}
\boldsymbol{\Gamma_n}(\boldsymbol{dx_n})~\boldsymbol{\Ka_n^{\flat}}(\boldsymbol{x_n},\boldsymbol{d\textsl{\bf x}_n})\\
\\
=\left\{\eta_0(d\textsl{x}_0)~\prod_{1\leq p\leq n}Q_p(\textsl{x}_{p-1}, d\textsl{x}_p)\right\}
\boldsymbol{\MM_{n}}(\boldsymbol{\textsl{\bf x}_n},\boldsymbol{dx_n})=\boldsymbol{\gamma_n}(\boldsymbol{d\textsl{\bf x}_n})~\boldsymbol{\MM_{n}}(\boldsymbol{\textsl{\bf x}_n},\boldsymbol{dx_n})
\end{array}
$$
 The proof of (\ref{key-decomposition-I}) is now completed.  This ends the proof of the Theorem.
\cqfd

The following Corollary is a direct consequence of (\ref{ident}) and (\ref{key-decomposition-I}). It provides an interpretation
of the conditional distribution of the dual process $\boldsymbol{\Xa_n}$ w.r.t. a given frozen trajectory as a conditional 
many body Feynman-Kac model w.r.t. a random path ${\bf \XXb^{\flat}_n}$ sampled with the backward distribution (\ref{def-Kab}).
\begin{cor}\label{cor-pGibbs}
For any $\boldsymbol{F_n}\in\Ba(\boldsymbol{\Sa_n})$, and for $(\boldsymbol{\eta_n}\otimes\eta_n)$-almost every paths $(\boldsymbol{\textsl{\bf x}_n},\textsl{x}_n)$
we have
\begin{eqnarray}
\EE\left({\bf F_n}(\boldsymbol{\Xa_n})~|~X_n={\textsl{x}_n}\right)&\propto&\EE\left(\boldsymbol{F_n}(\boldsymbol{\chi_n})~ \Za_n(\chi)~|~\XX_n=\textsl{x}_n\right)\label{cond-formula-direct}
\end{eqnarray}
For any $\left(f_n,\boldsymbol{f_n}\right)\in\left(\Ba(S_n)\times\Ba(\boldsymbol{S_n})\right)$, and for  
$\boldsymbol{\Pi_n}$-almost every path $\boldsymbol{x_n}$ we have
\begin{eqnarray*}
\EE\left({f_n}({\XX_n})~\left\vert~\boldsymbol{\chi_n}=\boldsymbol{x_n}\right.\right)&\propto&\EE\left({f_n}({X_n})~{Z_n}(X)~|~\boldsymbol{\Xa_n}=\boldsymbol{x_n}\right)\end{eqnarray*}
In addition, if the integral operators $Q_{n}$ have some densities $H_n$ w.r.t. some reference distributions  $\upsilon_n$, we have
\begin{eqnarray}
\EE\left({\bf F_n}(\boldsymbol{\Xa_n})~|~\boldsymbol{X_n}=\boldsymbol{\textsl{\bf x}_n}\right)&\propto&\EE\left(\boldsymbol{F_n}(\boldsymbol{\chi_n})~ \Za_n(\chi)~|~{\bf \XXb_n^{\flat}}=\boldsymbol{\textsl{\bf x}_n}\right)\label{cond-formula-backward}
\end{eqnarray}
\begin{eqnarray*}
\EE\left({\bf f_n}(\boldsymbol{\XXb^{\flat}_n})~\left\vert~\boldsymbol{\chi_n}=\boldsymbol{x_n}\right.\right)&\propto&\EE\left({\bf f_n}(\boldsymbol{X_n})~{Z_n}(X)~|~\boldsymbol{\Xa_n}=\boldsymbol{x_n}\right)
\end{eqnarray*}
\end{cor}

We end this section with an analytic description of the duality formulae (\ref{key-decomposition}) and (\ref{key-decomposition-I}) in terms of the
conditional distributions $\boldsymbol{\MM_{n}}$ and $\boldsymbol{\Ka_n^{\flat}}$ introduced in definition~\ref{def-Xa} and in (\ref{def-Kab}).
Using (\ref{key-decomposition})  we have 
$$
\forall \boldsymbol{\textsl{\bf x}_n}\in \boldsymbol{S_n}\qquad
\boldsymbol{\MM_{n}}(
\boldsymbol{\textsl{\bf x}_n},\point)\ll \boldsymbol{\eta_n}\boldsymbol{\MM_{n}}=\boldsymbol{\Pi_{n}}
$$ 
Thus, we can define the dual operator $\boldsymbol{\MM^{\star}_{n,\boldsymbol{\eta_n}}}$ of $\boldsymbol{\MM_{n}}$  from  $\LL_1(\boldsymbol{\eta_n})$ into $\LL_1(\boldsymbol{\Pi_n})$ given for any $\boldsymbol{f_n}\in \LL_1(\boldsymbol{\eta_n})$
by
$$
\boldsymbol{\MM^{\star}_{n,\boldsymbol{\eta_n}}}(\boldsymbol{f_n})=
\frac{{\bf d}\left(\boldsymbol{\eta_{n,\boldsymbol{f_n}}}\boldsymbol{\MM_{n}}\right)}{{\bf d}\left(\boldsymbol{\eta_n}\boldsymbol{\MM_{n}}\right)}=
\frac{{\bf d}\left(\boldsymbol{\eta_{n,\boldsymbol{f_n}}}\boldsymbol{\MM_{n}}\right)}{{\bf d}\boldsymbol{\Pi_{n}}}
\quad\mbox{\rm
with}\quad
 \boldsymbol{\eta_{n,\boldsymbol{f_n}}}(\boldsymbol{dx_n}):=\boldsymbol{\eta_n}(\boldsymbol{dx_n})~\boldsymbol{f_n}(\boldsymbol{x_n})~
$$
In addition, for any conjugate integers $\frac{1}{p}+\frac{1}{q}=1$, with $1\leq p,q\leq \infty$, and any pair of functions 
$\left(\boldsymbol{f_n},\boldsymbol{F_n}\right)\in\left(\LL_p(\boldsymbol{\eta_n})\times\LL_q\left(\boldsymbol{\Pi_n}\right)\right)$
we have
\begin{equation}\label{duality-form-star}
\boldsymbol{\Pi_n}\left(~\boldsymbol{F_n}~\boldsymbol{\MM^{\star}_{n,\boldsymbol{\eta_n}}}(\boldsymbol{f_n})~\right)=\boldsymbol{\eta_n}
\left(~\boldsymbol{\MM_{n}}(\boldsymbol{F_n})~
\boldsymbol{f_n}~
\right)
\end{equation}
These constructions shows that formula (\ref{key-decomposition-I}) holds true for general models (i.e. even if the  integral operators $Q_{n}$ don't have a density)
where ${\bf \XXb_n^{\flat}}$ stands for a random path with conditional distribution $\boldsymbol{\MM^{\star}_{n,\boldsymbol{\eta_n}}}\left(\boldsymbol{\chi_n},\point\right)$
given the historical process $\boldsymbol{\chi_n}$.

For a more detailed discussion on dual Markov transitions we refer the reader to~\cite{dlm-2003,revuz}.
In the reverse angle,
we have $$
\forall \boldsymbol{ x_n}\in \boldsymbol{\Sa_n}\qquad
\boldsymbol{\Ka_{n}^{\flat}}(
\boldsymbol{x_n},\point)\ll \boldsymbol{\Pi_n\Ka_{n}^{\flat}}=\boldsymbol{\eta_{n}}
$$
Thus (\ref{key-decomposition-I}) also  implies that $\boldsymbol{\MM_{n}}$ coincides with the dual operator
$\boldsymbol{\Ka^{\flat \star}_{n,\Pi_n}}$ of $\boldsymbol{\Ka_n^{\flat}}$  from  $\LL_1(\boldsymbol{\Pi_n})$ into $\LL_1(\boldsymbol{\eta_n})$; that is, we have
that
$$
(\ref{key-decomposition-I})~\Longrightarrow~
\boldsymbol{\Pi_n\Ka^{\flat}_{n}}=\boldsymbol{\eta_{n}}~~\Longrightarrow~~
\boldsymbol{\eta_{n}}\left(\boldsymbol{f_n}~\boldsymbol{\Ka^{\flat \star}_{n,\Pi_n}(F_n)}\right)=\boldsymbol{\Pi_n}\left(\boldsymbol{F_n}~\boldsymbol{\Ka_n^{\flat}}(\boldsymbol{f_n})\right)
$$
with
\begin{equation}\label{duality-form-star-2}
\boldsymbol{\Ka^{\flat \star}_{n,\Pi_n}}(\boldsymbol{\textsl{\bf z}_n},\boldsymbol{dx_n})=\boldsymbol{\Pi_n(dx_n)}~\frac{\boldsymbol{d\Ka^{\flat}_{n}}(
\boldsymbol{x_n},\point)}{d\boldsymbol{\Pi_n}\boldsymbol{\Ka^{\flat}_{n}}}(\boldsymbol{\textsl{\bf z}_n})=\boldsymbol{\MM_{n}}\left(\boldsymbol{\textsl{\bf z}_n},\boldsymbol{dx_n}\right)
\end{equation}
These formulations underline the duality between the random paths $\boldsymbol{\Xa_n}$ and $\XXb^{\flat}_n$ under the Feynman-Kac measures $\boldsymbol{\eta_n}$
and their many-body version $\boldsymbol{\Pi_n}$.

\subsection{Historical processes}\label{historical-sec}

Let us suppose that $(\eta_n,\gamma_n,\chi_n,\Ga_n)$ is the historical version of an auxiliary Feynman-Kac model
$(\gamma_n^{\prime},\eta_n^{\prime},\chi_n^{\prime},\Ga_n^{\prime})$ associated with some potential functions $G^{\prime}_n$ and some Markov chain $X^{\prime}_n$ transitions $M^{\prime}_n$ on some
state spaces $S_n^{\prime}$.  We let $\Phi^{\prime}_n$ be one step semigroup defined as $\Phi_n$ by replacing in (\ref{evolve-eq-intro-ref}), $(G_n,M_n)$ by $(G^{\prime},M^{\prime}_n)$.  We also let $Z^{\prime}_n(X^{\prime})=Z_n(X)$ and $\Za_n^{\prime}(\chi^{\prime})=\Za_n(\chi)$ the corresponding Radon-Nikodym derivatives defined in terms of $(G^{\prime}_n,\Ga^{\prime}_n,\chi^{\prime}_n,X^{\prime}_n)$.

When the integral operators $Q^{\prime}_{n}$ have some densities $H^{\prime}_n$ w.r.t. some reference distributions  $\upsilon^{\prime}_n$
 on $S^{\prime}_n$, the measure $\eta_n$ is expressed by a backward formula (\ref{backward-2}) with Markov transitions $\LL^{\prime}_{n+1,\eta_n^{\prime}}$
defined as in (\ref{def-LL})
 by replacing $(\eta_n,H_n)$ by $(\eta^{\prime}_n,H^{\prime}_n)$.

In this context, the Feynman-Kac models $(\eta_n,\gamma_n)$ and $(\Pi_n,\Gamma_n)$ defined in (\ref{FK-def-intro-ref}) and (\ref{def-Gamma-Pi}) are defined in terms of the  historical process $
X_n=\boldsymbol{X^{\prime}_n}=(X^{\prime}_0,\ldots,X_n^{\prime})$
of the chain $X^{\prime}_n$ and the ancestral lines $\chi_n$  of the particle model $\chi_n^{\prime}$. 

Notice that the dual process $\Xa_n$ associated with the particle model $\chi_n$ and the frozen path $\boldsymbol{X}_n=(X_0,\ldots,X_n)$
is a Markov chain on path space.  At each time step, given $\Xa_{n-1}$ we sample $N$ random trajectories  $\Xa_n=\left(\Xa^i_n\right)_{1\leq i\leq N}$.
One of them, say the first $\Xa^1_n=\textsl{x}_n$ takes the value of the frozen trajectory $X_n=\textsl{x}_n=(\textsl{x}^{\prime}_0,\ldots,\textsl{x}^{\prime}_n)\in S_n$.
The other ones are (conditionally) independent random paths with common distribution $\Phi_n\left(m(\Xa_{n-1})\right)$. This path-space chain can be interpreted
as an evolution of a genealogical tree with a frozen ancestral line.

We let $\MMb^{\natural}_n$ be the conditional expectation operator 
of the dual ancestral lines $\Xa_n$ given the frozen path $X_n$, that is
\begin{equation}\label{def-MMhist}
\EE\left(F_n(\Xa_n)~|~X_n\right):=\MMb^{\natural}_n(F_n)(X_n)
\end{equation}
for any function $F_n\in\Ba(\Sa_n)$. We also denote by $\left(\MMb_n,\Ka_n^{\flat}\right)$ the Markov transitions defined as $\left(\boldsymbol{\MMb_n,\Ka_n^{\flat}}\right)$ by replacing $(\boldsymbol{X_n,\Xa_n,\chi_n})$ by the historical processes $(\boldsymbol{X^{\prime}_n},\boldsymbol{\Xa_n^{\prime},\chi_n^{\prime}})$ of the chains
$(X^{\prime}_n,\Xa_n^{\prime},\chi^{\prime}_n)$. We let
$\XX_n^{\flat}$ be a random path on $S_n$ with conditional distributions $\Ka_n^{\flat}(\boldsymbol{\chi^{\prime}_n},d\textsl{x}_n)$
w.r.t. a complete ancestral tree $\boldsymbol{\chi^{\prime}_n}\in \boldsymbol{\Sa^{\prime}_n}=\prod_{0\leq k\leq n}\Sa^{\prime}_k$,
with   the product spaces $\Sa^{\prime}_k:=S^{\prime [N]}_k$, for any $k\geq 0$.

In this context, the duality formulae stated in Theorem~\ref{theo-transport} take the following form.

\begin{cor}
for any $\boldsymbol{F_n}\in\Ba(S_n\times\boldsymbol{\Sa_n})$, we have
\begin{equation}\label{key-decomposition-consequence}
\EE\left(\boldsymbol{F_n}(\XX_n,\boldsymbol{\chi_n})~{\Za_n^{\prime}(\chi^{\prime})}\right)=\EE\left(\boldsymbol{F_n}(X_n,
\boldsymbol{\Xa_{n}})~{Z_n^{\prime}(X^{\prime})}\right)
\end{equation}
In addition, when the integral operators $Q^{\prime}_{n}$ have some densities $H^{\prime}_n$ w.r.t. some reference distributions  $\upsilon^{\prime}_n$
 on $S^{\prime}_n$, we have
\begin{equation}\label{final-eq-hist}
\begin{array}{l}
\displaystyle\EE\left(F_n(\XX_n^{\flat},\boldsymbol{\chi^{\prime}_n})~\Za^{\prime}_n(\chi^{\prime})\right)
=
\displaystyle\EE\left(F_n(X_n,\boldsymbol{\Xa^{\prime}_n})~Z^{\prime}_n(X^{\prime})\right)=\EE\left(F_n(\XX_n,\boldsymbol{\chi^{\prime}_n})~\Za^{\prime}_n(\chi^{\prime})\right)
\end{array}
\end{equation}

\end{cor}

The formulae (\ref{key-decomposition-consequence}) and (\ref{final-eq-hist}) are direct consequences of (\ref{key-decomposition}) and (\ref{key-decomposition-I}).

In the further development of this section, we assume that  the integral operators $Q^{\prime}_{n}$ have some densities $H^{\prime}_n$ w.r.t. some reference distributions  $\upsilon^{\prime}_n$
 on $S^{\prime}_n$.

 Using (\ref{final-eq-hist}) we readily check that
\begin{prop} 
\begin{equation}\label{key-decomposition-consequence-2}
\mbox{\rm Law}\left(\XX_n^{\flat},\boldsymbol{\chi^{\prime}_n}\right)=\mbox{\rm Law}\left(\XX_n,\boldsymbol{\chi^{\prime}_n}\right)
\quad\mbox{\rm and}\quad
\mbox{\rm Law}\left(\XX_n|\boldsymbol{\chi^{\prime}_{n}}\right)=
\mbox{\rm Law}\left(\XX^{\flat}_n|\boldsymbol{\chi^{\prime}_{n}}\right)
\end{equation}
In particular, given the complete ancestral tree $\boldsymbol{\chi^{\prime}_{n}}$, the ancestral lines $\chi_n^i=(\chi^i_{k,n})_{0\leq k\leq n}$ are $N$ copies of the backward ancestral line
$\XX^{\flat}_n$ starting at $\chi^{\prime i}_n=\chi^i_{n,n}$, with $1\leq i\leq N$. 
\end{prop}

We emphasize that the path space model $\Xa_n$
is not equivalent to the genealogical tree evolution of the 
dual process $\Xa_n^{\prime}$ associated with the particle model $\chi^{\prime}_n$ and the frozen path $\boldsymbol{X^{\prime}_n}=(X^{\prime}_0,\ldots,X^{\prime}_n)=(\textsl{x}^{\prime}_0,\ldots,\textsl{x}^{\prime}_n)$. To be more precise, we observe that the Markov transitions of the ancestral lines $\Xa^{\#}_n=(\Xa^{\#, i}_n)_{1\leq i\leq N}$ of the dual process $\Xa_n^{\prime}$ are given by
$$
\Ma^{\#}_{\textsl{x}^{\prime}_n,n}(x_{n-1},dx_n)=\frac{1}{N}\left[\sum_{i=0}^{N-1}\Phi_{n}(m(x_{n-1})
)^{\otimes (i)}\otimes \delta_{(x^i_{n-1},\textsl{x}^{\prime}_n)}\otimes\Phi_{n}(m(x_{n-1})
)^{\otimes (N-i-1)}\right](dx_n)
$$
The ancestral line $\Xa^{\#, i}_n=(\Xa^{\#, i}_{k,n})_{0\leq k\leq n}$ 
encapsulates the genealogy of the $i$-th individual $\Xa_n^{\prime, i}$, including the frozen states at each level $k\leq n$.
By construction, the random states $\Xa_n^{\prime}$ are the $n$-th time marginal of random paths $\Xa_n$ and $\Xa^{\#}_n$.

The next result provides a more detailed description of the conditional distribution of the genealogical trees given the complete ancestral trees.
\begin{theo}\label{theo-ancestral-backward}
 Given the complete ancestral tree $\boldsymbol{\chi^{\prime}_{n}}$, the sequence of genealogical trees $(\chi_k)_{0\leq k\leq n}$ is a
Markov chain starting at $\chi_0=\chi^{\prime}_0$. The elementary transitions of the ancestral lines
$\chi_k\leadsto \chi_{k+1}$ given the population $\chi^{\prime}_{k+1}$ are defined  for any $f\in \Ba(\Sa_{k+1})$ 
$$
\begin{array}{l}
\displaystyle\EE\left(f(\chi_{k+1})~|~\chi_k,~\chi^{\prime}_{k+1}\right)\\
\propto
\displaystyle\int~\left\{\prod_{1\leq i\leq N}m(\chi_{k})(d\textsl{x}^i_{k})~H^{\prime}_{k+1}(\textsl{x}^i_{k,k},\chi^{\prime i}_{k+1})\right\}~f\left(
(\textsl{x}^l_{k},\chi^{\prime l}_{k+1})_{1\leq l\leq N}\right)
 \end{array}
 $$
 where $\textsl{x}^i_{k}:=\left(\textsl{x}^i_{l,k}\right)_{0\leq l\leq k}$ stands for an ancestral line of length $k$.

\end{theo}

\proof
By construction, for any $f_1,f_2\in\Ba(\Sa_{k+1})$  we have
$$
\begin{array}{l}
\EE\left(f_1(\chi_{k+1})f_2(\chi_{k},\chi^{\prime}_{k+1})~|~\chi_{k}\right)\\
\\
\propto
\int~\left\{\prod_{1\leq i\leq N}
m(\chi_{k})(d\textsl{x}^i_{k})~Q^{\prime}_{k+1}(\textsl{x}^i_{k,k},d\textsl{x}^{\prime i}_{k+1})\right\}~\\
\\\hskip7cm\times~
f_1((\textsl{x}^j_{k},\textsl{x}^{\prime j}_{k+1})_{1\leq j\leq N})~f_2\left(\chi_k,(\textsl{x}^{\prime j}_{k+1})_{1\leq j\leq N}\right)
\end{array}
$$
Using the fact that
$$
\begin{array}{l}
m(\chi_{k})(d\textsl{x}^i_{k})~Q^{\prime}_{k+1}(\textsl{x}^i_{k,k},d\textsl{x}^{\prime i}_{k+1})\\
\\
=
\displaystyle\frac{m(\chi_{k})(d\textsl{x}^i_{k}) H_{k+1}(\textsl{x}^i_{k,k},\textsl{x}^{\prime i}_{k+1})}{m(\chi^{\prime}_{k})(H^{\prime}_{k+1}(\point,\textsl{x}^{\prime i}_{k+1}))}\times m(\chi^{\prime}_{k})(H^{\prime}_{k+1}(\point,\textsl{x}^{\prime i}_{k+1}))~\nu_{k+1}(d\textsl{x}^{\prime i}_{k+1})\\
\\
\displaystyle\propto \frac{m(\chi_{k})(d\textsl{x}^i_{k}) H_{k+1}(\textsl{x}^i_{k,k},\textsl{x}^{\prime i}_{k+1})}{m(\chi^{\prime}_{k})(H^{\prime}_{k+1}(\point,\textsl{x}^{\prime i}_{k+1}))}\times \Phi^{\prime}_{k+1}\left(m(\chi^{\prime}_k)\right)(d\textsl{x}^{\prime i}_{k+1})
\end{array}
$$
we conclude that
$$
\begin{array}{l}
\EE\left(f_1(\chi_{k+1})f_2(\chi_{k},\chi^{\prime}_{k+1})~|~\chi_{k}\right)\\
\\
\displaystyle=\int~\left\{ \prod_{1\leq i\leq N} \Phi^{\prime}_{k+1}\left(m(\chi^{\prime}_k)\right)(d\textsl{x}^{\prime i}_{k+1})\right\}f_2\left(\chi_k,(\textsl{x}^{\prime j}_{k+1})_{1\leq j\leq N}\right)\\
\\
\hskip3cm
\displaystyle\int
\left\{\prod_{1\leq i\leq N}
\frac{m(\chi_{k})(d\textsl{x}^i_{k}) H_{k+1}(\textsl{x}^i_{k,k},\textsl{x}^{\prime i}_{k+1})}{m(\chi^{\prime}_{k})(H^{\prime}_{k+1}(\point,\textsl{x}^{\prime i}_{k+1}))}\right\}
f_1((\textsl{x}^j_{k},\textsl{x}^{\prime j}_{k+1})_{1\leq j\leq N})
\end{array}
$$

This ends the proof of the theorem.

\cqfd

\subsection{Genealogy and backward sampling models}\label{g+b-ref}

In the further development of this section, we assume that $(\eta_n,\gamma_n)$ is the historical version of an auxiliary Feynman-Kac model
$(\gamma_n^{\prime},\eta_n^{\prime})$.

\begin{defi}

We consider the Markov transitions from $S_n$ into itself defined by
 $\KK_n:=\MMb_{n}^{\natural}\Ka_n$,
with the operators $(\MMb_{n}^{\natural},\Ka_n)$ introduced in (\ref{def-MMhist})  and in (\ref{def-Kab-2}). In other words, for any function $f_n$ on $S_n$ and any frozen trajectory $\textsl{x}_n\in S_n$, we have
\begin{eqnarray*}
\KK_n(f_n)(\textsl{x}_n)
 &=&\EE\left(m(\Xa_n)(f_n)~|~X_n=\textsl{x}_n\right)
\end{eqnarray*}

When the integral operators $Q^{\prime}_{n}$ have some densities $H^{\prime}_n$ w.r.t. some distributions  $\upsilon^{\prime}_n$,
we consider the couple of Markov transitions from $S_n$ into itself defined by
 $ \KK_n^{\flat}:=\MMb_{n}\Ka_n^{\flat}
 $,
with the operators $(\MMb_{n},\Ka_n^{\flat})$ introduced in Section~\ref{historical-sec}.

\end{defi}

\begin{prop}\label{prop-equivalence}
The Markov transitions $\KK_n$, and $\KK_n^{\flat}$ are reversible w.r.t. the probability measures $\eta_n$. 
\end{prop}
The reversibility property is a direct consequence of the fact  that conditional SMC type PMCMC chains reduce to a standard
Gibbs sampler of a many-body Feynman-Kac distribution.

Next, under some rather strong regularity condition, we present an elementary proof of the ergodicity of the couple of conditional PMCMC transitions discussed above.
Sharp estimates of the contraction properties of $\KK_n$ and {\em its iterates $\KK^m_n$}, with $m\geq 1$, are developed in~section~\ref{taylor-invariant}.
These quantitative estimates are based on new 
Taylor type expansions of the PMCMC transitions  around the limiting invariant measure $\eta_n$ w.r.t. the precision parameter $1/N$.

\begin{prop}
We assume that the potential functions $G_n$ are lower and upper bounded by some positive constant, and we set $g_n:=\sup_{x,y}G_n(x)/G_n(y)$.
The measure   $\eta_n$ is the unique invariant measures of the Markov transitions $\KK_n$ and $\KK_n^{\flat}$.
In addition, we have
the estimates
\begin{equation}\label{dob-KK}
\beta\left(\KK_n\right)\vee \beta\left(\KK^{\flat}_n\right)\ \leq 1-\tau_n~\left(1-\frac{1}{N}\right)^{n}\quad
\mbox{with}\quad  \tau_n=1/\prod_{0\leq p<n}g_p.
\end{equation}
\end{prop}

The estimate of $\beta\left(\KK_n^{\flat}\right)$ is direct consequence of the following rather crude uniform minorization
condition
\begin{equation}\label{dob-KK-bis}
\KK_n(f_n)(\textsl{x}_n)\wedge \KK_n^{\flat}(f_n)(\textsl{x}_n)\geq \tau_n~\left(1-\frac{1}{N}\right)^{n}~\eta_n(f_n)
\end{equation}
 for any non-negative function
$f_n$ on $S_n$, and any path sequence $\textsl{x}_n=(\textsl{x}^{\prime}_p)_{0\leq p\leq n}$.
These lower bounds are easily checked by induction w.r.t. the time parameter. By construction, for any $\textsl{\bf x}_n=(\textsl{x}_{n-1},\textsl{x}^{\prime}_n)\in
S_n=\left(S_{n-1}\times S^{\prime}_n\right)$
we have
$$
\KK_n^{\flat}(f_n)(\textsl{x}_n)\geq  g_{n-1}^{-1}(1-1/N)~\KK_{n-1}^{\flat}(\overline{Q}_n(f_n))(\textsl{x}_{n-1})
$$
In much the same way, we prove that
\begin{eqnarray*}
\KK_n(f_n)(\textsl{x}_n)&\geq& g_{n-1}^{-1}(1-1/N) ~\EE\left(m(\Xa_{n-1})(\overline{Q}_n(f_n))~|~X_{n-1}=\textsl{x}_{n-1}\right) \\
&=&g_{n-1}^{-1}(1-1/N)~\KK_{n-1}(\overline{Q}_n(f_n))(\textsl{x}_{n-1})
\end{eqnarray*}
Iterating these estimates we check (\ref{dob-KK-bis}). 

\subsection{Taylor type expansions around the invariant measure}\label{taylor-invariant}

We assume in this paragraph that $(\eta_n,\gamma_n,\xi_n)$ is the historical version of an auxiliary Feynman-Kac model
$(\gamma_n^{\prime},\eta_n^{\prime},\xi_n^{\prime})$.
Our first objective is to find a Taylor type expansion of the Markov transition $\KK_n$ around its invariant measure $\eta_n$
w.r.t. powers of $1/N$.   We fix the time horizon $n$ and a frozen trajectory
 $\textsl{z}_n:=(\textsl{z}^{\prime}_0,\ldots,\textsl{z}^{\prime}_n)\in S_n=(S^{\prime}_0\times \ldots\times S^{\prime}_n)$, and for any $0\leq p\leq n$
 we set $\textsl{z}_p:=(\textsl{z}^{\prime}_0,\ldots,\textsl{z}^{\prime}_p)\in S_p$.

 We denote by $\Xa_{\textsl{z}_n,n}$ the dual
 mean field model associated with the Feynman-Kac particle model $\chi_n$ and the 
frozen path $X_n=\textsl{z}_n$. Using the exchangeability properties of the dual particles, there is no loss of generality
to assume that only the first one $\Xa_{\textsl{z}_n,n}^1=X_n$ is frozen.
With this convention, for any function $f_n\in\Ba(S_n)$
we have
$$
\KK_n(f_n)(\textsl{z}_n)=\EE\left(m(\Xa_{\textsl{z}_n,n})(f_n)\right)=
\frac{1}{N}~f_n(\textsl{z}_n)+\left(1-\frac{1}{N}\right)~\EE\left(m(\Xa_{\textsl{z}_n,n}^-)(f_n)\right)
$$
where  $m(\Xa_{\textsl{z}_n,n}^-)$ stands for the occupation measure of the non frozen particles
$
m(\Xa_{\textsl{z}_n,n}^-):=\frac{1}{N-1}~\sum_{1<i\leq N}\delta_{\Xa_{\textsl{z}_n,n}^i}.
$
This shows that whenever they exists these Taylor expansions are related to the bias and the fluctuations of the 
measures $m(\Xa_{\textsl{z}_n,n}^-)$. 
To analyze these properties we observe that
$$
\EE\left(m(\Xa_{\textsl{z}_n,n})(f_n)~|~\Xa_{\textsl{z}_{n-1},n-1}\right)=\Phi_{\textsl{z}_n,n}\left(m(\Xa_{\textsl{z}_{n-1},n-1})\right)(f_n)
$$
with the one step transformations $\Phi_{\boldsymbol{\textsl{\bf z}_n},n}$ defined as $\Phi_n$ by replacing the Markov transitions
$M_n$ by
$$
M_{\textsl{z}_n,n}(\textsl{x}_{n-1},d\textsl{x}_n)=\frac{1}{N}~\delta_{\textsl{z}_n}(d\textsl{x}_n)+\left(1-\frac{1}{N}\right)~M_n(\textsl{x}_{n-1},d\textsl{x}_n)
$$
In addition, the occupation measures $m(\Xa_{\textsl{z}_n,n}^-) $ of all the particles but the first frozen ones
are based on $(N-1)$ conditionally independent random states with common law $\Phi_{n}\left(m(\Xa_{\textsl{z}_{n-1},n-1})\right)$.
Thus, the local fluctuations of $m(\Xa_{\textsl{z}_n,n})$ around  
$\Phi_{\textsl{z}_n,n}\left(m(\Xa_{\textsl{z}_{n-1},n-1})\right)$ can be expressed in terms of the local sampling
random fields
$$
V^N_n:=\sqrt{N-1}~\left[m(\Xa_{\textsl{z}_n,n}^-)-\Phi_{n}\left(m(\Xa_{\textsl{z}_{n-1},n-1})\right)\right]
$$
with the formula
$$
m(\Xa_{\textsl{z}_n,n})(f_n)=\Phi_{\boldsymbol{\textsl{\bf z}_n},n}\left(m(\Xa_{\textsl{z}_{n-1},n-1})\right)+\left(1-\frac{1}{N}\right)~\frac{1}{\sqrt{N-1}}~V^N_n
$$

\begin{prop}
 Let $X_{\textsl{z}_n,n}$ stand for a Markov chain on $S_n$, with initial distribution 
$\eta_{\textsl{z}_0,0}=\frac{1}{N}\delta_{\textsl{z}_0}+\left(1-\frac{1}{N}\right)~\eta_0$ 
and Markov transitions $M_{\textsl{z}_n,n}$ from $S_{n-1}$ into $S_n$. We have
\begin{equation}\label{def-FK-order-1}
\EE\left(
m(\Xa_{\textsl{z}_n,n})(f_n)~\prod_{0\leq p<n}m(\Xa_{\textsl{z}_p,p})(G_p)
\right)=\EE\left(f_n(X_{\textsl{z}_n,n})~\prod_{0\leq p<n}G_p(X_{\textsl{z}_p,p})\right)
\end{equation}
\end{prop}

The proof is similar to the one that $\boldsymbol{\gamma}_n^{(N,1)}$ is an unbiased approximation of $\boldsymbol\gamma_n$ and omitted, see \cite{d-2004}.

The r.h.s. Feynman-Kac measure  in (\ref{def-FK-order-1}) can be expressed in terms of powers of the precision parameter $1/N$.
To describe these models, we let $\epsilon_n$ be a sequence of independent $\{0,1\}$-valued random
variables with $\PP(\epsilon_n=1)=1/N$. For any $\epsilon=(\epsilon_p)_{0\leq p\leq n}\in\{0,1\}^{n+1}$ we set 
$X_{\textsl{z}_n,n}^{(\epsilon)}$ be a Markov chain on $S_n$, with initial distribution $\eta^{(\epsilon)}_{\textsl{z}_0,0}$
and Markov transitions $M^{(\epsilon)}_{\textsl{z}_n,n}$ defined by 
\begin{eqnarray*}
\eta^{(\epsilon)}_{\textsl{z}_0,0}&=&\epsilon_0~\delta_{\textsl{z}_0}+\left(1-\epsilon\right)~\eta_0\\
M^{(\epsilon)}_{\textsl{z}_n,n}(\textsl{x}_{n-1},d\textsl{x}_n)&=&\epsilon_n~\delta_{\textsl{z}_n}(d\textsl{x}_n)+\left(1-\epsilon_n\right)~M_n(\textsl{x}_{n-1},d\textsl{x}_n)
\end{eqnarray*}
In this notation, we readily check that
$$
\begin{array}{l}
\EE\left(f_n(X_{\textsl{z}_n,n})~\prod_{0\leq p<n}G_p(X_{\textsl{z}_p,p})\right)
=\left(1-\frac{1}{N}\right)^{(n+1)}\displaystyle\gamma_n(f_n)\\
\\+\sum_{1\leq p\leq n+1}\left(\frac{1}{N}\right)^{p}
\left(1-\frac{1}{N}\right)^{(n+1)-p}~\sum_{\epsilon_0+\ldots+\epsilon_n=p
}  \EE\left[f_n(X_{\textsl{z}_n,n}^{(\epsilon)})~\prod_{0\leq p<n}G_p(X_{\textsl{z}_p,p}^{(\epsilon)})\right]
\end{array}
$$
These decompositions can be easily turned into Taylor's type polynomial expansions in power of $1/N$. The Taylor expansion of 
the normalized Feynman-Kac measures with the $0$-th order measure $\eta_n$
follows standard arguments on quotient power series.

The next proposition is easily proved using rather standard stochastic perturbation techniques (cf. for instance~\cite{d-2004,d-2013}).
\begin{prop}
The random fields $\sqrt{N}[m(\Xa_{\textsl{z}_n,n})-\eta_n]$ and $\sqrt{N}[m(\xi_n)-\eta_n]$
converge in law as $N\uparrow\infty$ to the same Gaussian and centered random fields. The same property holds true for the
random fields associated with the  unnormalized particle measures.
In addition, for any function $f_n\in\Ba(S_n)$ s.t. $\eta_n(f_n)=0$, and any frozen trajectory 
$\textsl{z}_n=(\textsl{z}_p^{\prime})_{0\leq p\leq n}\in S_n=\prod_{0\leq p\leq n}S^{\prime}_p$ we have the asymptotic bias expansion
\begin{equation}\label{first-order}
\begin{array}{l}
\lim_{N\uparrow\infty}N~\KK_n(f_n)(\textsl{z}_n)
=\sum_{0\leq p\leq n}\eta_p\left(\overline{Q}_{p,n}(1)~\left[\overline{Q}_{p,n}\left(f_n\right)(\textsl{z}_{p})-\overline{Q}_{p,n}(f_n)\right]\right)
\end{array}
\end{equation}
with $\textsl{z}_p:=(z^{\prime}_0,\ldots,z^{\prime}_p)\in S_p$, for any $p\leq n$.
\end{prop}

To get one step further, we need to analyze the propagation properties of the non frozen particles. 

\begin{theo}\label{theo-taylor-intro}
 {We assume that the condition ($H^\prime$) stated in (\ref{Hprime}) is satisfied. There exists some finite constant $c>0$ such that for any $n\geq 0$, $m\geq 1$  and 
$N>cn(m \overline{g}_n)^2$ we have 
 \begin{equation}\label{Taylor-KK1}
 \KK_n(\textsl{z}_n,d\textsl{y}_n)=\eta_n(d\textsl{y}_n)+\sum_{1\leq k\leq m}\frac{1}{N^k}~d^{(k)}\KK_n(\textsl{z}_n,d\textsl{y}_n)+\frac{1}{N^{m+1}}~\partial^{(m+1)} \KK_n(\textsl{z}_n,d\textsl{y}_n)
 \end{equation}
In the above displayed formula $d^{(k)}\KK_n$, and $\partial^{(m+1)} \KK_n$ stands 
 for  some sequence of signed and bounded integral operators  such that
 \begin{equation}\label{key-properties}
 d^{(k)}\KK_n(1)(\textsl{z}_n)=\partial^{(m+1)} \KK_n(1)(\textsl{z}_n)=0= \eta_n\left(d^{(k)}\KK_n(f_n)\right)=\eta_n\left(\partial^{(m+1)} \KK_n(f_n)\right)
 \end{equation}
  for any function $f_n$ on the path space $S_n$, and
 \begin{equation}\label{key-properties-2}
 \beta\left(d^{(k)}\KK_n\right)\leq \left[cn(k \overline{g}_n)^2\right]^{k}~~\mbox{and}~~\beta\left(\partial^{(m+1)} \KK_n\right)\leq  \left({cn((m+1)
\overline{g}_{n})^2}\right)^{m+1}
 \end{equation}
 In addition, when the Feynman-Kac model $(\gamma_n^{\prime},\eta^{\prime}_n)$ satisfies the regularity condition {\rm (H)} stated in (\ref{H}), the above estimate remains valid by replacing $\overline{g}_n$ by $\overline{g}:=\sup_{n\geq 0}\overline{g}_n<\infty$}
\end{theo}
 {
This Theorem is a particular case of the more general Theorem 
\ref{final-theo}, that can basically be stated as follows.
We let 
\begin{equation}\label{nonfrozen}
\PP^{(N,q)}_{\textsl{z}_n,n}=\mbox{\rm Law}\left(\Xa^{2}_{\textsl{z}_n,n},\Xa^{3}_{\textsl{z}_n,n},\ldots,\Xa^{q+1}_{\textsl{z}_n,n}\right)
\end{equation} be the distribution of the first $q$ random non frozen particles $\Xa^{i+1}_{\textsl{z},n}$ $i=1,\ldots,q$.
In this notation, and under the regularity conditions stated in Theorem~\ref{theo-taylor-intro}, there exists some finite constant $c>0$ such that for any $n\geq 0$, $r>m \geq 1$  and 
$N>cn((r+q) \overline{g}_n)^2$ we have the Taylor expansion
\begin{equation}\label{nonfrozen-2}
\PP^{(N,q)}_{\textsl{z}_n,n}=\eta_n^{\otimes q}+\sum_{1\leq k\leq m}\frac{1}{N^k}~d^{(k)}\PP^{(q)}_{\textsl{z}_n,n}+\frac{1}{N^{m+1}}~\partial^{(m+1)}\PP^{(N,q)}_{\textsl{z}_n,n}
\end{equation}
for some signed and bounded measures $d^{(k)}\PP^{(q)}_{\textsl{z}_n,n}$ with null mass $d^{(k)}\PP^{(q)}_{\textsl{z}_n,n}(1)=0$ whose values don't depend on the population size $N$ and such that
$$
\left\Vert d^{(k)}\PP^{(q)}_{\textsl{z},n} \right\Vert_{\tiny tv}\leq \left[cn((q+2k) \overline{g}_n)^2\right]^{k}
\quad\mbox{\rm and}\quad
 \left\Vert\partial^{(m+1)}\PP^{(N,q)}_{\textsl{z}_n,n} \right\Vert_{\tiny tv}\leq  b(q)~\left(cn((q+m) \overline{g}_{n})^2\right)^{m+1}
$$
with some finite constant $b(q)<\infty$ whose values only depend on the parameters $q$.
A more precise description of the derivatives operators is provided in Theorem~\ref{final-theo}.}

 {
We end this section with some direct consequences of these expansions around the fixed point Feynman-Kac measures. To illustrate our result we assume that the 
 regularity condition ({\rm H})  is met and the size of the system $N$ is chosen s.t. $N>cn\overline{g}^2$ for some finite constant $c<\infty$.}
 
 {
$\bullet$ These expansions can also be used to estimate of 
the behavior of the particle measures $m(\xi_{\textsl{z}_n,n})$ as $N\uparrow\infty$. For instance, we have the bias and the variance estimates
 $$
 \EE\left(m(\Xa_{\textsl{z}_n,n})(f_n)\right)=\eta_n(f_n)+\frac{1}{N}~\left([f_n(\textsl{z}_n)-\eta_n(f)]+d^{(1)}\PP_{\textsl{z}_n,n}^{(1)}(f)\right)+\frac{1}{N^2}~r^{(N,1)}_{\textsl{z}_n,n}(f)
 $$
 and
 $$
 \mbox{\rm Var}\left(m(\Xa_{\textsl{z}_n,n})(f_n)\right)=\frac{1}{N}~\left([\eta_n(f^2_n)-\eta_n(f_n)^2]+d^{(1)}\PP_{\textsl{z}_n,n}^{(2)}((f-\eta_n(f))^{\otimes 2})\right)+\frac{1}{N^2}~r^{(N,2)}_{\textsl{z}_n,n}(f)
 $$
 with a couple of remainder terms such that
 $$
\sup_{i=1,2}\left\vert r^{(N,i)}_{\textsl{z}_n,n}(f)\right\vert\leq c~(n\overline{g}^2)^2 \quad\mbox{\rm for some finite}\quad c<\infty
 $$
 The last estimate is related to the variance of the particle measures $m(\Xa_{\textsl{z}_n,n})$ delivered by the PMCMC model. In much the same way, the variance of a function of the trajectory delivered by the PMCMC model is computed using the expansion of $ \EE\left(m(\Xa_{\textsl{z}_n,n})(f_n^2)\right)$.}

 {
$\bullet$  Using the first order expansion (\ref{Taylor-KK1}), for any $\mu_n,\nu_n\in\Pa(S_n)$ we readily check that
$$
N^2~\left\Vert(\mu_n-\nu_n)\KK_n-\frac{1}{N}~(\mu_n-\nu_n)d^{(1)}\KK_n\right\Vert_{\tiny tc}\leq c~\left(n\overline{g}^2\right)^2
$$
with some finite constant $c<\infty$, and the first order integral operator $d^{(1)}\KK_n$ defined in (\ref{first-order})  and given by
$$
d^{(1)}\KK_n(f_n)(\textsl{z}_n)
=\sum_{0\leq p\leq n}\eta_p\left(\overline{Q}_{p,n}(1)~\left[\overline{Q}_{p,n}\left(f_n\right)(\textsl{z}_{p})-\overline{Q}_{p,n}(f_n)\right]\right)
$$
This implies  that
\begin{equation}\label{dob-estimates}
 \left\vert\beta\left(\KK_n\right)-\frac{1}{N}~\beta\left(d^{(1)}\KK_n\right)\right\vert\leq c~\left(n\overline{g}^2/N\right)^2
\end{equation}
Using (\ref{first-order}), we also have the crude estimate
$$
\beta\left(d^{(1)}\KK_n\right)\leq 2\sum_{0\leq k\leq n}~
\left\Vert\overline{Q}_{k,n}(1)\right\Vert\leq 2 (n+1)~\overline{g}
$$
To check that r.h.s. linear estimate w.r.t. the time horizon $n$ is sharp, we choose unit potential functions $G_n=1$ and
a function $f_n(\textsl{z}_n)=\varphi(\textsl{z}^{\prime}_0)$ that only depend on the initial state of the path $\textsl{z}_n=(\textsl{z}^{\prime}_k)_{0\leq k\leq n}\in S_n=\prod_{0\leq k\leq n}S^{\prime}_n$. In this situation, we have
$$
d^{(1)}\KK_n(f_n)(\textsl{y}_n)-d^{(1)}\KK_n(f_n)(\textsl{z}_n)
=(n+1)\left(\varphi(\textsl{y}^{\prime}_0)-\varphi(\textsl{z}^{\prime}_0)\right)\Rightarrow \beta\left(d^{(1)}\KK_n\right)\geq (n+1)
$$   

These estimates (\ref{dob-estimates}) ensure that the Markov chain with transitions
$\KK_n$ converge exponentially fast to $\eta_n$ with a rate that can be made arbitrary large when the precision parameter and the size
 of the particle population model
$N\uparrow\infty$.\\
$\bullet$ Using the properties (\ref{key-properties}) we readily prove Taylor expansions of any $m$-th iterate $ \KK_n^m=\KK_n^{m-1}\KK_n$ of the PMCMC transition $\KK_n$.  For instance, for any $m\geq 1$, we have
\begin{equation}\label{Taylor-m}
\KK_n^m(\textsl{y}_n,d\textsl{z}_n)=\eta_n(d\textsl{z}_n)+\frac{1}{N^m}\left[d^{(1)}\KK_n\right]^m(\textsl{y}_n,d\textsl{z}_n)+\frac{1}{N^{m+1}}~\partial^{(m+1)}\KK^m(\textsl{y}_n,d\textsl{z}_n)
\end{equation}
with 
 the remainder integral operator $\partial^{(m+1)}\KK_n^m$ such that
$$
\partial^{(m+1)}\KK^m_n(1)(\textsl{y}_n)=0\quad\mbox{\rm and}\quad\beta\left(\partial^{(m+1)}\KK_n^m\right)\leq m~
\left(cn\overline{g}^2\right)^{m+1}~\left(1+{cn\overline{g}^2}/{N}\right)^{m-1}
$$
This result shows that
the distribution of the random state of the Markov chain with transition $\KK_n$ after $m$ iteration is 
equal to $\eta_n$ up to some remainder measure with total variation norm of order $N^{-m}$. In addition, arguing as above we find that
$$
N^{m+1}\left\vert \beta\left(\KK_n^m\right)-\frac{1}{N^m}~\beta\left(\left[d^{(1)}\KK_n\right]^m\right)\right\vert\leq \beta\left(\partial^{(m+1)}\KK_n^m\right)
$$
$\bullet$ The decompositions (\ref{Taylor-m}) can be used to derive without any additional work the $\LL_p$-norms 
between the distributions of the random states of the conditional PMCMC model and the invariant measures. For instance,
for any $p\geq 1$ we have
$$
\left\Vert \KK_n^m(f_n)-\eta_n(f_n)\right\Vert_{\LL_p(\eta_n)}=\frac{1}{N^m}~
\left\Vert  \left[d^{(1)}\KK_n\right]^m\!\!\!(f_n)\right\Vert_{\LL_p(\eta_n)}+\frac{1}{N^{m+1}}\left\Vert \partial^{(m+1)}\KK^m(f_n)\right\Vert_{\LL_p(\eta_n)}
$$
}
$\bullet$ The proof of the Taylor expansions (\ref{nonfrozen}) is based on renormalization techniques and a differential calculus
on the measures $\Upsilon^{(N,q)}_{\textsl{z}_n,n}$ on $S_n^q$ defined   for any $F_n\in\Ba(S^q_n)$ by
\begin{equation}\label{def-Upsilon}
\Upsilon^{(N,q)}_{\textsl{z}_n,n}(F_n):=\EE\left(m(\Xa_{\textsl{z}_n,n})^{\otimes q}(F_n)~\prod_{0\leq p<n}m(\Xa_{\textsl{z}_p,p})(\overline{G}_p)^q\right)
\end{equation}
We will show that $\Upsilon^{(N,q)}_{\textsl{z}_n,n}$ are differentiable at any order with $d^{(0)}\Upsilon^{(N,q)}_{\textsl{z}_n,n}=\eta_n^{\otimes q}$.
On the other hand, formula (\ref{key-decomposition}) implies that
\begin{equation}\label{key-transfert}
\int~\eta_n(d\textsl{z}_n)~\Upsilon^{(N,q-1)}_{\textsl{z}_n,n}(F_n)=\Upsilon^{(N,q)}_n(F_n\otimes 1)
\end{equation}
 for any $F_n\in\Ba(S^{q-1}_n)$ , with the measure $\Upsilon^{(N,q)}_n$ defined as $\Upsilon^{(N,q)}_{\textsl{z}_n,n}$ by replacing $(\Xa_{\textsl{z}_p,p})_{0\leq p\leq n}$
by $(\Xa_{n})_{0\leq p\leq n}$.
This formula can be used to compute Taylor type expansions for the occupation measures of the process $\Xa_n$,
including the $(q+1)$-moments of the unnormalized particle normalizing constants $\prod_{0\leq p<n}m(\chi_{p})(G_p)$.

 In this connexion,  the transfer formula (\ref{key-transfert}) also
 shows that the particle
approximation $\prod_{0\leq p<n}m(\Xa_p)(G_p)$ of the normalizing constants associated with the 
particle model with a frozen trajectory
is {\em biased even if the particle Markov chain model starts with the desired target measure}. For instance for $q=1$ and $F_n=1$
formula (\ref{key-transfert}) implies that
$$
\EE\left(\prod_{0\leq p<n}m(\Xa_p)(\overline{G}_p)\right)=1+\EE\left(\left[\prod_{0\leq p<n}m(\chi_{p})(\overline{G}_p)-1\right]^2\right)\not=1
$$
Running a Markov chain with one of the transitions $\KK_n$, we design
a asymptotically unbiased estimate using the easily checked formula
$$
\EE\left(\left[\prod_{0\leq p<n}m(\Xa_p)(G_p)\right]^{-1}\right)=\left[\prod_{0\leq p<n}\eta_p(G_p)\right]^{-1}
$$

\section{Propagation of chaos expansions}\label{general-second-order}

This section, as its name indicates, will focus on the fine analysis of the size $N$  dependency of PMCMC samplers and related problems such as asymptotic independency of $q<<N$ subsets of the particle models investigated in the first sections of the article --that is, propagation of chaos properties.

\subsection{Combinatorial preliminaries}\label{scp-sec}
We let $X=\left(X^i\right)_{2\leq i\leq N}$ be a sequence of random variables on some state space $S$, and $z\in S$ a given fixed state. For any $q<N$ we set
$$
m(X)^{\odot q}=\frac{1}{(N-1)_q}\sum_{a\in I_q^N}\delta_{\left(X^{a(1)},\ldots,X^{a(q)}\right)}
$$
where $I_q^N$ stands the set of of all $(N-1)_q=\frac{(N-1)!}{((N-1)-q)!}$ multi-indexes 
 $a=(a(1),\ldots,a(q))\in\{2,\ldots,N\}^q$ with different values, or equivalently one to one mappings from $[q]:=\{1,\ldots,q\}$ into
 $\{2,\ldots,N\}=[N]-\{1\}$. The link between these measures and tensor product measures 
 is expressed in terms of the  Markov transitions
 $\AA_{a}^{(q)}$ indexed by the set of  mappings $a$ from $[q]$ into itself and defined for any $x=(x^1,\ldots,x^q)\in S^q$
 by 
 $$
 \AA_{a}^{(q)}(F)(x)=F(x^a)\quad\mbox{\rm with}\quad
 x^a:=\left(x^{a(1)},\ldots,x^{a(q)}\right)
 $$
 for any function $F$ on $\Ba(S^q)$, and any $(x^1,\ldots,x^q)\in S^q$. The connection between these measures 
 is described in the following technical lemma taken from~\cite{dpr-2009}.
 
 We emphasize that the tensor product measures discussed above are symmetry-invariant by construction. In the further development of this section,
 it is assumed without restrictions that these measures act on symmetric functions $F$; that is
 $
 F=\frac{1}{q!}\sum_{\sigma\in\Ga_q} \AA_{\sigma}^{(q)}(F)
 $,
where $\Ga_q$ stands for the symmetric group of all permutations of $[q]$.
 
 \begin{lem}
 For any $q<N$ we have the formula
 $$
m(X)^{\otimes q}= m(X)^{\odot q}\AA^{(N,q)}
\quad\mbox{with}\quad
 \AA^{(N,q)}=\frac{1}{(N-1)^q}\sum_{a\in [q]^{[q]}}~\frac{(N-1)_{\vert a\vert}}{(q)_{\vert a\vert}}~\AA_a^{(q)}
 $$
where $\vert a\vert$ for the cardinality of the set $a([q])$, and $(m)_p=m!/(m-p)!$ stands for the number of one to one mappings
from $[p]$ into $[m]$.
\end{lem}

\begin{defi}
For any $\textsl{z}\in S$ we consider the random measures
$$
m_{\textsl{z}}(X)=\frac{1}{N}\delta_{\textsl{z}}+\left(1-\frac{1}{N}\right)~m(X)\qquad
m^{(1)}_{\textsl{z}}(X)=\delta_{z}\quad\mbox{and}\quad m^{(0)}_{\textsl{z}}(X)=m(X)
$$
For any $b\in \{0,1\}^{[q]}$, we denote by $\BB^{(q)}_{\textsl{z},b}$ the Markov transitions
defined for any $x=(x^1,\ldots,x^q)\in S^q$ by
$$
\BB^{(q)}_{\textsl{z},b}(F)(x)=F\left(x^{b}_{\textsl{z}}\right)\quad\mbox{
with}\quad
x^{b}_{\textsl{z}}:=\left(b(1)\textsl{z}+(1-b(1))x^1,\ldots,b(q)\textsl{z}+(1-b(q))x^q\right)
$$
\end{defi}

We observe that
$$
m_{\textsl{z}}(X)^{\otimes q}=\sum_{b\in \{0,1\}^{[q]}}
\frac{1}{N^{\vert b\vert_1}}~\left(1-\frac{1}{N}\right)^{q-\vert b\vert_1}~m^{(b)}_{\textsl{z}}(X)
$$
with $\vert b\vert_1=\sum_{1\leq p\leq q}b(p)$ and
$$
m^{(b)}_{\textsl{z}}(X)=m^{(b(1))}_{\textsl{z}}(X)\otimes\ldots\otimes m^{(b(q))}_{\textsl{z}}(X)
$$

\begin{lem}
For any $q<N$, and $b\in \{0,1\}^{[q]}$ we have $m^{(b)}_{\textsl{z}}(X)=m^{\otimes q}(X)\BB^{(N,q)}_{\textsl{z},b}$ and
$$
 m_{\textsl{z}}(X)^{\otimes q}=m^{\otimes q}(X)\BB^{(N,q)}_{\textsl{z}}
\quad\mbox{with}\quad
\BB^{(N,q)}_{\textsl{z}}=\sum_{b\in \{0,1\}^{[q]}}
\frac{1}{N^{\vert b\vert_1}}~\left(1-\frac{1}{N}\right)^{q-\vert b\vert_1}~\BB^{(q)}_{\textsl{z},b}
$$
as well as
$$
m_{\textsl{z}}(X)^{\otimes q}=m(X)^{\odot q}\CC^{(N,q)}_{\textsl{z}}\quad\mbox{with}\quad \CC^{(N,q)}_{\textsl{z}}:=\AA^{(N,q)}\BB^{(N,q)}_{\textsl{z}}
$$
\end{lem}

\begin{defi}
We let $(p_1,p_2)$ be a couple of integers s.t. $0\leq p_1\leq q-1$ and $0\leq p_2\leq q$.
\begin{itemize}
\item We consider the collection of sets
\begin{eqnarray*}
\Ia_{q}&:=&\{0,\ldots,q-1\}\times \{0,\ldots,q\}\qquad
~[r]_{q-p_1}^{[q]}:=\{a\in [r]^{[q]}~:~\vert a\vert=q-p_1 \}\\
\{0,1\}_{1,p_2}^{[q]}&:=&\{b\in \{0,1\}^{[q]}~:~\vert b\vert_1=p_2 \}
\quad\mbox{and}\quad I_q(p_1,p_2)=[q]_{q-p_1}^{[q]}\times \{0,1\}_{1,p_2}^{[q]}
\end{eqnarray*}
\item We let $\Aa^{(q)}_{p_1}$, and resp. $\Ba^{(q)}_{p_2}$ be the uniform distributions on $[q]_{q-p_1}^{[q]}$, and resp. on $\{0,1\}_{1,p_2}^{[q]}$. We also denote by 
$
\Ca^{(q)}_{p_1,p_2}=\Aa^{(q)}_{p_1}\otimes \Ba^{(q)}_{p_2}
$ the uniform measure on $I_q(p_1,p_2)$.

\item For any $c=(a,b)\in I_q(p_1,p_2)$, we let $\CC^{(q)}_{\textsl{z},(a,b)}$ be the coalescent operator  defined for any $x=(x^1,\ldots,x^q)\in S^q$ by
$$
\CC^{(q)}_{\textsl{z},(a,b)}(F)(x):=F(x_{\textsl{z}}^{(a,b)})
$$
with
$$ x_{\textsl{z}}^{(a,b)}=\left(
b(1)z+(1-b(1))x^{a(1)},\ldots,b(q)z+(1-b(q))x^{a(q)}\right),
$$
so that $\CC^{(q)}_{\textsl{z},(a,b)}=\AA^{(q)}_a\BB^{(q)}_{\textsl{z},b}$.

 \end{itemize}
\end{defi}

 \begin{rem}When maps in $[q]^{[q]}$ are represented graphically,
  the parameter $p_1$ in $[q]_{q-p_1}^{[q]}$represents the number of coalescences of the change of index mapping $a$. The
 $p_2$ is the number of $b(i)$ such that $b(i)=0$ or $x_{\textsl{z}}^{(a,b),i}=\textsl{z}$; it will be referred to as
 the number of $\textsl{z}$-infections of the mapping $b$.
\end{rem}

We recall that the Stirling numbers 
of the second kind $S(q,p)$ is the number of partitions of $[q]$ into $p$ sets, so that
$$
\#\left([r]_p^{[q]}\right)=S(q,p)~(r)_p\quad\mbox{\rm and}\quad r^q=\sum_{1\leq p\leq q}S(q,p)~(r)_p
$$ 
for any $p\leq q\leq r$. We also recall that the Stirling numbers of the first kind $s(q,p)$ provide the coefficients of the polynomial expansion
of $(r)_q$
\begin{equation}\label{SFK}
(r)_q=\sum_{1\leq p\leq q}s(q,p)~r^p
\end{equation}

We also use the conventions $(r)_{q}=0$ and $(r)_0=1=(-r)_0$ for any $q>r\geq 0$, as well as $s(q,0)=s(0,-q)=S(0,-q)=S(q,0)=0$
except $s(0,0)=S(0,0)=1$, for $q=0$.

These formulae can be found in any textbook on combinatorial analysis, see for instance~\cite{berge,comtet}.

\begin{defi}
 We also consider the sequence of probabilities $\Pa^{(N,q)}=\Pa^{[N,q,1]}\otimes \Pa^{[N,q,2]}$ on the set $\Ia_{q} $ defined by
\begin{equation}\label{statcoinf}
\Pa^{(N,q)}(p_1,p_2):=\underbrace{
\frac{1}{(N-1)^{q}}~S(q,q-p_1)~(N-1)_{q-p_1}
}_{\Pa^{[N,q,1]}(p_1)}\times
~\underbrace{\left(\begin{array}{c}
q\\p_2
\end{array}\right)~\left(1-\frac{1}{N}\right)^{q-p_2}\frac{1}{N^{p_2}}}_{\Pa^{[N,q,2]}(p_2)}
\end{equation}
\end{defi}

Notice that $\Pa^{[N,q,1]}(p_1)=\#\left([N-1]_{q-p_1}^{[q]}\right)/\#[N-1]^{[q]}$ is a statistics for the number of coalescences, whereas $\Pa^{[N,q,2]}(p_2)$ is the proportion of infested mappings with $p_2$ infections.
By construction, we have the following lemma.

\begin{lem}
For any $q<N$, we have the formula
\begin{eqnarray*}
\CC^{(N,q)}_{\textsl{z}}&=&\sum_{p\in \Ia_{0,q}}~\Pa^{(N,q)}(p)~~\widehat{\CC}_{\textsl{z},p}^{(q)}\quad\mbox{with}\quad
\widehat{\CC}_{\textsl{z},p}^{(q)}=\sum_{c\in I_q(p)}
\Ca^{(q)}_{p}(c)~
\CC^{(q)}_{\textsl{z},c}
\end{eqnarray*}

\end{lem}
We end this section with a Taylor expansion of the measure $\Pa^{(N,q)}$ introduced above.

\begin{prop}
For any $q<N$, the mapping $N\mapsto \Pa^{(N,q)}$ is differentiable at any order $m\geq 0$. The $m$-order derivative is supported by $$
  \Ta^{(m)}_{q,n}:=\{(p_1,p_2)\in \Ia_{q}~:~0\leq p_1+p_2\leq m\}
  .$$
\end{prop}

Indeed, Fla (\ref{statcoinf}) shows that the fraction in the variable $N$, $\Pa^{(N,q)}(p_1,p_2)$, can be expanded as a formal power series in $\frac{1}{N}$ (or, more precisely, as an analytic function in the neighborhood of $0$) with leading term in 
$\frac{1}{N^{p_1+p_2}}$.  The Proposition follows.

Expanding the formula for $\Pa^{(N,q)}(p_1,p_2)$ using (\ref{SFK}) and the Taylor expansion
$$\frac{1}{(1-x)^n}=\sum\limits_{0\leq k}(n+k-1)_k~\frac{x^k}{k!}=\sum\limits_{0\leq k}{n+k-1\choose k}~x^k$$
with $(n-1)_0:=1$, we get an explicit formula for the derivatives.

\begin{prop}\label{key-combinatorial-lemma}
The $m$-th order derivative is given by the signed measure (with total null mass) supported on the set
  $\Ta^{(m)}_{q,n}$:
\begin{equation}\label{key-formula-PN}
d^{(m)}\Pa^{(q)}:=\sum_{ (p_1,p_2)\in   \Ta^{(m)}_{q,n}}\tau^{(m)}_{q,p_1,p_2}~\delta_{(p_1,p_2)},
\end{equation}
with
$$
\tau^{(m)}_{q,p_1,p_2}=
\sum_{\boldsymbol{k}\in \Ka_q^{(m)}(p_1,p_2)}\alpha_{q,p_1,p_2}(\boldsymbol{k}),
$$

\begin{eqnarray}
\Ka_q^{(m)}(p_1,p_2)\hskip-.3cm&:=&\hskip-.3cm\left\{(k_1,k_2,k_3)\in [0,q-p_1[\times [0,q-p_2]\times\NN~:~\sum_{1\leq i\leq 2}p_i+\sum_{1\leq i\leq 3}k_i=m\right\},\nonumber\\
\alpha_{q,p_1,p_2}(k_1,k_2,k_3)\hskip-.4cm&=&\hskip-.3cm
S(q,q-p_1)\left(\begin{array}{c}
q\\p_2
\end{array}\right)~\nonumber\\
&&\hskip-.2cm\times
~s(q-p_1,q-p_1-k_1)~
(-1)^{k_2}~
\left(\begin{array}{c}
q-p_2\\k_2
\end{array}\right)
{(p_1+k_1)+k_3-1\choose k_3}\label{ref-bound-tau}
\end{eqnarray}

\end{prop}

\begin{rem}\label{key-rem-v2}
We observe that $\tau^{(0)}_{q,p_1,p_2}=1_{(0,0)}(p_1,p_2)$. As will appear later on, this identity encodes the propagation of chaos properties (i.e. asymptotic independency) of PMCMC samplers.
We also mention that $\alpha_{q,p_1,p_2}(\boldsymbol{k})=0=\tau^{(m)}_{q,p_1,p_2}$ as soon as $p_1>q$ or $p_2>q$.
\end{rem}

\begin{rem}\label{key-rem-v2followup}
The $m$-th order signed measure $d^{(m)}\Pa^{(q)}$ and the mapping $(p_1,p_2)\mapsto \tau^{(m)}_{q,p_1,p_2}$ in formula (\ref{key-formula-PN}) only charge
couple of integers $(p_1,p_2)\in \left([1,q]\times [0,q]\right)$ s.t. $0\leq p_1+p_2\leq m$. 
The first coordinate $0\leq p_1<q$ can be interpreted as the number of coalescent states, while
 $p_2$ can be interpreted as the the number of $z$-infected states.
 
By construction, the mapping $(p_1,p_2)\mapsto \tau^{(m)}_{q,p_1,p_2}$ can also be seen as a measure with null total mass supported on the  set $0\leq p_1+p_2\leq m$. For instance, for $m=1,2$, recalling that $s(q,q-1)=-q(q-1)/2=-S(q,q-1)$, $s(q,q-2)=\frac{q!}{3!(q-3)!}~\frac{3q-1}{4}$, and
$S(q,q-2)=\frac{q!}{3!(q-3)!}~(3q-5)/4$, we have
\begin{equation}\label{ref-beta-m}
\begin{array}{rclcrcll}
\tau^{(2)}_{q,2,0}&=&\frac{q!}{3!(q-3)!}~\frac{3q-5}{4}                              &\tau^{(2)}_{q,0,2}&=&\frac{q(q-1)}{2}&\\
\tau^{(2)}_{q,0,0}&=&\frac{q^2(q-1)}{2}+\frac{q!}{3!(q-3)!}~\frac{3q-1}{4}&\tau^{(2)}_{q,1,0}&=&-\left(\frac{q(q-1)}{2}\right)^2&\\
\tau^{(2)}_{q,0,1}&=&-\frac{q^2(q-1)}{2}-q(q-1)&
\tau^{(2)}_{q,1,1}&=&q~\frac{q(q-1)}{2}&\\
 \tau^{(1)}_{q,1,0}&=&\frac{q(q-1)}{2}&
 \tau^{(1)}_{q,0,1}&=&q&\tau^{(1)}_{q,0,0}=-(\tau^{(1)}_{q,1,0}+\tau^{(1)}_{q,0,1})
\end{array}
\end{equation}
 
 \end{rem}
 \begin{defi}
We denote by $\boldsymbol{p_n}:=(p_0,\ldots,p_n)$ a given multi-index
in $\boldsymbol{\Ia_{n,q}}:=(\Ia_{q})^{n+1}$, with $p_k=(p_k^1,p_k^2)\in \Ia_{q}$ for any $0\leq k\leq n$. We also denote by $\boldsymbol{c_n}=(c_0,\ldots,c_n)$ a sequence of mappings in the set
$$
\boldsymbol{\Ja_{q,n}}=\cup_{\boldsymbol{p_n}\in\boldsymbol{\Ia_{n,q}}}\boldsymbol{I_{q}(p_n)}
\quad \mbox{with} \quad
\boldsymbol{I_{q}(p_n)}:=\prod_{0\leq k\leq n}I_q(p_k)
$$
For any 
$\boldsymbol{m_n}=\left(m_0,\ldots,m_n\right)\in\NN^{n+1}$, we set $\vert\boldsymbol{m_n}\vert=\sum_{0\leq k\leq n}m_k$, 
and we use the multi-index notation
$$\boldsymbol{\tau^{(m_n)}_{q,p_n}}=\prod_{0\leq k\leq n}\tau^{(m_k)}_{q,p_k^1,p^2_k},\ \
\boldsymbol{\tau}^{(m)}_{\boldsymbol{q,p_n}}:=
\sum_{\vert\boldsymbol{m_n}\vert =m}\boldsymbol{\tau^{(m_n)}_{q,p_n}},\ 
\ \boldsymbol{\Ta^{(m)}_{q,n}}:=\coprod_{\vert\boldsymbol{m_n}\vert =m}\prod_{0\leq k\leq n}{\Ta^{(m_k)}_{q,n}}
$$
and
$$
\boldsymbol{\Ca^{(q)}_{p_n}(c_n)}:=\prod_{0\leq k\leq n} \Ca^{(q)}_{p_k}(c_k)
\qquad
\boldsymbol{\Pa^{(N,q)}_n}(\boldsymbol{p_n}):=\prod_{0\leq k\leq n}
\Pa^{(N,q)}(p_k)
$$
\end{defi}
In this notation,  and recalling that $p_n^1+p_n^2>m_n\Rightarrow\boldsymbol{\tau^{(m_n)}_{q,p_n}}=0$, we readily prove the following extension
of lemma~\ref{key-combinatorial-lemma}
\begin{prop}\label{key-combinatorial-lemma-ext}
For any $q<N$ and $n\geq 0$, the mapping $N\mapsto \boldsymbol{\Pa^{(N,q)}_n}$ is differentiable at any order. In addition,
the $m$-th order derivative is  the signed measure with null mass
$$
\displaystyle d^{(m)}\boldsymbol{\Pa^{(q)}_n}=\sum_{\boldsymbol{p_n}\in \boldsymbol{\Ta^{(m)}_{q,n}}}
\boldsymbol{\tau}^{(m)}_{\boldsymbol{q,p_n}}~\delta_{\boldsymbol{p_n}}
$$
In addition, we have
\begin{equation}\label{estimate-1}
 \sum_{\boldsymbol{p_n}\in \boldsymbol{\Ta^{(m)}_{q,n}}}\left\vert \boldsymbol{\tau}^{(m)}_{\boldsymbol{q,p_n}}\right\vert\leq \frac{(m+n)!}{m!~n!}~(cq)^{2m}
 \end{equation}
 for some finite constant $c<\infty$.
\end{prop}
 {
\proof
 By Theorem 2 in~\cite{arratia}, for any $p\leq q$ we have the rather crude estimates
 $$
 S(q,q-p)\leq c~\frac{q^{2p}}{2^{p} p!}\leq c~q^{2p}
 \quad\mbox{\rm and}\quad \left\vert s(q,q-p)\right\vert\leq c~\left(\begin{array}{c}
 q\\
 q-p
 \end{array}\right)~\left(\frac{q-p}{2}\right)^{p}\leq c~q^{2p}
 $$
 for some finite constant $c<\infty$.  We also notice that
 \begin{equation}\label{some-crude-est}
  \left(\begin{array}{c}
q\\p_{2}
\end{array}\right)\leq q^{p_2}\quad \left(\begin{array}{c}
q-p_2\\k_2
\end{array}\right)\leq q^{k_2}\quad \mbox{\rm and}\quad
{(p_1+k_1)+k_3-1\choose k_3}\leq (2eq)^{k_3}
 \end{equation}
 To prove the r.h.s. estimate we use Stirling approximation to check that $${(p_1+k_1)+k_3-1\choose k_3}\leq (q+k_3)^{k_3}/k_3!\leq ({e(q+k_3)}/{k_3})^{k_3}\leq (e(q+1))^{k_3}\leq (2eq)^{k_3} $$
 
 Combining (\ref{some-crude-est}) and (\ref{ref-bound-tau}) with the fact that
 $$
\left\vert \alpha_{q,p_1,p_2}(\boldsymbol{k})\right\vert\leq (cq)^{2p_1+p_2+2k_1+k_2+k_3}\leq (cq)^{2m}
 $$
 for any $\boldsymbol{k}\in \Ka_q^{(m)}(p_1,p_2)$ and some finite constant $c<\infty$, 
 We conclude that
 $$
 \left\vert \tau^{(m)}_{q,p_1,p_2}\right\vert\leq \left(m-(p_1+p_2)\right)^2 \times (cq)^{2m}\leq m^2 (cq)^{2m}\leq (c^{\prime}q)^{2m}
 $$
 for some finite constant $c^{\prime}<\infty$. In much the same way
 $$
 \sum_{0\leq p_1+p_2\leq m}\left\vert \tau^{(m)}_{q,p_1,p_2}\right\vert\leq  m^2 (c^{\prime}q)^{2m}\leq (cq)^{2m}
 $$
 for some finite constant $c<\infty$. 
 This yields the rather crude estimates
 $$
 \forall m_0+\ldots+m_n=m\qquad
 \sum_{0\leq p^1_0+p^2_0\leq m}\ldots \sum_{0\leq p^1_n+p^2_n\leq m} \left\vert \boldsymbol{\tau}^{(m)}_{\boldsymbol{q,p_n}}\right\vert\leq  (cq)^{2m}
 $$
 The estimate (\ref{estimate-1})  comes from the fact that the cardinality of the set $\{(m_0,\ldots,m_n)~:~m_0+\ldots+m_n=m\}$
 coincides with the number $\frac{(m+n)!}{m!~n!}$ of finite multisets of size $m$ whose elements are drawn from a set of $(n+1)$ elements.
 \cqfd}

\begin{defi}\label{freeness}
 For further use, let ${\boldsymbol c}=(c_0,...,c_n)$, $c_i=(a_i,b_i)$ be a sequence of mappings in the set
$\boldsymbol{\Ja_{q,n}}$, and let
 us say that
 \begin{itemize}
  \item the $p$-th trajectory, $1\leq p\leq q$  of ${\boldsymbol c}$ is free if $\forall i\leq n, \forall m\not=p, $
  $$a_i\circ \ldots\circ a_n(p)\not=a_i\circ \ldots\circ a_n(m)\ \mbox{\rm{and}}\ b_i(a_{i+1}\circ\ldots\circ a_n(p))\not=1$$
  \item the $p$-th trajectory is coalescent if $\exists i\leq n, \exists m\not=p, a_i\circ \ldots\circ a_n(p)=a_i\circ \ldots\circ a_n(m)$
  \item the $p$-the trajectory is infected if $\exists i\leq n, b_i(a_{i+1}\circ\ldots\circ a_n(p))=1$.
 \end{itemize}
\end{defi}

\subsection{Unnormalized tensor product measures}\label{un-norm-sec}

Let us apply now these combinatorial results to PMCMC samplers. Our first result is concerned with tensor product measures.
Given a frozen trajectory $\textsl{z}:=\left(\textsl{z}_n\right)_{n\geq 0}\in \prod_{n\geq 0}S_n$, we denote by 
$\Xa_{\textsl{z},n}$ the dual
 mean field model associated with the Feynman-Kac particle model $\chi_n$ and the 
frozen path $X_n=\textsl{z}_n$. 

We also set $$\eta_{\textsl{z},n}^N:=m(\Xa_{\textsl{z},n})=m_{\textsl{z}_n}(\Xa_{\textsl{z},n}^-), \ \gamma_{\textsl{z},n}^N(1):=
\prod\limits_{0\leq p< n}\eta_{\textsl{z},p}^N(G_p),\ 
\gamma_{\textsl{z},n}^N:= \gamma_{\textsl{z},n}^N(1)\cdot \eta_{\textsl{z},n}^N,$$
and finally, for any function $F$ on $S^q_n$
$$\Upsilon_{\textsl{z},n}^{(N,q)}(F):=\EE\left((\gamma_{\textsl{z},n}^N)^{\otimes q}(F)\right)/{ \gamma}_n(1)^q.$$

\begin{defi}
We consider the tensor product measures
\begin{eqnarray}\label{tensormeas}
\Delta_{\textsl{z},\boldsymbol{p_n}}^{(q)}&=&\left(\eta^{\otimes q}_0\widehat{\CC}^{(q)}_{\textsl{z}_{0},p_0}\right)\left(\overline{Q}^{\otimes q}_{1}\widehat{\CC}^{(q)}_{\textsl{z}_{1},p_1}\right)\ldots
\left(\overline{Q}^{\otimes q}_{n}\widehat{\CC}^{(q)}_{\textsl{z}_{n},p_n}\right)=\sum_{\boldsymbol{c_n}\in \boldsymbol{I_{q}(p_n)}}~\boldsymbol{\Ca^{(q)}_{p_n}(c_n)}~\Delta_{\textsl{z},\boldsymbol{c_n}}^{(q)}
\end{eqnarray}
with the conditional expectation operators
$$
\Delta_{\textsl{z},\boldsymbol{c_n}}^{(q)}:=
\left(\eta^{\otimes q}_0\CC^{(q)}_{\textsl{z}_{0},c_0}\right)\left(\overline{Q}^{\otimes q}_{1}\CC^{(q)}_{\textsl{z}_{1},c_1}\right)\ldots
\left(\overline{Q}^{\otimes q}_{n}\CC^{(q)}_{\textsl{z}_{n},c_n}\right)
$$
\end{defi}

\begin{theo}\label{theo-v2-tensor}
For any $q<N$, $n\geq 0$, we have
\begin{eqnarray*}
\Upsilon_{\textsl{z},n}^{(N,q)}&=&
\sum_{\boldsymbol{p_n}\in \boldsymbol{\Ia_{n,q}}}~
\sum_{\boldsymbol{c_n}\in \boldsymbol{I_{q}(p_n)}}~\left[\boldsymbol{\Pa^{(N,q)}_n}(\boldsymbol{p_n})~\boldsymbol{\Ca^{(q)}_{p_n}(c_n)}\right]~\Delta_{\textsl{z},\boldsymbol{c_n}}^{(q)}
\end{eqnarray*}
\end{theo}

\proof

By construction, we have $\eta_{\textsl{z},n}^N:=m_{\textsl{z}_n}(\Xa_{\textsl{z},n}^-)$ and
$$
m_{\textsl{z}_n}(\Xa_{\textsl{z},n}^-)^{\otimes q}=m(\Xa_{\textsl{z},n}^-)^{\odot q}\CC^{(N,q)}_{\textsl{z}_n}
$$
On the other hand, for any function $F$ on $S^q_n$ we have
$$
\EE\left(m(\Xa_{\textsl{z},n+1}^-)^{\odot q}(F)~\left\vert~\Fa_{n}\right.\right)={\left(\eta^N_{\textsl{z},n}\right)^{\otimes q}\left(Q^{\otimes q}_{n+1}(F)\right)}/{\eta^N_{\textsl{z},n}(G_n)^q}
$$
This implies that
\begin{eqnarray*}
\EE\left(\left(\gamma^N_{\textsl{z},n+1}\right)^{\otimes q}(F)~\left\vert~\Fa_{n}\right.\right)&=&
\gamma^N_{\textsl{z},n}(1)^q\times \left(\eta^N_{\textsl{z},n}\right)^{\otimes q}\left(Q^{\otimes q}_{n+1}\CC^{(N,q)}_{\textsl{z}_{n+1}}(F)\right)\\
&=&\left(\gamma^N_{\textsl{z},n}\right)^{\otimes q}\left(Q^{\otimes q}_{n+1}\CC^{(N,q)}_{\textsl{z}_{n+1}}(F)\right)
\end{eqnarray*}
from which we conclude that
\begin{eqnarray*}
\Upsilon_{\textsl{z},n}^{(N,q)}(F)&=&\left(\eta^{\otimes q}_0\CC^{(N,q)}_{\textsl{z}_{0}}\right)\left(\overline{Q}^{\otimes q}_{1}\CC^{(N,q)}_{\textsl{z}_{1}}\right)\ldots
\left(\overline{Q}^{\otimes q}_{n}\CC^{(N,q)}_{\textsl{z}_{n}}\right)(F).
\end{eqnarray*}
The Theorem follows by expanding the $\CC^{(N,q)}_{\textsl{z}_i}$ in terms of the $\CC^{(q)}_{\textsl{z}_i,c_i}$.
\cqfd

The next corollary is a direct consequence of the proof of theorem~\ref{theo-v2-tensor}. It provides a more probabilistic description of the measure $\Upsilon_n^{(N,q)}$ in terms of expectation operators.
\begin{cor}
For any $q<N$, $n\geq 0$,  $\Upsilon_{\textsl{z},n}^{(N,q)}$  is differentiable at any order. In addition, its derivatives are for any $n\geq 0$ given by
the recursion
$$
d^{(m)}\Upsilon_{\textsl{z},n}^{(q)}(F)=\sum_{r_1+r_2=m}~\sum_{p\in \Ia_{q}}~d^{(r_1)}\Pa^{(q)}(p)~~d^{(r_2)}\Upsilon_{\textsl{z},n-1}^{(q)}\left(\overline{Q}_{n}^{\otimes q}\widehat{\CC}^{(q)}_{\textsl{z}_n,p}(F)\right)
$$
with the conventions $\Upsilon_{\textsl{z},-1}^{(q)}\overline{Q}_{0}^{\otimes q}=\eta_0^{\otimes q}$ 
and $d^{(r_2)}\Upsilon_{\textsl{z},-1}^{(q)}\overline{Q}_{0}^{\otimes q}=0$ for $r_2>0$.
In particular we get the expansions
\begin{eqnarray}\label{expansion-dmQ-ref}
d^{(m)}\Upsilon_{\textsl{z},n}^{(q)}&=&\sum_{\boldsymbol{p_n}\in\boldsymbol{\Ta^{(m)}_{q,n}}}\boldsymbol{\tau}^{(m)}_{\boldsymbol{q,p_n}}
\times\Delta_{\textsl{z},\boldsymbol{p_n}}^{(q)}.
\end{eqnarray}
\end{cor}

For further use, let us study further the action of the operators $\Delta_{\textsl{z},{\boldsymbol c_n}}^{(q)}$. We already know that they contribute to $d^{(m)}\Upsilon_{\textsl{z},n}^{(q)}$ only if the total number of coalescences and infections of ${\boldsymbol c_n}$, written $Tot({\boldsymbol c_n})$ is less than $m$.

\begin{lem}\label{letruc}
 {For any $\boldsymbol{p_n}\in\boldsymbol{\Ta^{(m)}_{q,n}}$, $m\geq 1$, and $n^{\prime}\geq n$  we have
\begin{equation}\label{estimate-2}
\left\Vert  \Delta_{\textsl{z},\boldsymbol{p_n}}^{(q)}\overline{Q}^{\otimes q}_{n,n^{\prime}}\right\Vert_{\tiny tv}\leq \overline{g}_{n^{\prime}}^{2m}\quad\mbox{\rm and}\quad 
\left\Vert d^{(m)}\Upsilon_{\textsl{z},n}^{(q)}\overline{Q}^{\otimes q}_{n,n^{\prime}}\right\Vert_{\tiny tv}\leq~(cn(q \overline{g}_{n^{\prime}})^2)^{m}
 \end{equation}
for some finite constant $c<\infty$.}
In addition, let $f$ a $\eta_n$-centered function on $S_n$ ($\eta_n(f)=0$). Then, for any sequence of mappings ${\boldsymbol c_n}$,
 $$Tot({\boldsymbol c_n})< \frac{q}{2}\Rightarrow \Delta_{\textsl{z},{\boldsymbol c_n}}^{(q)}(f^{\otimes q})=0.$$
 In particular, $d^{(m)}\Upsilon_{\textsl{z},n}^{(q)}(f^{\otimes q})=0$ whenever $m<\frac{q}{2}.$
\end{lem}
\proof
 {The first assertion comes from the fact that for any ${\boldsymbol c_n}$ with $p^1_k$-coalescences and $p^2_k$ infections at levels $0\leq k\leq n$ we have
the rather crude estimates
$$
\left\Vert  \Delta_{\textsl{z},\boldsymbol{c_n}}^{(q)}(\overline{Q}^{\otimes q}_{n,n^{\prime}}(F))\right\Vert\leq \prod_{0\leq k\leq n}\overline{g}_{n^{\prime}}^{~2p^1_k+p^2_k}\leq \overline{g}_{n^{\prime}}^{~2\left\vert\boldsymbol{p_n}\right\vert}$$
with
$$
\left\vert\boldsymbol{p_n}\right\vert:=\sum_{0\leq k\leq n}(p^1_k+p^2_k)\leq \sum_{0\leq k\leq n}m_k=m
$$
for any function $F$ on $S_{n^{\prime}}^q$ s.t. $\Vert F\Vert\leq 1$.
The end of the proof of (\ref{estimate-2}) is a direct consequence of (\ref{estimate-1}), (\ref{tensormeas}) and
$$
\frac{(m+n)!}{m!n! n^m}\leq \frac{(n+m)^m}{n^mm!}\leq \frac{m^m}{m!}~\left(\frac{1}{m}+\frac{1}{n}\right)\leq e^{2m}
$$}
To prove the second assertion, let us assume that $Tot({\boldsymbol c_n})< \frac{q}{2}$. It follows immediately that one trajectory is free in the sense of Definition 
\ref{freeness}. Because of the symmetry of the problem (which, as usual, is invariant by permutation of the particles), we may assume without restriction that the particles of this free trajectory all have the same index $q$ ($a_i(q)=q\ \forall i\leq n$).
Let us write $\hat{\boldsymbol c}_n$ for the sequence of mappings obtained by restricting each $a_i$ to a map from $[q-1]$ to itself (this process is well-defined because of the freeness asumption) and by restricting similarly $b_i$ to $[q-1]$. It follows then from the very definition of $\Delta_{\textsl{z},{\boldsymbol c_n}}^{(q)}(f^{\otimes n})$ that
$$\Delta_{\textsl{z},{\boldsymbol c_n}}^{(q)}(f^{\otimes q})=\Delta_{\textsl{z},\hat{\boldsymbol c}_n}^{(q-1)}(f^{\otimes q-1})\cdot \eta_n(f)=0.$$
This ends the proof of the Lemma.
\cqfd

\begin{cor}\label{lautretruc}
 We have for an arbitrary $q\leq N$:
 $$\EE[(\gamma_{\textsl{z},n}^N(G_n-\eta_n(G_n))^{q}]=\EE[(\gamma_{\textsl{z},n}^N)^{\otimes q}((G_n-\eta_n(G_n))^{\otimes q})]=O(N^{-q/2}).$$
\end{cor}

\begin{cor}\label{encoreautretruc}
 We have for an arbitrary $q\leq N$:
 $$\EE[(\gamma_{\textsl{z},n}^N(G_n)-\gamma_n(G_n))^{q}]=O(N^{-q/2}).$$
\end{cor}

Indeed, 
$$ \gamma_{\textsl{z},n}^N(G_n)-\gamma_n(G_n)=\prod_{0\leq p\leq n}\eta_{\textsl{z},p}^N(G_p)-\prod_{0\leq p\leq n}\eta_{p}(G_p)$$
$$=\gamma_{\textsl{z},n}^N(G_n-\eta_{n}(G_n))+[\gamma_{\textsl{z},n-1}^N(G_{n-1})-\gamma_{n-1}(G_{n-1})]\eta_{n}(G_n)$$
$$=\sum_{i=0}^n[\gamma_{\textsl{z},i}^N(G_{i}-\eta_{i}(G_{i})]\prod\limits_{j=i+1}^n\eta_{i}(G_i).$$
The proof follows from the previous Corollary and the Minkowski identity.

\subsection{Normalized tensor product measures}\label{norm-FK-ref}

In the present paragraph, we show that the distribution $\PP^{(N,q)}_{\textsl{z},n+1}$  of the first $q$ random non frozen particles (see definition \ref{nonfrozen}) has derivatives at all orders.

We recall the intrumental identity:
for any $u\not= 1,\ q\geq 0$ and $m\geq 1$
\begin{equation}\label{key-fraction}
\frac{1}{(1-u)^{q+1}}=\sum_{0\leq k\leq m} (q+k)_k~\frac{u^k}{k!}+ u^m~\sum_{1\leq k\leq q+1}\left(
\begin{array}{c}
(q+1)+m\\
k+m
\end{array}
\right)~\left(\frac{u}{1-u}\right)^k
\end{equation}
A detailed proof of this result can be found in~\cite{dpr-2009} (cf. lemma 4.11 on page 820).

Using the identity ${n+1\choose k}=\sum\limits_{k\leq l\leq n}{n\choose l}$ (following e.g. from $1-(1-x)^{n+1}=x\sum\limits_{0\leq k\leq n}(1-x)^k$), we get 
\begin{equation}\label{key-fraction2}
 \frac{1}{x^q}=\frac{(q+r)!}{(q-1)!}\sum\limits_{0\leq l\leq r}\frac{1}{(q+l)}\frac{(-1)^l}{l!(r-l)!}~x^l+\sum_{1\leq k\leq q}\left(
 \begin{array}{c}
q+r\\
k+r
\end{array} \right)~\frac{(1-x)^{r+k}}{x^k}
\end{equation}

\begin{theo}\label{final-theo}
 {There exists some finite constant $c>0$ such that for any $n\geq 0$, $r>m \geq 1$  and 
$N>cn((r+q) \overline{g}_{n+1})^2$ we have 
$$
\left\Vert\PP^{(N,q)}_{\textsl{z},n+1}-\eta_{n+1}^{\otimes q}-\sum_{1\leq k\leq m} \frac{1}{N^k}~d^{(k)}\PP^{(q)}_{\textsl{z},n+1}\right\Vert_{\tiny tv}\leq 
b(q)~\left(\frac{cn((q+r) \overline{g}_{n+1})^2}{N}\right)^{m+1}
$$
with some finite constant $b(q)<\infty$ whose values only depend on the parameters $q$,} and the $k$-th order derivatives given for any function $F$ on $S_n^q$ by 
\begin{equation}\label{final-formula}
d^{(k)}\PP^{(q)}_{\textsl{z},n+1}(F)=\frac{(q+2k)!}{(q-1)!}~\sum_{0\leq l\leq 2k} \frac{(-1)^l}{(q+l)}~\frac{1}{l!~(2k-l)!}
~d^{(k)}\Upsilon^{(l+q)}_{\textsl{z},n}\left[\overline{Q}^{\otimes (l+q)}_{n,n+1}(1^{\otimes l}\otimes F)\right]
\end{equation}
 {In addition, we have
\begin{equation}\label{estimate-3}
\left\Vert d^{(k)}\PP^{(q)}_{\textsl{z},n+1} \right\Vert_{\tiny tv}\leq \left[cn((q+2k) \overline{g}_{n+1})^2\right]^{k}
 \end{equation}
for some finite constant $c<\infty$.}

\end{theo}

\proof

 {The proof of (\ref{estimate-3}) is a direct consequence of (\ref{estimate-2}). We set
$\overline{\gamma}^N_{z,n}(f)={\gamma}^N_{z,n}(f)/\gamma_n(1)$. In this notation, we have
$$
{\gamma}^N_{z,n}(\overline{G}_n)-\gamma_n(\overline{G}_n)=\gamma_n(1)\left(\overline{\gamma}^N_{z,n}(\overline{G}_n)-1\right)
$$
By Corollary~\ref{encoreautretruc}, for any even integer $q$ we have
\begin{eqnarray*}
\EE\left(\left(\overline{\gamma}^N_{z,n}(\overline{G}_n)-1\right)^q\right)&=&\sum_{k\geq q/2} \frac{1}{N^k}   \sum_{0\leq p\leq q}\left(\begin{array}{c}
q\\p\end{array}\right)(-1)^{q-p}~d^{(k)}\Upsilon^{(p)}_{\textsl{z},n}\left(\overline{G}_n^{\otimes p}\right)\\
&\leq &2^q\sum_{k\geq q/2}\left(\frac{cn(q\overline{g}_{n+1})^2}{N}\right)^{k}
\end{eqnarray*}
for some finite constant $c<\infty$ and for any $N\geq cn(q \overline{g}_{n+1}^2)^2$. The r.h.s. estimate is readily checked recalling that
$\overline{Q}_{n,n+1}^{\otimes p}(1)=\overline{G}_n^{\otimes p}$, and applying  (\ref{estimate-2}).

Thus, there exists some finite universal constant $c<\infty$ such that
$$
 \EE\left(\left(\overline{\gamma}^N_{z,n}(\overline{G}_n)-1\right)^q\right)^{1/q}
\leq \left({cn(q \overline{g}_{n+1})^2}/{N}\right)^{1/2}
$$
and
$$
 \EE\left(\left(\overline{\gamma}^N_{z,n}(\overline{G}_n)\right)^q\right)^{1/q}
\leq 1+\left({cn(q \overline{g}_{n+1})^2}/{N}\right)^{1/2}\leq 2
$$
as soon as $N\geq cn(q \overline{g}_{n+1})^2$.
}
Following the proof of theorem~\ref{theo-v2-tensor} we find that
$$
\begin{array}{l}
\EE\left(m(\Xa^-_{\textsl{z},n+1})^{\odot q}(F)\right)=\EE\left[\overline{\gamma}^N_{\textsl{z},n}(\overline{G}_n)^{-q}\times 
\left(\overline{\gamma}^N_{\textsl{z},n}\right)^{\otimes q}\left(\overline{Q}^{\otimes q}_{n,n+1}(F)\right)\right]
\end{array}$$
 {On the other hand, for any $1\leq k\leq q$, $r\geq 1$ and $\Vert F\Vert\leq 1$ we have
$$
\begin{array}{l}
\left\vert\EE\left[\overline{\gamma}^N_{\textsl{z},n}(\overline{G}_n)^{-k}\times \left(1-\overline{\gamma}^N_{\textsl{z},n}(\overline{G}_n)\right)^{r+k}
\left(\overline{\gamma}^N_{\textsl{z},n}\right)^{\otimes q}\left(\overline{Q}^{\otimes q}_{n,n+1}(F)\right)\right]\right\vert\\
\\
\leq \EE\left[\overline{\gamma}^N_{\textsl{z},n}(\overline{G}_n)^{q-k}\times \left\vert1-\overline{\gamma}^N_{\textsl{z},n}(\overline{G}_n)\right\vert^{r+k}\right]\\
\\
\displaystyle\leq 2^{q-k}
\left(
\EE\left[ \left(\overline{\gamma}^N_{\textsl{z},n}(\overline{G}_n)-1\right)^{2(r+k)}\right]\right)^{1/2}\leq 2^{q}\left({cn((r+k)
 \overline{g}_{n+1})^2}/{N}\right)^{(r+k)/2}
\end{array}
$$
as soon as  $N\geq cn((r+q) \overline{g}_{n+1})^2$.
Recalling that
$$
\overline{G}_n^{\otimes l}\otimes \left(\overline{Q}^{\otimes q}_{n,n+1}(F)\right)= \overline{Q}^{\otimes (q+l)}_{n,n+1}(1^{\otimes l}\otimes F)
$$
and combining (\ref{key-fraction2}) with Corollary \ref{encoreautretruc} we find that
$$
\begin{array}{l}
\displaystyle\left\vert\PP^{(N,q)}_{\textsl{z},n+1}(F)-\frac{(q+r)!}{(q-1)!}\sum_{0\leq l\leq r} \frac{1}{(q+l)}~\frac{(-1)^l}{l!~(r-l)!}~
\Upsilon^{(N,l+q)}_{\textsl{z},n}\left[\overline{Q}^{\otimes (l+q)}_{n,n+1}(1^{\otimes l}\otimes F)\right]
\right\vert\\
\\
\displaystyle\leq a(q)~\left({cn((r+q) \overline{g}_{n+1})^2}/{N}\right)^{(r+1)/2}
\end{array}$$
for any $r\geq 0$ and some finite constant $a(q)<\infty$. We obtain the $k$-th differential operator formulae (\ref{final-formula}) by choosing $r=2k$. After some elementary manipulations, this yields the estimate
$$
\left\Vert\PP^{(N,q)}_{\textsl{z},n+1}-\eta_{n+1}^{\otimes q}-\sum_{1\leq k\leq r} \frac{1}{N^k}~d^{(k)}\PP^{(q)}_{\textsl{z},n+1}\right\Vert_{\tiny tv}\leq 
b(q)~\left({cn((q+r) \overline{g}_{n+1})^2}/{N}\right)^{r+1/2}
$$
with some finite constant $b(q)<\infty$ and the derivatives operators $d^{(k)}\PP^{(q)}_{\textsl{z},n+1}$ given in (\ref{final-formula}). The end of  the proof of the theorem is now easily completed.
 }
\cqfd

It is instructive to derive explicit expressions for the derivatives --this will be one of the topics addressed in the forthcoming paragraphs.
Let us anticipate on these developments and make explicit the first order derivative in a simple case.
For $k=q=1$, and any function $f$ on $S_n$, with $\eta_n(f)=0$, 
using the first order expansions that will be stated in corollary~\ref{cor-Q-expansions} 
it is readily checked that
\begin{eqnarray*}
d^{(1)}\PP^{(1)}_{\textsl{z},n+1}(f)=\sum_{0\leq k\leq n} \overline{Q}_{k,n+1}(f)(\textsl{z}_k)
-\sum_{0\leq k\leq n} 
\eta_k\left(\overline{Q}_{k,n+1}(1)\overline{Q}_{k,n+1}(f)\right) .
\end{eqnarray*}

\subsection{Infected forest expansions}\label{ifexp-ref}
We know that $\PP^{(q,N)}_{\textsl{z},n}$ has derivatives at all orders and can be expanded in terms of the derivatives of $\Upsilon^{(N)}_{\textsl{z},n}$. In turn, these last derivatives can be expanded in terms of the elementary integral operators $\Delta_{\textsl{z},n,\boldsymbol{c}}^{(q)}$. However, because of the symmetries of Feynman-Kac models, many of these operators coincide and this expansion is not efficient, neither computationally nor theoretically. The present paragraph aims at clarifying these questions and get rid of redundancies in combinatorial expansions of derivatives.

The results in this paragraph build largely on \cite{dpr-2009}. We will therefore skip the details of the arguments that follow closely the ones in \cite{dpr-2009} and refer simply the reader to that article for further details on the definitions, proofs, reasonings and so on on trees, forests and jungles.

\subsubsection{Forests and jungles}

We start with recalling some more or less classical terminology on trees and forests introduced in~\cite{dpr-2009}.

A tree is a (isomorphism class of) finite
non-empty oriented connected 
graph ${\bf t}$ without loops such that any
vertex of ${\bf t}$ has at most one outgoing
edge. Paths are oriented from the vertices to the root. The height of a tree is the maximum lenght of a path. Similarly, the level of a vertex in a tree is the length of the path that connects it to the root. These definition will extend in a straightforward way to the objects to be introduced below (forests and jungles).

A forest ${\bf f} $ is a multi-set of trees, that is a set of trees, with repetitions of the same tree allowed, or equivalently an element of the commutative monoid $\langle \Ta\rangle$ on $\Ta$,
with the empty graph  $T_0=\emptyset$ as a unit. Since the algebraic notation is the most convenient, we write
${\bf f} ={\bf t}_1^{m_1}...{\bf t}_k^{m_k}$,
for the forest with the tree ${\bf t}_i$ appearing with multiplicity $m_i$, $i\leq k$. When ${\bf t}_i\not= {\bf t}_j$ for $i\not= j$, we say that $\bf f$ is written in normal form.

The sets of forests with height $(n+1)$, and with $q$ vertices at each level set
is written $\Fo_{q,n}$. 

To a sequence ${\bf a}=(a_0,\ldots,a_n)\in \Aa_{q,n}:=([q]^{[q]})^{n+1}$ is naturally associated a forest $F({\bf a})$: the one with one vertex for each element of $[q]^{n+1}$, and a edge for each pair $(i,a_k(i)), i\in [q]$. The sequence can also be represented graphically uniquely by a planar graph $J({\bf a})$, where however the edges between vertices at level $k+1$ and $k$ are allowed to cross.
 We call such a planar graph, where paths between vertices are entangled, {\em a jungle}. The set of such jungles is written $\Ju_{q,n}$. 
Here is the graphical representation of a jungle (for consistency with the probabilistic interpretation of heights and levels as time-indices, we represent trees, forest and jungles \it horizontally \rm and from left to right --roots are on the left !).
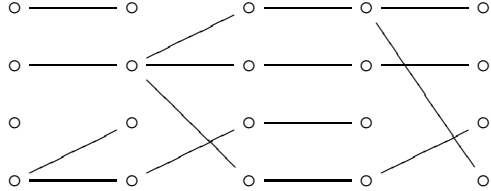
\begin{center}
\hskip.2cm
\xymatrix@C=3em@R=1em{ 
\circ\ar@{-}[r]&\circ&\circ\ar@{-}[r]\ar@{-}[ld]&\circ\ar@{-}[r]\ar@{-}[rddd]&\circ&\\ 
\circ\ar@{-}[r]&\circ\ar@{-}[r]\ar@{-}[rdd] &\circ\ar@{-}[r]&\circ\ar@{-}[r]&\circ\\  
\circ&\circ\ar@{-}[ld]&\circ\ar@{-}[r]&\circ&\circ\ar@{-}[ld]\\ 
\circ\ar@{-}[r]&\circ\ar@{-}[ur]&\circ\ar@{-}[r]&\circ&\circ
                                 }
\end{center}

 The group $\boldsymbol{\Ga_{q,n}}:=\Ga_q^{n+2}$ also acts naturally on sequences of maps ${\bf a\in
{\Aa}_{q,n}}$, and on jungles $J({\bf a})\in \Ju_{q,n}$ by permutation of the vertices at each level. More precisely, for all ${\bf a}\in \Aa_{n,q}$ and all ${\boldsymbol \sigma}=(\sigma_0,...,\sigma_{n+1})\in \boldsymbol{\Ga_{q,n}}$ by the pair of formulae
\begin{equation}\label{action}
{\boldsymbol \sigma} ({\bf a}):=( \sigma_0a_0\sigma^{-1}_1,\sigma_1a_1\sigma^{-1}_2,...,\sigma_{n}a_n\sigma^{-1}_{n+1})
\quad\mbox{\rm and}\quad
{\boldsymbol \sigma}J({\bf a}):=J({\boldsymbol\sigma}({\bf a}))
\end{equation}
Notice that if two sequences ${\bf a}$ and $\bf a^{\prime}\in\Aa_{q,n}$ differ only by the order of the vertices in $J({\bf a})$ and $J({\bf  a^{\prime}})$ (i.e. by the action of an element of $\boldsymbol{\Ga_{q,n}}$) then the associated forests are identical: $F({\bf a})=F({\bf  a^{\prime}})$.  
The converse is also true: if $F({\bf a})=F({\bf a^{\prime}})$, then $J({\bf a})$ and $J({\bf a^{\prime}})$ differ only by the ordering of the vertices, since they have the same underlying non planar graph. In this situation, ${\bf a}$ and ${\bf a^{\prime}}$ belong to the same orbit 
$$
\left[{\bf a}\right]:=\left\{\boldsymbol{\sigma}({\bf a})~:~\boldsymbol{\sigma}\in \boldsymbol{\Ga_{q,n}}\right\}
$$
under the action of $\boldsymbol{\Ga_{q,n}}$. In particular,
the set of equivalence classes of jungles in ${\Ju}_{q,n}$ under the action of the permutation groups $\boldsymbol{\Ga_{q,n}}$ is in bijection 
with both the set of $\boldsymbol{\Ga_{q,n}}$-orbits of maps in $\Aa_{q,n}$ and the set of forests ${\Fo}_{q,n}$. 
Writing $$\mbox{\rm Stab}\left(\boldsymbol{a}\right):=\left\{\boldsymbol{\tau}\in \boldsymbol{\Ga_{q,n}}~:~\boldsymbol{\tau}(\boldsymbol{a})=\boldsymbol{a}\right\}$$ for the stabilizer subgroup of $\boldsymbol{a}$, the class formula yields
$$\#[{\bf a}]=\#\boldsymbol{\Ga_{q,n}}/\# \mbox{\rm Stab}\left(\boldsymbol{a}\right)=(q!)^{n+2}/\# \mbox{\rm Stab}\left(\boldsymbol{a}\right).$$

To compute the cardinality of the set $\mbox{\rm Stab}\left(\boldsymbol{a}\right)$ in terms of forests and trees, we denote by $\mbox{\rm Cut}({\bf t})$ the forest deduced from cutting the root of the tree
${\bf t}$; that is, removing its root vertex, and all its incoming edges. In the reverse angle, we denote by $\mbox{\rm Cut}^{-1}({\bf f})$ the tree deduced from 
the forest ${\bf f}$ by adding a common root to its rooted tree. The symmetry multiset  ${\bf S}({\bf t})$ of a tree ${\bf t}=\mbox{\rm Cut}^{-1}({\bf t}_1^{m_1}\ldots {\bf t}_k^{m_k})$ associated with a forest written in normal form, is 
defined by
${\bf S}({\bf t}):=(m_1,\ldots,m_k)
$.
The symmetry multiset of a forest is given by
$$
{\bf S}({\bf t}_1^{m_1}\ldots {\bf t}_k^{m_k}):=\left(\underbrace{{\bf S}({\bf t}_1),\ldots,{\bf S}({\bf t}_1)}_{m_1-\mbox{terms}},\ldots,
\underbrace{{\bf S}({\bf t}_k),\ldots,{\bf S}({\bf t}_k)}_{m_k-\mbox{terms}}
\right)
$$ 
We also extend $\mbox{\rm Cut}({\bf f})$ to forests ${\bf f}={\bf t}_1^{m_1}\ldots {\bf t}_k^{m_k}$ by setting
\begin{equation}\label{cut-forests}
\mbox{\rm Cut}({\bf f})=\mbox{\rm Cut}({\bf t}_1)^{m_1}\ldots \mbox{\rm Cut}({\bf t}_k)^{m_k}
\end{equation}
where $\mbox{\rm Cut}({\bf t}_i)^{m_i}$  stands for the forest deduced from $\mbox{\rm Cut}({\bf t}_i)$ repeated $m_i$ times.
Combining the class formula with a recursion
with respect to the height parameter, we obtain
\begin{equation}\label{exact-formula-Fa}
\#\left([\boldsymbol{a}]\right)={(q!)^{n+2}}/\#\left(\mbox{\rm Stab}(\boldsymbol{a})\right)
\quad\mbox{\rm with}\quad
\#\left(\mbox{\rm Stab}(\boldsymbol{a})\right)=\prod\limits_{i={-1}}^{n}{\bf S}(\mbox{\rm Cut}^i(F(\boldsymbol{a})))!
\end{equation}
where we have used the multi-index factorial notation $(n_1,\ldots,n_k)=n_1!\ldots n_k!$, for any $n_k\in\NN$ , with $k\geq 0$.
A detailed proof of this closed  formula is provided in~\cite{dpr-2009}.

We let the reader check that, for example, for $\bf a$ as in the above graphical representation, $\#\left(\mbox{\rm Stab}(\boldsymbol{a})\right)=1\cdot 1\cdot 2!\cdot 2!=4$ and
$\#\left([\boldsymbol{a}]\right)=(4!)^4\cdot 3!$.

\subsubsection{Infected forests}
Recall that the study of PMCMC samplers requires the introduction of sequences of mappings $\boldsymbol c=(\boldsymbol a,\boldsymbol b)\in \boldsymbol{\Ja_{q,n}}$, where the maps $b_k$ can be thought of as ``infections'' (using the terminology previously introduced).
The infection of a jungle $J({\boldsymbol a})$ (or of the associated sequence of maps $\boldsymbol a$)  is defined accordingly by a sequence of functions
  $
 {\bf b}=(b_0,\ldots,b_n)\in (\{0,1\}^{[q]})^{n+1}$.
 
Graphically, the infection is represented by the label $1$, and the non infection by the label $0$ on the edges of the jungle.
 The diagram below provides an example of infected planar forest of height $4$ with $5$ trees and $6$ leaves, 
and the corresponding sequence of infection mappings. 
 \begin{center}
\hskip.3cm
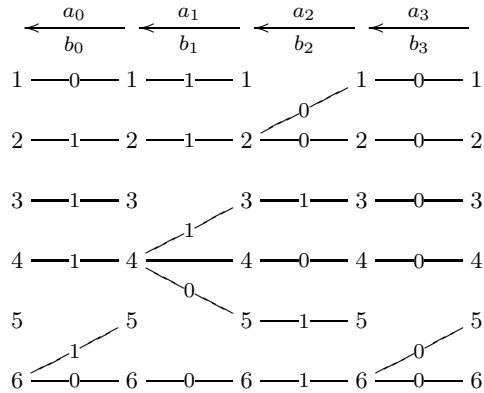
\xymatrix@C=3em@R=1em{     
& \ar[l]_{a_0}^{b_0}&\ar[l]_{a_1}^{b_1}&\ar[l]_{a_2}^{b_2}  & \ar[l]_{a_3}^{b_3}&\\
\text{\footnotesize 1}\ar@{-}[r]|{0}&\text{\footnotesize 1}\ar@{-}[r]|{1}
&\text{\footnotesize 1}&\text{\footnotesize 1}\ar@{-}[r]|{0} &\text{\footnotesize 1}\\
\text{\footnotesize 2}\ar@{-}[r]|{1}&\text{\footnotesize 2}\ar@{-}[r]|{1} &\text{\footnotesize 2}\ar@{-}[r]|{0} \ar@{-}[ur]|{0}
&\text{\footnotesize 2}\ar@{-}[r]|{0}&\text{\footnotesize 2} \\     
\text{\footnotesize 3}\ar@{-}[r]|{1}&\text{\footnotesize 3}&\text{\footnotesize 3}\ar@{-}[r]|{1}\ar@{-}[ld]|{1}
&\text{\footnotesize 3}\ar@{-}[r]|{0}&\text{\footnotesize 3}&\\ 
\text{\footnotesize 4}\ar@{-}[r]|{1}&\text{\footnotesize 4}\ar@{-}[r]\ar@{-}[rd]|{0} &\text{\footnotesize 4}\ar@{-}[r]|{0}&\text{\footnotesize 4}\ar@{-}[r]|{0}&\text{\footnotesize 4}\\  
\text{\footnotesize 5}&\text{\footnotesize 5}\ar@{-}[ld]|{1}&\text{\footnotesize 5}\ar@{-}[r]|{1}&\text{\footnotesize 5}&\text{\footnotesize 5}\ar@{-}[ld]|{0}\\ 
\text{\footnotesize 6}\ar@{-}[r]|{0}&\text{\footnotesize 6}\ar@{-}[r]|{0}&\text{\footnotesize 6}\ar@{-}[r]|{1}&\text{\footnotesize 6}\ar@{-}[r]|{0}&\text{\footnotesize 6}
                                 }
\end{center}

By construction, there are  $\prod_{0\leq k\leq n}\left(\begin{array}{c}
q\\
i_k
\end{array}\right)$ ways of infecting a given forest with $0\leq i_k\leq q$ infections at each level $0\leq k\leq n$.
Some of them are clearly equivalent. To be more precise, we consider the following equivalence relation on infected jungles
\begin{eqnarray*}
(\boldsymbol{a},\boldsymbol{b})\sim (\boldsymbol{a^{\prime}},\boldsymbol{b^{\prime}})&\Longleftrightarrow&
\exists  \boldsymbol{\sigma}\in\boldsymbol{\Ga_{q,n}}~:~ \boldsymbol{\sigma} (\boldsymbol{a},\boldsymbol{b})=
(\boldsymbol{a^{\prime}},\boldsymbol{b^{\prime}})
\end{eqnarray*} 
The equivalence classes are denoted by
\begin{eqnarray*}
\left[\boldsymbol{a},\boldsymbol{b}\right]&:=&\left\{\boldsymbol{\sigma} (\boldsymbol{a},\boldsymbol{b})~:~\boldsymbol{\sigma}\in\boldsymbol{\Ga_{q,n}}~\right\}
=\left\{(\boldsymbol{\sigma} (\boldsymbol{a}),\boldsymbol{b}\boldsymbol{\overline{\sigma}}^{-1})~:~\boldsymbol{\sigma}\in\boldsymbol{\Ga_{q,n}}~\right\}
\end{eqnarray*}
with
$$
\overline{\boldsymbol{\sigma}}:=\left(\sigma_1,\ldots,\sigma_{n+1}\right)\quad\mbox{\rm and}\quad
\overline{\boldsymbol{\sigma}}^{-1}=\left(\sigma_1^{-1},\ldots,\sigma_{n+1}^{-1}\right)
$$

The definitions of forests and jungles discussed in the previous section extend also in a straightforward way to the infected case (edges being colored by 0 or 1). To a sequence ${\bf (a,b)}$ is then naturally associated an infected forest $F({\bf a,b})$: the one with one vertex for each element of $[q]^{n+1}$, and an infected edge for each triplet $(i,b_k(i),a_k(i)), i\in [q]$. 
The set of infected forests is in bijection with the 
 set of $\boldsymbol{\Ga_{q,n}}$-orbits  of maps  in $\boldsymbol{\Ja_{q,n}}$. 

The class formula yields once again a way to compute the cardinals of the classes $[\boldsymbol{a},\boldsymbol{b}]$ from the action of the symmetry group $\boldsymbol{\Ga_{q,n}}$.
\begin{lem}\label{card-jungles}
The number of infected jungles in $[\boldsymbol{a},\boldsymbol{b}]$ is given by
$$
\#\left[\boldsymbol{a},\boldsymbol{b}\right]=(q!)^{n+2}/\mbox{\rm Stab}_{\boldsymbol{a}}(\boldsymbol{b})=
\#[\boldsymbol{a}]\times \frac{\#\left(\mbox{\rm Stab}(\boldsymbol{a})\right)}{\#\left(\mbox{\rm Stab}_{\boldsymbol{a}}(\boldsymbol{b})\right)}
$$
with
$$
\quad
\mbox{\rm Stab}_{\boldsymbol{a}}\left(\boldsymbol{b}\right):=\left\{
\boldsymbol{\tau}\in \mbox{\rm Stab}\left(\boldsymbol{a}\right)~:~\boldsymbol{b}\boldsymbol{\overline{\tau}}=\boldsymbol{b}
\right\}.
$$
\end{lem}

As for the non infected case, $\#(\mbox{\rm Stab}_{\boldsymbol{a}}\left(\boldsymbol{b}\right))$ can be computed inductively, following essentially  the same principles. We describe briefly how this can be done.

Let $\bf t_1, ..., t_n$ and $\bf t_1',...,t_m'$ be two families of distinct infected trees and $l_i,\ i=1...n,\ p_j,\ j=1...m$ two sequences of positive integers.
We write $\bf t_1^{\rm l_1} ... t_n^{\rm l_n}\circledast  t_1^{'\rm p_1} ... t_n^{'\rm p_m}$ for the infected tree obtained by joigning, for $i=1...n$, $l_i$ copies of $\bf t_i$ to a common root with infection index 0 and for $i=1...m$, $p_i$ copies of $\bf t_i'$ to the same common root with infection index 1.
Any infected tree $\bf t$ can be written uniquely in that way: we write $\bf S'(t)=\rm (l_1,...,l_n,p_1,...,p_m)$ for the corresponding multiset and call it the symmetry multiset of $\bf t$. 

Cuts of infected trees and infected forests are infected forests that are defined as in the non infected case by removing the root and erasing all infected edges connected to the root.
A (right only) inverse operation $Cut^{-1}$ acting on an infected forest $\bf t_1^{\rm k_1} ... t_n^{\rm k_n}$ is defined by linking all the infected trees to a common root with non infected edges.

Mimicking the inductive arguments for counting jungles using cardinals of stabilizers in~\cite{dpr-2009}, we get
\begin{equation}\label{stabil}
\mbox{\rm Stab}_{\boldsymbol{a}}\left(\boldsymbol{b}\right)=\prod_{i=-1}^n{\bf S^{\prime}}(\mbox{\rm Cut}^i([\boldsymbol{a},\boldsymbol{b}]))!
\end{equation}

 \subsubsection{Expectation operators on infected forests}

 Recall that $\boldsymbol{\Ja_{q,n}}$ is the set of $(n+1)$ mappings $\boldsymbol{c}=\left(\boldsymbol{a},\boldsymbol{b}\right)=(c_0,\ldots,c_n)$ with $c_k=(a_k,b_k)\in I_{q}(p_k^1,p^2_k)$, for any $0\leq k\leq n$. 
 
For any symmetric function $F$ on $S_n^q$,  and any $  \boldsymbol{c}=\left(\boldsymbol{a},\boldsymbol{b}\right)$ 
and $ \boldsymbol{c^{\prime}}:=\left(\boldsymbol{a^{\prime}},\boldsymbol{b^{\prime}}\right)$ we have
$$
\boldsymbol{c}\sim \boldsymbol{c^{\prime}}
\Longrightarrow
\Delta_{\textsl{z},\boldsymbol{c}}^{(q)}(F)=\Delta_{\textsl{z},\boldsymbol{c^{\prime}}}^{(q)}(F)
$$
We check this claim using the fact that for any $a_1,a_2\in [q]^{[q]}$, and any $b\in\{0,1\}^{[q]}$, and $\sigma\in\Ga_q$ we have
$$
\AA_{a_1}\AA_{a_2}=\AA_{a_1a_2}\quad\mbox{\rm and}\quad \BB_{\textsl{z},b}=\AA_{\sigma}\BB_{\textsl{z},b\sigma}\AA_{\sigma^{-1}}
$$

 Thus, for any
${\bf f}\in\boldsymbol{\Fa_{q,n}}$ we can define unambiguously  $\Delta_{\textsl{z},\boldsymbol{f}}^{(q)}=\Delta_{\textsl{z},\boldsymbol{c}}^{(q)}$
 for any choice ${\bf c}$ of a representative of $\bf f$ in $\boldsymbol{\Ja_{q,n}}$. 

We also denote by $\boldsymbol{\Fa_{q}}(\boldsymbol{p_n})$ the set of forests with $p^1_k$-coalescences 
and $p^2_k$ infections at each level $0\leq k\leq n$. By construction, these forests are associated with
the mappings $\boldsymbol{c_n}\in \boldsymbol{I_q(p_n)}$. 
In this notation, the operators (\ref{tensormeas}) can be rewritten in terms of the expectations operators on the set of infected forests
\begin{equation}\label{conditional-expectation-forests}
\Delta_{\textsl{z},\boldsymbol{p_n}}^{(q)}=\sum_{\boldsymbol{c_n}\in \boldsymbol{I_q(p_n)}}
~
\boldsymbol{\Ca^{(q)}_{(p_n)}(c_n)}~\Delta_{\textsl{z},\boldsymbol{c_n}}^{(q)}
=
\sum_{\boldsymbol{f}\in \boldsymbol{\Fa_q(p_n)}}~~\lambda_{\boldsymbol{q,p_n}}\left(\boldsymbol{f}\right)~
\Delta_{\textsl{z},\boldsymbol{f}}^{(q)}
\end{equation}
with the probability measure $\lambda_{\boldsymbol{q,p_n}}$  given by
$$
\lambda_{\boldsymbol{q,p_n}}\left(\boldsymbol{f}\right)={\#\left(\boldsymbol{f}\right)}/{
\#\left(\boldsymbol{I_q(p_n)}\right)},
$$
where we used the shortcut notation $\#(\bf{f}):=\#[{\bf c}]$ for an arbitrary representative of $\bf f$ in $\boldsymbol{\Ja_{q,n}}$.
We summarize the above discussion with the following theorem.
\begin{theo}\label{theo-second-order-forests}
For any $m\geq 0$ we have
$$
d^{(m)}\Upsilon_{\textsl{z},n}^{(q)}
=\sum_{\boldsymbol{p_n}\in\boldsymbol{\Ta^{(m)}_{q,n}}}\boldsymbol{\tau}^{(m)}_{\boldsymbol{q,p_n}}~\left(\sum_{\boldsymbol{f}\in\boldsymbol{\Fa_q(p_n)}}\lambda_{\boldsymbol{q,p_n}}\left(\boldsymbol{f}\right)
\Delta_{\textsl{z},\bf f}^{(q)}\right).$$
\end{theo}

 \subsubsection{Infected forests}

The first order derivative is expressed in terms of two classes of infected forests.
The explicit description of the second derivative depends on $20$ different types of infected forests. We investigate them in this paragraph.

Let us fix $3<q<N$ and the time horizon $n$. There exists a single forest  $\boldsymbol{f_0}$ with trivial trees with no infection. 
There is also a single non infected forest $\boldsymbol{f^k_{1,0}}$ with only one coalescence
 at level $k$. We also have a single forest $\boldsymbol{f^k_{0,1}}$ with trivial trees and an infection at level $k$.
A synthetic description of these forests is given below.

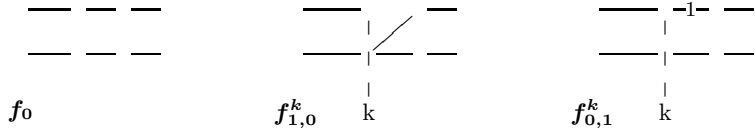
\begin{center}
\hskip.3cm
\xymatrix@C=1em@R=1em{\ar@{-}[r]  & \ar@{-}[r]   & \ar@{-}[r]  &&\\
                                                  \ar@{-}[r]&  \ar@{-}[r] & \ar@{-}[r] && \\    
               \text{\footnotesize $\boldsymbol{f_0}$}&&&&&
                                 } 
\xymatrix@C=1em@R=1em{            \ar@{-}[r]                                          &\ar@{--}[dd]         &    \ar@{-}[r]    &&&\\
                           \ar@{-}[r]                                         &  \ar@{-}[r]      \ar@{-}[ur]   &   \ar@{-}[r]    &&&\\    
           \text{\footnotesize $\boldsymbol{f^k_{1,0}}$}     &\text{\footnotesize k}&&&&
                                 }   \xymatrix@C=1em@R=1em{      \ar@{-}[r]               &          \ar@{--}[dd]  \ar@{-}[r]|{1}        &    \ar@{-}[r]   & &\\
                                                   \ar@{-}[r]               &           \ar@{-}[r]    &  \ar@{-}[r]       &    &&\\         
           \text{\footnotesize $\boldsymbol{f^k_{0,1}}$}         & \text{\footnotesize k}&&
                                 }
\end{center}
The corresponding measures are given by $\Delta_{\textsl{z},\boldsymbol{f_{0}}}^{(q)}=\eta^{\otimes q}_n$, and
the pair of measures
\begin{equation}\label{max-one-infec-coal}
\Delta_{\textsl{z},\boldsymbol{f^k_{1,0}}}^{(q)}=\eta_{n}^{\otimes (q-2)}\otimes \left[\int\eta_k(dw)~
(\delta_w\overline{Q}_{k,n})^{\otimes 2}\right]
\quad\mbox{\rm and}\quad
\Delta_{\textsl{z},\boldsymbol{f^k_{0,1}}}^{(q)}=\eta_n^{\otimes (q-1)}\otimes \delta_{\textsl{z}_k}\overline{Q}_{k,n}
\end{equation}
It is also immediate to check using (\ref{stabil}) that $$
\#\left(\boldsymbol{f_{0}}\right)=q!^{n+1}\qquad
\#\left(\boldsymbol{f^k_{1,0}}\right)=q!^{n+2}/((q-2)!2!)\quad\mbox{\rm and}\quad
\#\left(\boldsymbol{f^k_{0,1}}\right)=q!^{n+1}~q
$$

There are two non infected forests $\boldsymbol{f_{2,0}^{k,1}}$ and $\boldsymbol{f_{2,0}^{k,2}}$  
with two coalescences at level $k$. The first one has a non trivial tree with three leaves, the second one has two trees with two leaves.

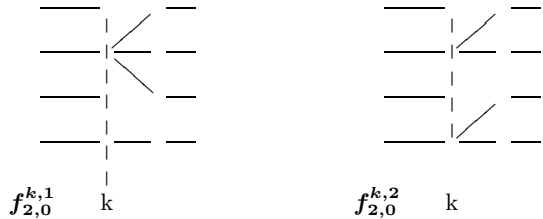
\begin{center}
\hskip.3cm
\xymatrix@C=1em@R=1em{            \ar@{-}[r]                                          &\ar@{--}[dddd]         &    \ar@{-}[r]    &&&&\\
                    \ar@{-}[r]                                         &  \ar@{-}[r]      \ar@{-}[ur] \ar@{-}[dr]  &   \ar@{-}[r]    &&&&\\     
           \ar@{-}[r]                                                &         &    \ar@{-}[r]    &&&&\\ 
         \ar@{-}[r]                                                &  \ar@{-}[r]        &    \ar@{-}[r]   &&& &     \\          
           \text{\footnotesize $   \boldsymbol{f_{2,0}^{k,1}}$}       &\text{\footnotesize k}&&&&&
                                 }\xymatrix@C=1em@R=1em{           \ar@{-}[r]                                          &\ar@{--}[ddd]         &    \ar@{-}[r]    &\\
                                 \ar@{-}[r]                                         &  \ar@{-}[r]      \ar@{-}[ur] &   \ar@{-}[r]    &\\     
       \ar@{-}[r]                                          &       &    \ar@{-}[r]    &\\
      \ar@{-}[r]                                         &  \ar@{-}[r]      \ar@{-}[ur]  &   \ar@{-}[r]    &\\                                  
           \text{\footnotesize $   \boldsymbol{f_{2,0}^{k,2}}$}       &\text{\footnotesize k}&&
                                 }

\end{center}

The corresponding measures are given by
\begin{eqnarray}
\Delta_{\textsl{z},  \boldsymbol{f_{2,0}^{k,1}} }^{(q)}
&=&\eta_{n}^{\otimes (q-3)}\otimes\left[\int\eta_k(dw)~\left(\delta_w\overline{Q}_{k,n}^{\otimes 3}\right)\right] \nonumber\\
\Delta_{\textsl{z},\boldsymbol{f_{2,0}^{k,2}} }^{(q)}(F)&
=&\eta_{n}^{\otimes (q-4)}\otimes\left\{\int \eta_k(dw_1)~\eta_k(dw_2)\left[
\left(\delta_{w_1}\overline{Q}_{k,n}\right)^{\otimes 2}\otimes\left(\delta_{w_2}\overline{Q}_{k,n}\right)^{\otimes 2}\right]\right\}\label{only-two-coal}
\end{eqnarray}
and we have $
\#\left(    \boldsymbol{f_{2,0}^{k,1}}\right)={(q!)^{n+2}}/{((q-3)!3!)}$, and
$\#\left(    \boldsymbol{f_{2,0}^{k,2}}\right)={(q!)^{n+2}}/{((q-4)!2^3)}
$.

There is also a single non coalescent forest $\boldsymbol{f_{0,2}^{k}} $ with two trivial infected trees at level $k$. 
There are two forests $\boldsymbol{f_{1,1}^{k,i}} $, $i=1,2$, with one infection and one coalescence at level $k$.
The first one  has a single coalescent tree with only one
 infected leaf. The last one has a non infected coalescent tree and a single infected trivial tree. 
 
 \begin{center}
\hskip.3cm
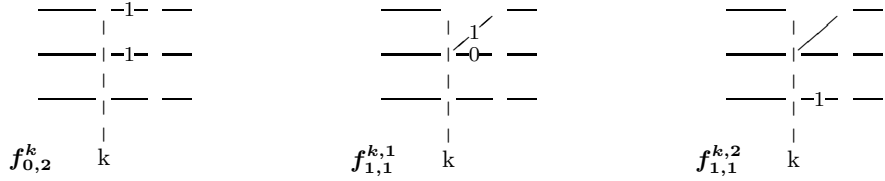
\xymatrix@C=1em@R=1em{              \ar@{-}[r]                                          &\ar@{--}[ddd]    \ar@{-}[r]|{1}     &    \ar@{-}[r]    &&&&\\
                       \ar@{-}[r]                                         &  \ar@{-}[r]|{1}      &    \ar@{-}[r]    &&&&\\     
         \ar@{-}[r]                                                &     \ar@{-}[r]      &    \ar@{-}[r]    &&&&\\      
          \text{\footnotesize $  \boldsymbol{f_{0,2}^{k}}   $}        & \text{\footnotesize k}&&&&&
                                 }
                                 \xymatrix@C=1em@R=1em{            \ar@{-}[r]                                          &\ar@{--}[ddd]        &    \ar@{-}[r]    &&&&\\
                   \ar@{-}[r]                                         &  \ar@{-}[r]|{0}   \ar@{-}[ur]|{1}   &   \ar@{-}[r]    &&&&\\     
        \ar@{-}[r]                                                &     \ar@{-}[r]      &    \ar@{-}[r]    &&&&\\      
              \text{\footnotesize $  \boldsymbol{f_{1,1}^{k,1}}   $}        & \text{\footnotesize k}&&&&&
                                 }
\xymatrix@C=1em@R=1em{          \ar@{-}[r]                                          &\ar@{--}[ddd]        &    \ar@{-}[r]    &\\
                \ar@{-}[r]                                         &  \ar@{-}[r]    \ar@{-}[ur]   &   \ar@{-}[r]    &\\     
          \ar@{-}[r]                                                &     \ar@{-}[r]|{1}      &    \ar@{-}[r]    &\\        
            \text{\footnotesize $  \boldsymbol{f_{1,1}^{k,2}} $}         & \text{\footnotesize k}&&
                                 }

\end{center}
 The corresponding measures are given by 
 \begin{eqnarray}
 \Delta_{\textsl{z}, \boldsymbol{f_{0,2}^{k,1}}  }^{(q)}&=&\eta_{n}^{\otimes (q-2)}\otimes \left(\delta_{\textsl{z}_k}\overline{Q}_{k,n}\right)^{\otimes 2}\qquad\quad
 \Delta_{\textsl{z},   \boldsymbol{f_{1,1}^{k,1}}  }^{(q)}=\Delta_{\textsl{z},n,   \boldsymbol{f_{0,1}^{k}}  }^{(q)}\nonumber\\
 \Delta_{\textsl{z},   \boldsymbol{f_{1,1}^{k,2}}   }^{(q)}&=&\eta_{n}^{\otimes (q-3)}\otimes\left[\int\eta_k(dw)~\left(\delta_{w}\overline{Q}_{k,n}\right)^{\otimes 2}\right]\otimes \left(\delta_{\textsl{z}_k}\overline{Q}_{k,n}\right)\label{one-infec-one-coal}
 \end{eqnarray}
One checks that $\#(\boldsymbol{f_{0,2}^{k}}) =q!^{n+1}~q(q-1)/2$,
$
  \#(\boldsymbol{f_{1,1}^{k,1}})=q!^{n+1}~q(q-1)$
and $ \#\left(    \boldsymbol{f_{1,1}^{k,2}}\right)=\frac{(q!)^{n+2}}{2(q-3)!}.$

We also have the traditional four non infected forests $\boldsymbol{f_{1,1}^{k,l,i}} $, $i=1,2,3,4$ with two coalescences at level $k$ and $l$~\cite{dpr-2009}.
The first one has two coalescent trees with all the leaves at level $n$. The second one also
 has two coalescent trees but one has two leaves at level $n$, the other has a leaf at level $l$ and another at level $n$.
The third one has a single coalescent tree with three leaves at level $n$, and a coalescent branch  at level $l$. The last one
has a single coalescent tree with two leaves at level $n$ and a coalescent branch  at level $l$.

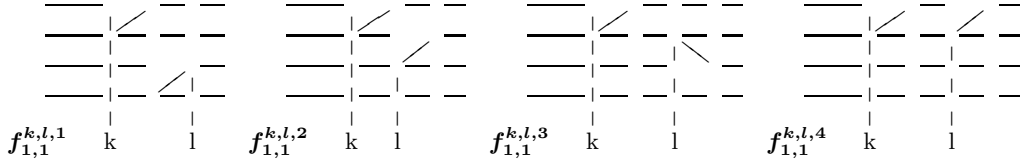
\begin{center}
\hskip.3cm
\xymatrix@C=.7em@R=.5em{\ar@{-}[r]&\ar@{--}[dddd]         &\ar@{-}[r]&\ar@{-}[r]&\\
\ar@{-}[r] &\ar@{-}[r] \ar@{-}[ur] &\ar@{-}[r]&\ar@{-}[r]&\\     
 \ar@{-}[r] &\ar@{-}[r]&&\ar@{--}[dd]  \ar@{-}[r]&\\ 
     \ar@{-}[r] &\ar@{-}[r]  & \ar@{-}[ur]   \ar@{-}[r] &\ar@{-}[r]&\\        
             \text{\footnotesize $   \boldsymbol{f_{1,1}^{k,l,1}}   $}         &\text{\footnotesize k}&&\text{\footnotesize l}&
                                 }
   \xymatrix@C=.7em@R=.5em{
 \ar@{-}[r] &\ar@{--}[dddd] &\ar@{-}[r]&\ar@{-}[r]&\\
 \ar@{-}[r]&\ar@{-}[r]\ar@{-}[ur] &&\ar@{-}[r]&\\     
 \ar@{-}[r]&\ar@{-}[r]&\ar@{--}[dd]  \ar@{-}[ur] \ar@{-}[r]  &  \ar@{-}[r]&\\ 
 \ar@{-}[r] &  \ar@{-}[r]& \ar@{-}[r] &\ar@{-}[r]&\\          
          \text{\footnotesize $\boldsymbol{f_{1,1}^{k,l,2}}$}         &\text{\footnotesize k}&\text{\footnotesize l}&&
                                 }
                              \xymatrix@C=.7em@R=.5em{            \ar@{-}[r]                                          &\ar@{--}[dddd]        &    \ar@{-}[r]    &\ar@{-}[r]&\ar@{-}[r]&\\
\ar@{-}[r]                                         &  \ar@{-}[r]      \ar@{-}[ur]   &  \ar@{-}[r]&\ar@{--}[ddd] \ar@{-}[dr] \ar@{-}[r]&\ar@{-}[r]&\\     
        \ar@{-}[r]                                                &     \ar@{-}[r]      &    \ar@{-}[r]    &   &  \ar@{-}[r]&\\ 
        \ar@{-}[r]                                                &  \ar@{-}[r]        &    \ar@{-}[r]    &   \ar@{-}[r] &\ar@{-}[r]&\\          
             \text{\footnotesize $\boldsymbol{f_{1,1}^{k,l,3}}$}     &\text{\footnotesize k}&&\text{\footnotesize l}&&
                                 }
                                 \xymatrix@C=.7em@R=.5em{           \ar@{-}[r]                                          &\ar@{--}[dddd]        &    \ar@{-}[r]    &&\ar@{-}[r]&\\
                  \ar@{-}[r]                                         &  \ar@{-}[r]      \ar@{-}[ur]   &  \ar@{-}[r]&\ar@{--}[ddd] \ar@{-}[ur] \ar@{-}[r]&\ar@{-}[r]&\\     
         \ar@{-}[r]                                                &     \ar@{-}[r]      &    \ar@{-}[r]    &  \ar@{-}[r]  &  \ar@{-}[r]&\\
     \ar@{-}[r]                                                &     \ar@{-}[r]      &    \ar@{-}[r]    &  \ar@{-}[r]  &  \ar@{-}[r]&\\      
       \text{\footnotesize $ \boldsymbol{f_{1,1}^{k,l,4}}$}        &  \text{\footnotesize k}&&\text{\footnotesize l}&&
                                 }
\end{center}
In this case, we readily check that
$$
\text{\footnotesize $\#\left(    \boldsymbol{f_{1,1}^{k,l,1}}\right)=\frac{q!^{n+2}}{4(q-4)!}\quad
\#\left(    \boldsymbol{f_{1,1}^{k,l,2}}\right)=\frac{q!^{n+2}}{(q-3)!2!}\quad
\#\left(    \boldsymbol{f_{1,1}^{k,l,3}}\right)=\frac{q!^{n+2}}{(q-3)!2!}\quad\#\left(    \boldsymbol{f_{1,1}^{k,l,4}}\right)=\frac{q!^{n+2}}{(q-2)!2!}$}
$$
and the corresponding measures are given by
\begin{eqnarray}
\Delta_{\textsl{z}, \boldsymbol{f_{1,1}^{k,l,1}}   }^{(q)}&=& 
\eta_{n}^{\otimes (q-4)}\otimes 
\left[\int~\eta_k(du)~\left(\delta_u\overline{Q}_{k,n}\right)^{\otimes 2}\right]\otimes \left[\int~\eta_l(dv)~\left(\delta_v\overline{Q}_{l,n}\right)^{\otimes 2}\right]\nonumber\\
\Delta_{\textsl{z},           \boldsymbol{f_{1,1}^{k,l,2}}    }^{(q)}&=&
\eta_{n}^{\otimes (q-3)}\otimes\left[\int \eta_k(du)~\overline{Q}_{k,l}(1)(u)~\delta_u\overline{Q}_{k,n}\right]\otimes\left[
\int~\eta_l(dv)~\left(\delta_v\overline{Q}_{l,n}\right)^{\otimes 2}\right]\nonumber\\
\Delta_{\textsl{z},             \boldsymbol{f_{1,1}^{k,l,3}}    }^{(q)}&=&
\eta_{n}^{\otimes (q-3)}\otimes\left[\int \eta_k(du)~\left(\left\{\int~\overline{Q}_{k,l}(u,dv)\left(\delta_v\overline{Q}_{l,n}\right)^{\otimes 2}\right\}\otimes \delta_u\overline{Q}_{k,n}\right)\right]\nonumber\\
\Delta_{\textsl{z},    \boldsymbol{f_{1,1}^{k,l,4}}    }^{(q)}&=&\eta_{n}^{\otimes (q-2)}\otimes\left[\int~\eta_k(du)~\overline{Q}_{k,l}(1)(u)~\overline{Q}_{k,l}(u,dv)~\left(\delta_v\overline{Q}_{l,n}\right)^{\otimes 2}\right]\label{coales-2-ref}
\end{eqnarray}

We also have two non coalescent forests  $\boldsymbol{f_{0,1,1}^{k,l,i}}$, $i=1,2$, with two infections at level $k$ and $l$.
The first one has two infected trivial trees. The second one has a trivial tree with two infections. 
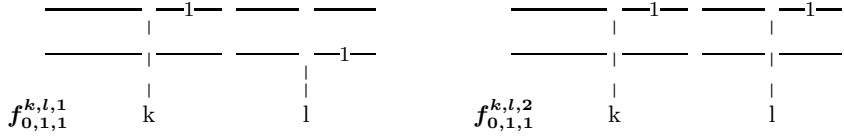
\begin{center}
\hskip.3cm
\xymatrix@C=2em@R=1em{  \ar@{-}[r]  &\ar@{--}[dd]    \ar@{-}[r]|{1}    &    \ar@{-}[r]    &   \ar@{-}[r] &\\
\ar@{-}[r]  &  \ar@{-}[r]     &  \ar@{-}[r]&\ar@{--}[d]  \ar@{-}[r]|{1}&\\            
 \text{\footnotesize $\boldsymbol{f_{0,1,1}^{k,l,1}}$}      &\text{\footnotesize k}&&\text{\footnotesize l}&&
                                 }
                                 \xymatrix@C=2em@R=1em{
 \ar@{-}[r]&\ar@{--}[dd]     \ar@{-}[r]|{1}     &    \ar@{-}[r]    &\ar@{--}[dd]   \ar@{-}[r]|{1} & \\
 \ar@{-}[r]  &  \ar@{-}[r] &\ar@{-}[r]& \ar@{-}[r]&\\              
 \text{\footnotesize $\boldsymbol{f_{0,1,1}^{k,l,2}}$} &\text{\footnotesize k}&&\text{\footnotesize l}&
                                 }
\end{center}
In this case, we have $\#\left(    \boldsymbol{f_{0,1,1}^{k,l,1}}\right)=q!^{n+1} q(q-1)$ and $
\#\left(    \boldsymbol{f_{0,1,1}^{k,l,2}}\right)=q!^{n+1} q$, and 
\begin{equation}\label{infect-2-ref}
\Delta_{\textsl{z},      \boldsymbol{f_{0,1,1}^{k,l,1}}      }^{(q)}=\eta_{n}^{\otimes (q-2)}\otimes\delta_{\textsl{z}_k}\overline{Q}_{k,n}\otimes
\delta_{\textsl{z}_l}\overline{Q}_{l,n}\quad\mbox{\rm and}\quad
\Delta_{\textsl{z},  \boldsymbol{f_{0,1,1}^{k,l,2}}    }^{(q)}=\overline{Q}_{k,l}(1)(\textsl{z}_k)~\left[\eta_{n}^{\otimes (q-1)}\otimes \delta_{\textsl{z}_l}\overline{Q}_{l,n} \right]
\end{equation}

We also have two forests $\boldsymbol{f_{1,0,1}^{k,l,i}}$, $i=1,2$, with a coalescence at level $k$ and an infection at level $l>k$.
The first one has a coalescent tree with an infection. The second one has a non infected coalescent tree and an infected trivial tree.
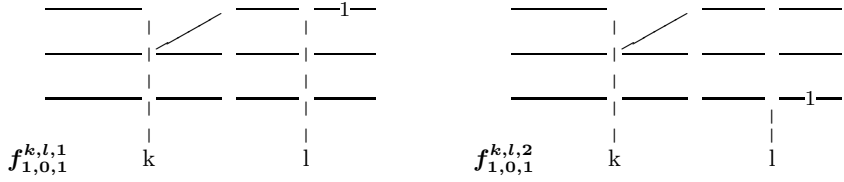
\begin{center}
\hskip.3cm
\xymatrix@C=2em@R=1em{ 
\ar@{-}[r]&\ar@{--}[ddd]        & \ar@{-}[r]&\ar@{--}[ddd]\ar@{-}[r]|{1}&&\\
\ar@{-}[r]  &  \ar@{-}[r] \ar@{-}[ur]   &\ar@{-}[r] \ar@{-}[r]&\ar@{-}[r]&&\\     
\ar@{-}[r] & \ar@{-}[r] &  \ar@{-}[r]  &\ar@{-}[r]&&\\        
 \text{\footnotesize $\boldsymbol{f_{1,0,1}^{k,l,1}}$} &\text{\footnotesize k}&&\text{\footnotesize l}&&
                                 }
 \xymatrix@C=2em@R=1em{     
 \ar@{-}[r]&\ar@{--}[ddd] &\ar@{-}[r]&\ar@{-}[r]&\\
 \ar@{-}[r]  &  \ar@{-}[r]      \ar@{-}[ur]   &  \ar@{-}[r] \ar@{-}[r]&\ar@{-}[r]&\\     
 \ar@{-}[r] &     \ar@{-}[r]      &    \ar@{-}[r]  & \ar@{--}[d]\ar@{-}[r]|{1}&\\       
          \text{\footnotesize $  \boldsymbol{f_{1,0,1}^{k,l,2}} $}         &\text{\footnotesize k}&&\text{\footnotesize l}&
                                 }
\end{center}
In this case we have $ \text{\footnotesize $\#\left(    \boldsymbol{f_{1,0,1}^{k,l,1}}\right)=q!^{n+2}/(q-2)!$}$, and
$ \text{\footnotesize $\#\left(    \boldsymbol{f_{1,0,1}^{k,l,2}}\right)=q!^{n+2}/(2(q-3)!)$}$. 
The corresponding measures are given by
\begin{eqnarray}
\Delta_{\textsl{z},         \boldsymbol{f_{1,0,1}^{k,l,1}}        }^{(q)}&=&\eta_{n}^{\otimes (q-2)}\otimes\left[\int\eta_k(du)~\overline{Q}_{k,l}(1)(u)~\delta_u\overline{Q}_{k,n}\right]\otimes \delta_{\textsl{z}_l}\overline{Q}_{l,n}\nonumber\\
\Delta_{\textsl{z},        \boldsymbol{f_{1,0,1}^{k,l,2}}        }^{(q)}&=&
\eta_{n}^{\otimes (q-3)}\otimes\left[\int\eta_k(du)~\left(\delta_u\overline{Q}_{k,n}\right)^{\otimes 2}\right]\otimes \delta_{\textsl{z}_l}\overline{Q}_{l,n}\label{coal+infect-ref}
\end{eqnarray}

Finally, there are three forests $\boldsymbol{f_{0,1,0,1}^{k,l,i}}$, $i=1,2,3$, with an infection at $k$ and a coalescence at level $l>k$.
The first one has a infected tree with a leaf at level $n$ and a non infected coalescent tree. The second one
 has a infected tree with a leaf at level $l$ and a non infected coalescent tree. And finally, the last one has an infected coalescent tree.
 
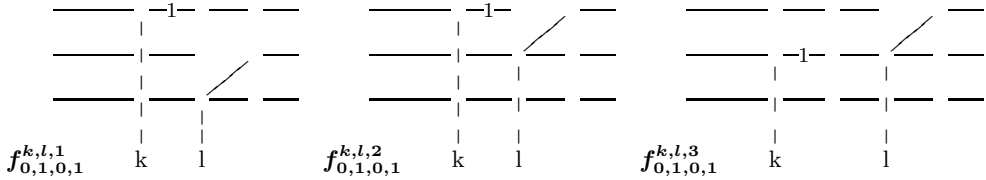
\begin{center}
\hskip.3cm
\xymatrix@C=1.2em@R=1em{ 
\ar@{-}[r]& \ar@{-}[r]|{1}   \ar@{--}[ddd]  &\ar@{-}[r]&\ar@{-}[r]&\\
\ar@{-}[r] &\ar@{-}[r]\ar@{-}[r]&&\ar@{-}[r]&\\     
 \ar@{-}[r]&\ar@{-}[r]&\ar@{--}[d]  \ar@{-}[ur] \ar@{-}[r]&\ar@{-}[r]&\\      
\text{\footnotesize $  \boldsymbol{f_{0,1,0,1}^{k,l,1}} $} &\text{\footnotesize k}&\text{\footnotesize l}&
                                 } 
                                 \xymatrix@C=1.2em@R=1em{         
 \ar@{-}[r] & \ar@{--}[ddd] \ar@{-}[r]|{1} & &\ar@{-}[r]&\\     
 \ar@{-}[r] & \ar@{-}[r] & \ar@{-}[ur] \ar@{--}[dd] \ar@{-}[r]&\ar@{-}[r]&\\ 
\ar@{-}[r] &  \ar@{-}[r]  &\ar@{-}[r]&\ar@{-}[r]&\\          
 \text{\footnotesize $\boldsymbol{f_{0,1,0,1}^{k,l,2}}  $}&\text{\footnotesize k}&\text{\footnotesize l}&&
                                 } 
                                 \xymatrix@C=1.2em@R=1em{      
                \ar@{-}[r]      &  \ar@{-}[r]      & \ar@{-}[r] \ar@{-}[r]&&\ar@{-}[r]&\\     
         \ar@{-}[r]   &   \ar@{--}[dd]    \ar@{-}[r]|{1}       & \ar@{-}[r]  &   \ar@{-}[ur]\ar@{--}[dd]     \ar@{-}[r]&\ar@{-}[r]&\\ 
        \ar@{-}[r]    &     \ar@{-}[r]      &   \ar@{-}[r] &\ar@{-}[r]&\ar@{-}[r]&\\          
          \text{\footnotesize $  \boldsymbol{f_{0,1,0,1}^{k,l,3}}   $}        &\text{\footnotesize k}&&\text{\footnotesize l}&&
                             }
\end{center}
In this case we have
$\#\left(\boldsymbol{f_{0,1,0,1}^{k,l,1}}\right)=q!^{n+2}/(2(q-3)!)$ and for any $i\in\{2,3\}$
$\#\left(    \boldsymbol{f_{1,0,1}^{k,l,i}}\right)=q!^{n+2}/(2(q-2)!)
$
In addition, the corresponding measures are given by
\begin{eqnarray}
\Delta_{\textsl{z},       \boldsymbol{f_{0,1,0,1}^{k,l,1}}       }^{(q)}&=&\eta_{n}^{\otimes (q-3)}\otimes\left[\int\eta_l(du)\left(\delta_u\overline{Q}_{l,n}\right)^{\otimes 2}\right]\otimes \delta_{\textsl{z}_k}\overline{Q}_{k,n}\nonumber\\
\Delta_{\textsl{z},           \boldsymbol{f_{0,1,0,1}^{k,l,2}}       }^{(q)}&=&\overline{Q}_{k,l}(1)(\textsl{z}_k)\left[
\eta_{n}^{\otimes (q-2)}\otimes\left\{\int\eta_l(du)\left(\delta_u\overline{Q}_{l,n}\right)^{\otimes 2}\right\}\right]\nonumber\\
\Delta_{\textsl{z},         \boldsymbol{f_{0,1,0,1}^{k,l,3}}       }^{(q)}&=&
\eta_{n}^{\otimes (q-2)}\otimes\left[\int \overline{Q}_{k,l}(\textsl{z}_k,du)\left(\delta_u\overline{Q}_{l,n}\right)^{\otimes 2}\right]\label{infect+coal-ref}
\end{eqnarray}

For any multi-index $\boldsymbol{\kappa}$, and any integer $i$ we set
$$
\overline{\Delta}^{(q)}_{\textsl{z},\boldsymbol{f_{\boldsymbol{\kappa}}^{\point,\point,i}}}
:=\sum_{0\leq k<l\leq n}\overline{\Delta}^{(q)}_{\textsl{z},\boldsymbol{f_{\boldsymbol{\kappa}}^{k,l,i}}}
\quad\mbox{\rm
with}
\quad\overline{\Delta}^{(q)}_{\textsl{z},n,\boldsymbol{f_{\boldsymbol{\kappa}}^{k,l,i}}}:=\Delta^{(q)}_{\textsl{z},n,\boldsymbol{f_{\boldsymbol{\kappa}}^{k,l,i}}}-\eta^{\otimes q}_n
$$

\subsubsection{First and second derivatives}

To describe with some precision the first two order derivatives of  the mapping $N\mapsto \Upsilon_{\textsl{z},n}^{(q)}$ we need 
 to compute the expectation operators on random infected forests defined in
 (\ref{conditional-expectation-forests}). 
 The ones associated with forests with at most one infection or one coalescence at some level
 only depend one one class of forests. Thus, using (\ref{max-one-infec-coal}) their description is immediate. Using (\ref{only-two-coal}),
 the centered operator associated with non infected forests with a couple of coalescence at some level is given by 
 $$
{{\overline{\Delta}}}^{(q)}_{\textsl{z},\boldsymbol{f_{2,0}^{\point,\star}}}:=
\frac{1}{1+\frac{3}{4}(q-3)}~{{\overline{\Delta}}}^{(q)}_{\textsl{z},\boldsymbol{f_{2,0}^{\point,1}}}+\left(1-\frac{1}{1+\frac{3}{4}(q-3)}\right)
~
{{\overline{\Delta}}}^{(q)}_{\textsl{z},\boldsymbol{f_{2,0}^{\point,2}}}
 $$
 In much the same way, by (\ref{one-infec-one-coal}) the one associated with forests with a single coalescence and a single infection at 
 some level is given by
 $$
{{\overline{\Delta}}}^{(q)}_{\textsl{z},\boldsymbol{f_{1,1}^{\point,\star}}}:=
\frac{2}{q}~{{\overline{\Delta}}}^{(q)}_{\textsl{z},\boldsymbol{f_{1,1}^{\point,1}}}+\left(1-\frac{2}{q}~\right)~{{\overline{\Delta}}}^{(q)}_{\textsl{z},\boldsymbol{f_{1,1}^{\point,2}}}
 $$
 In view of (\ref{coales-2-ref}), the centered expectation operator associated with forests with a single coalescence at two different levels is given by
  
\begin{eqnarray*}
{{\overline{\Delta}}}^{(q)}_{\textsl{z},\boldsymbol{f_{1,1}^{\point,\point,\star}}}&
:=&\frac{(q-2)(q-3)}{(q-2)(q-3)+4(q-2)+2}~
{{\overline{\Delta}}}^{(q)}_{\textsl{z},\boldsymbol{f_{1,1}^{\point,\point,1}}}\\
&&\hskip3cm+\frac{2(q-2)}{(q-2)(q-3)+4(q-2)+2}
\left[{{\overline{\Delta}}}^{(q)}_{\textsl{z},\boldsymbol{f_{1,1}^{\point,\point,2}}}+{{\overline{\Delta}}}^{(q)}_{\textsl{z},\boldsymbol{f_{1,1}^{\point,\point,3}}}\right]\\
&&
\hskip5cm
+\frac{2}{(q-2)(q-3)+4(q-2)+2}~
{{\overline{\Delta}}}^{(q)}_{\textsl{z},\boldsymbol{f_{1,1}^{\point,\point,4}}}
\end{eqnarray*}

Using (\ref{infect-2-ref}) the one associated with non coalescent forests with a single infection at two different levels is given by
$$
{{\overline{\Delta}}}^{(q)}_{\textsl{z},\boldsymbol{f_{0,1,1}^{\point,\point,\star}}}:=
\left(1-\frac{1}{q}\right)~{{\overline{\Delta}}}^{(q)}_{\textsl{z},\boldsymbol{f_{0,1,1}^{\point,\point,1}}}+
\frac{1}{q}~{{\overline{\Delta}}}^{(q)}_{\textsl{z},\boldsymbol{f_{0,1,1}^{\point,\point,2}}}
$$ 
 Finally, using (\ref{coal+infect-ref}
) and (\ref{infect+coal-ref}) the operator associated with a single coalescence and a single infection at two different levels are given by
\begin{eqnarray*}
{{\overline{\Delta}}}^{(q)}_{\textsl{z},\boldsymbol{f_{1,0,1}^{\point,\point,\star}}}&:=&\frac{2}{q}~ 
{{\overline{\Delta}}}^{(q)}_{\textsl{z},\boldsymbol{f_{1,0,1}^{\point,\point,1}}}+\left(1-\frac{2}{q}\right)
{{\overline{\Delta}}}^{(q)}_{\textsl{z},\boldsymbol{f_{1,0,1}^{\point,\point,2}}}
\end{eqnarray*}
and
\begin{eqnarray*}
{{\overline{\Delta}}}^{(q)}_{\textsl{z},\boldsymbol{f_{0,1,0,1}^{\point,\point,\star}}}&:=&
\left(1-\frac{2}{q}\right)~
{{\overline{\Delta}}}^{(q)}_{\textsl{z},\boldsymbol{f_{0,1,0,1}^{\point,\point,1}}}+\frac{1}{q}~
{{\overline{\Delta}}}^{(q)}_{\textsl{z},\boldsymbol{f_{0,1,0,1}^{\point,\point,2}}}
+\frac{1}{q}~
{{\overline{\Delta}}}^{(q)}_{\textsl{z},\boldsymbol{f_{0,1,0,1}^{\point,\point,3}}}
\end{eqnarray*}

Expanding the formulae stated in theorem~\ref{theo-second-order-forests}, extending the combinatorial
methods developed in~\cite{dpr-2009} for computing the cardinals $\#\left(\boldsymbol{f}\right)$ we
prove the following expansions.

\begin{cor}\label{cor-Q-expansions}
The first three derivatives of $\Upsilon_{\textsl{z},n}^{(N,q)}$ are given by
$$
\begin{array}{l}
d^{(0)}\Upsilon_{\textsl{z},n}^{(q)}=\eta_n^{\otimes q}\\
\\
d^{(1)}\Upsilon_{\textsl{z},n}^{(q)}=  \tau^{(1)}_{q,1,0}
{{\overline{\Delta}}}_{\textsl{z},\boldsymbol{f^{\point}_{1,0}}}^{(q)}+ \tau^{(1)}_{q,0,1}~{{\overline{\Delta}}}_{\textsl{z},\boldsymbol{f^{\point}_{0,1}}}^{(q)}\\
\\
d^{(2)}\Upsilon_{\textsl{z},n}^{(q)}\\
\\
=
\tau^{(2)}_{q,1,0}~{{\overline{\Delta}}}^{(q)}_{\textsl{z},\boldsymbol{f_{1,0}^{\point}}}+\tau^{(2)}_{q,0,1}~~{{\overline{\Delta}}}^{(q)}_{\textsl{z},\boldsymbol{f_{0,1}^{\point}}}+\tau^{(2)}_{q,1,1}~{{\overline{\Delta}}}^{(q)}_{\textsl{z},\boldsymbol{f_{1,1}^{\point,\star}}}
+\tau^{(2)}_{q,2,0}~{{\overline{\Delta}}}^{(q)}_{\textsl{z},\boldsymbol{f_{2,0}^{\point,\star}}}
+\tau^{(2)}_{q,0,2}~
{{\overline{\Delta}}}^{(q)}_{\textsl{z},\boldsymbol{f_{0,2}^{\point}}}\\
\\
+\left(\tau^{(1)}_{q,1,0}\right)^2~{{\overline{\Delta}}}^{(q)}_{\textsl{z},\boldsymbol{f_{1,1}^{\point,\point,\star}}}
+\left(\tau^{(1)}_{q,0,1}\right)^2 ~{{\overline{\Delta}}}^{(q)}_{\textsl{z},\boldsymbol{f_{0,1,1}^{\point,\point,\star}}}
+\tau^{(1)}_{q,1,0}\tau^{(1)}_{q,0,1}\left\{
{{\overline{\Delta}}}^{(q)}_{\textsl{z},\boldsymbol{f_{1,0,1}^{\point,\point,\star}}}
+
{{\overline{\Delta}}}^{(q)}_{\textsl{z},\boldsymbol{f_{0,1,0,1}^{\point,\point,\star}}}
\right\}\\
\\
\hskip6cm+n~\tau^{(1)}_{q,0,0}\left[\tau^{(1)}_{q,1,0}~{{\overline{\Delta}}}^{(q)}_{\textsl{z},\boldsymbol{f_{1,0}^{\point}}}+\tau^{(1)}_{q,0,1}~{{\overline{\Delta}}}^{(q)}_{\textsl{z},\boldsymbol{f_{0,1}^{\point}}}\right]
\\
\end{array}
$$
with the parameters $\tau^{(m)}_{q,p_1,p_2}$ given in (\ref{ref-beta-m}).
\end{cor}

When $q=1$, all the terms are null except $ \tau^{(1)}_{1,0,1}=1=-\tau^{(1)}_{1,0,0}$. In this case, we find that
\begin{eqnarray*}
d^{(1)}\Upsilon_{\textsl{z},n}^{(1)}
&=&  \sum_{0\leq k\leq n}\left[\Delta_{\textsl{z},\boldsymbol{f^{k}_{0,1}}}^{(1)}-\eta_n\right]
=\sum_{0\leq k\leq n}\delta_{\textsl{z}_k}(\overline{Q}_{k,n}-\eta_n)\\
d^{(2)}\Upsilon_{\textsl{z},n}^{(1)}&
=&~\overline{\Delta}^{(1)}_{\textsl{z},\boldsymbol{f_{0,1,1}^{\point,\point,2}}}-n~\overline{\Delta}^{(1)}_{\textsl{z},\boldsymbol{f_{0,1}^{\point}}}\\
&=&\sum_{0\leq k<l\leq n}~\left[\overline{Q}_{k,l}(1)(\textsl{z}_k)~\delta_{\textsl{z}_l}\overline{Q}_{l,n}-\eta_n\right]-n\sum_{0\leq k\leq n}\left[\delta_{\textsl{z}_k}\overline{Q}_{k,n}-\eta_n\right]
\end{eqnarray*}

\section{Some extensions and open questions}\label{last-section}
\subsection{Island type methodologies}\label{islands-sec}

Particle MCMC methods are computationally intensive sampling techniques. As discussed in~\cite{durham,vddm-2014}, parallel and distributed computations
provide an appealing solution to tackle these issues.
The central idea of Island models is run in parallel $N_2$ particle models with $N_1$ individuals, instead of running a single particle model with $N_1N_2$ particles.
These $N_2$ batches are termed islands in reference to dynamic population models. Within each island the $N_1$ individuals evolve as a standard genetic 
type particle model. In this interpretation, island particle models can be thought as a parallel implementation of particle models. In the further development of this section, we show that these methodologies can also be used in a natural way to design island type particle MCMC samplers. 

 To design these models, we consider a collection of bounded and non-negative potential functions
$\mathfrak{G_n}$ on some  measurable state spaces $\mathfrak{E_{n}}$, with $n\in\NN$. We let $\mathfrak{X_n}$ be a Markov chain on $\mathfrak{E_{n}}$ with initial distribution $\mu_0\in\Pa(\mathfrak{E_{0}})$ and
some Markov transitions $\mathfrak{M_n}$ from $\mathfrak{E_{n-1}}$ into $\mathfrak{E_{n}}$. The Feynman-Kac measures $(\mu_n,\nu_n)$  associated with the parameters
$(\mathfrak{G_n},\mathfrak{M_n})$ are defined for any $\mathfrak{f_n}\in\Ba(\mathfrak{E_n})$ by the formulae
\begin{equation}\label{FK-def-intro-ref-frack}
\mu_n(\mathfrak{f_n})
:={\nu_n(\mathfrak{f_n})}/{\nu_n(1)}\quad\mbox{\rm with}\quad\nu_n(\mathfrak{f_n}):=
\EE\left(\mathfrak{f_n}(\mathfrak{X_n})~\prod_{0\leq p<n}\mathfrak{G_p}(\mathfrak{X_p})\right)
\end{equation} 
The mean field $N^{\prime}$-particle approximation 
 $$
X^{\prime}_n=\left(X^{\prime i}_n\right)_{1\leq i\leq N^{\prime}}\in S^{\prime}_n:=\mathfrak{E_{n}^{\rm [N^{\prime}]}}$$ 
of these Feynman-Kac
models is defined as in (\ref{mean-field-intro}) by considering the evolution semigroup of the Feynman-Kac model $\mu_n$.

We let $M_n^{\prime}$ be Markov transitions of $X^{\prime}_n$ and we 
consider the potential functions $G^{\prime}_n$ on $S^{\prime}_n$ defined by
\begin{equation}\label{potential-islands}
G^{\prime}_n(X^{\prime}_n)=m(X^{\prime}_n)(\mathfrak{G_n})=\frac{1}{N^{\prime}}\sum_{1\leq i\leq N^{\prime}}\mathfrak{G_n}\left(X^{\prime i}_n\right)
\end{equation}
We let $(\eta_n^{\prime},\gamma^{\prime}_n)$ be the Feynman-Kac measures   associated with the parameters
$(G_n^{\prime},M^{\prime}_n)$. In this framework, the unbiasedness properties of the unnormalized Feynman-Kac particle measures takes the form
\begin{equation}\label{transfert-islands-prime}
\begin{array}{l}
f_n^{\prime}(X_n^{\prime})=m(X^{\prime}_n)(\mathfrak{f_n})\\
\\
~\Longrightarrow~
\nu_n(\mathfrak{f_n})=\EE\left(\mathfrak{f_n}(\mathfrak{X_n})~\prod_{0\leq p<n}\mathfrak{G_p}(\mathfrak{X_p})\right)=
\EE\left(f_n^{\prime}(X_n^{\prime})~\prod_{0\leq p<n}G^{\prime}_p(X^{\prime}_p)\right)=\gamma_{n}^{\prime}(f^{\prime}_n)
\end{array}
\end{equation}
The path space version $(\eta_n,\gamma_n)$ 
of these measures are defined by the Feynman-Kac measures associated with the historical process $X_n$ and the potential function $G_n$
given by
$$
X_n=\left(X^{\prime}_0,\ldots,X^{\prime}_n\right)\in S_n=(S^{\prime}_0\times\ldots\times S^{\prime}_n)
\quad
\mbox{\rm and}\quad
G_n(X_n)=G^{\prime}_n(X^{\prime}_n)
$$
The mean field $N$-particle approximations $\xi_n^{\prime}=\left(\xi^{\prime i}_n\right)_{1\leq i\leq N}$ 
of the measures $(\eta_n^{\prime},\gamma^{\prime}_n)$ can be interpreted as a genetic type model
 island type particles 
$$
\forall 1\leq i\leq N\qquad
\xi^{\prime i}_n=\left(\xi^{\prime i,j}_n\right)_{1\leq j\leq N^{\prime}}\in S^{\prime}_n:=\mathfrak{E_{n}^{\rm [N^{\prime}]}}
$$
with mutation transitions $M^{\prime}_n$ and the selection potential functions $G^{\prime}_n$ given in (\ref{potential-islands}).
By construction, the $N$-particle approximation $\xi_n$ of the path space measures $(\eta_n,\gamma_n)$ is an
a genealogical tree based model in the space of islands. Each particle 
$$
\xi_n^i=\left(\xi_{0,n}^{\prime i},\ldots,\xi_{n,n}^{\prime i}\right)\in S_n=\left(
\mathfrak{E_{0}^{\rm [N^{\prime}]}}\times\ldots\times
\mathfrak{E_{n}^{\rm [N^{\prime}]}}\right)
$$
represents the line of island ancestor $\xi^{\prime i}_{p,n}\in \mathfrak{E_{p}^{\rm [N^{\prime}]}} $ of the $i$-th island $\xi^{\prime i}_{n,n}=\xi^{\prime i}_n\in \mathfrak{E_{n}^{\rm [N^{\prime}]}}$ 
at time $n$, at every level $0\leq p\leq n$, with $1\leq i\leq N$. In other words, 
$(\eta_n,\gamma_n,\xi_n)$ is the historical version of the Feynman-Kac model
$(\gamma_n^{\prime},\eta_n^{\prime},\xi_n^{\prime})$. 
In this case,the dual mean field particle model $\Xa_{n}$
 evolves on the state spaces  $\Sa_n=S_n^{[N]}$, with a frozen trajectory of islands $X_n$.

This model can be interpreted as the evolution of $N$ interacting islands
 $$
\forall 1\leq i\leq N
\qquad \Xa_{n}^i=\left(\Xa_{n}^{i,j}\right)_{1\leq j\leq N^{\prime}}\in  
\mathfrak{E_{n}^{\rm [N^{\prime}]}}$$ with $N^{\prime}$ individuals in each island. The conditional particle Markov chain
models discussed in section~\ref{g+b-ref} can be used without further work to design island type particle Markov chain
models with the target measure $\eta_n$. Using (\ref{transfert-islands-prime}), we see that the $S^{\prime}_n$-marginal of $\eta_n$
can be used to compute any Feynman measures of the form (\ref{FK-def-intro-ref-frack}). Similar constructions can be developed to design
a backward-sampling based particle MCMC model.

 Of course, we can iterate these Russian nesting doll type constructions at any level. For a more thorough discussion on these island type particle methodologies we refer the reader to~\cite{dr-2010,d-2013}, and the recent article~\cite{vddm-2014}. 
An important question is to analyze the convergence properties of the islands type particle models presented above 
in terms of the number of individual and the number of islands.

\end{document}